\documentclass[11pt]{amsart}

\usepackage{amsmath}
\usepackage{amsfonts}
\usepackage{amssymb}
\usepackage{graphicx}
\usepackage{mathrsfs}
\usepackage{pb-diagram}
\usepackage{epstopdf}
\usepackage{amscd}
\usepackage{color}
\usepackage{verbatim}
\usepackage[all]{xy}
\newtheorem{theorem}{Theorem}[section]
\newtheorem{lemma}[theorem]{Lemma}
\newtheorem{corollary}[theorem]{Corollary}
\newtheorem{setup}[theorem]{Setup}
\newtheorem{prop}[theorem]{Proposition}
\newtheorem{defn}[theorem]{Definition}

\newtheorem{remark}[theorem]{Remark}
\newtheorem{conjecture}[theorem]{Conjecture}
\newtheorem{assumption}[theorem]{Assumption}

\numberwithin{equation}{section}

\newcommand{\Z}{\mathbb{Z}}
\newcommand{\R}{\mathbb{R}}
\newcommand{\C}{\mathbb{C}}
\newcommand{\bS}{\mathbb{S}}

\newcommand{\bx}{\boldsymbol{x}}
\newcommand{\WH}[1]{\widehat{#1}}
\newcommand{\CM}{\mathcal{M}}

\newcommand{\CF}{\mathcal{F}}

\newcommand{\CC}{\mathcal{C}}
\newcommand{\CL}{\mathcal{L}}
\newcommand{\bL}{\mathbf{L}}
\newcommand{\BL}{\mathbb{L}}
\newcommand{\CP}{\mathbb{C}P}

\newcommand{\Hom}{\mathrm{Hom}}

\newcommand{\Jac}{\mathrm{Jac}}

\newcommand{\WT}[1]{\widetilde{#1}}
\newcommand{\OL}[1]{\overline{#1}}
\newcommand{\oi}{\overline{\iota}}

\newcommand{\integer}{\mathbb{Z}}

\newcommand{\real}{\mathbb{R}}
\newcommand{\cpx}{\mathbb{C}}

\newcommand{\consti}{\mathbf{i}\,}
\newcommand{\conste}{\mathbf{e}}
\newcommand{\ev}{\mathrm{ev}}

\newcommand{\proj}{\mathbb{P}}
\newcommand{\bP}{\mathbb{P}}

\newcommand{\moduli}{\mathcal{M}}
\newcommand{\cM}{\mathcal{M}}

\newcommand{\AI}{A_\infty}
\newcommand{\Tw}{\mathrm{Tw}}
\newcommand{\LM}{\mathcal{LM}}
\newcommand{\functor}{\LM^\BL}
\newcommand{\Fuk}{\mathcal{F}uk}
\newcommand{\MF}{\mathcal{MF}}
\newcommand{\HF}{\mathrm{HF}}
\newcommand{\Id}{\mathrm{Id}}
\newcommand{\Mor}{\mathrm{Mor}}

\newcommand{\arc}[1]{\stackrel{\frown}{#1}}

\begin{document}

\author[Cho]{Cheol-Hyun Cho}
\address{Department of Mathematical Sciences, Research institute of Mathematics\\ Seoul National University\\ San 56-1, 
Shinrimdong\\ Gwanakgu \\Seoul 47907\\ Korea}
\email{chocheol@snu.ac.kr}

\author[Hong]{Hansol Hong}
\address{Department of Mathematics/The Institute of Mathematical Sciences	 \\ The Chinese University of Hong Kong\\Shatin\\New Territories \\ Hong Kong}
\email{hhong@ims.cuhk.edu.hk, hansol84@gmail.com}

\author[Lau]{Siu-Cheong Lau}
\address{Department of Mathematics\\ Harvard University\\ One Oxford Street\\ Cambridge \\ MA 02138\\ USA}
\email{s.lau@math.harvard.edu}

\title[Localized mirror functor and HMS for $\mathbb{P}^1_{a,b,c}$]{Localized mirror functor for Lagrangian immersions, and homological mirror symmetry for $\mathbb{P}^1_{a,b,c}$}

\begin{abstract}
This paper gives a new way of constructing Landau-Ginzburg mirrors using deformation theory of Lagrangian immersions motivated by the works of Seidel, Strominger-Yau-Zaslow and Fukaya-Oh-Ohta-Ono.  Moreover we construct a canonical functor from the Fukaya category to the mirror category of matrix factorizations.  This functor derives homological mirror symmetry under some explicit assumptions.

As an application, the construction is applied to spheres with three orbifold points to produce their quantum-corrected mirrors and derive homological mirror symmetry.  Furthermore we discover an enumerative meaning of the (inverse) mirror map for elliptic curve quotients.
\end{abstract}
%\subjclass{57R18}
%\keywords{orbifold, groupoid, equivariant immersion}
\maketitle

\tableofcontents

\section{Introduction}
Homological mirror symmetry \cite{kontsevich94} and the SYZ program \cite{SYZ96} have led to deep developments in symplectic and algebraic geometry.  Both of them stem from the same idea that geometry of Lagrangian branes should correspond to geometry of coherent sheaves on the mirror.  Homological mirror symmetry focuses on the interactions among all geometric objects and the equivalence between two sides of the mirror, while SYZ emphasizes on the geometric origin of the mirror transform which induces the equivalence.

Recently homological mirror symmetry for the quintic Calabi-Yau threefold has been proved \cite{Sheridan11,Nohara-Ueda,Sh}.  On the other hand, wall-crossing and scattering which occur up to infinite order in the Gross-Siebert program \cite{GS07} are extremely complicated for the quintic.  Writing down the mirror equation from the current SYZ program is already highly non-trivial, and homological mirror symmetry is hidden behind the order-by-order quantum corrections.  Apparently there is a gap between homological mirror symmetry and the SYZ program.

This motivates the following question:\\
\textbf{Question}: \emph{Can we modify and generalize the SYZ approach such that homological mirror symmetry naturally follows from construction?}

One purpose of this paper is to give a positive answer to this question.

An important geometric ingredient in the proof of homological mirror symmetry for genus two curves by Seidel \cite{Se} and Fermat-type Calabi-Yau hypersurfaces by Sheridan \cite{Sh} is the use of specific Lagrangian immersions.  Their works motivates our current program: instead of taking Lagrangian torus fibrations in the original SYZ approach, we take a finite set of Lagrangians with mild singularities (namely, immersions) in order to construct the mirror.  In this paper we will focus on the construction when we take exactly one Lagrangian $\BL$.  We will formulate the construction with more than one Lagrangians in a forthcoming paper.

More precisely, we use one Lagrangian $\BL$ which is possibly \emph{immersed} and its deformation theory to construct a Landau-Ginzburg model $W$ as a Lagrangian Floer potential.  The flexibility of being immersed allows the Lagrangian deformation theory to capture more information -- in good cases the Lagrangian immersion split-generates the whole derived Fukaya category.  Our constructive approach can be regarded as a generalized formulation of SYZ in which general Lagrangian branes are used instead of tori.  We call this construction to be \emph{generalized SYZ}.

\begin{remark} \label{rmk:genSYZ}
The terminology `generalized SYZ\rq{} was first used by Aganagic-Vafa \cite{AV12}.  They cooked up non-compact Lagrangian branes in the resolved conifold from knots inside $\bS^3$.  Then they used these Lagrangians, instead of compact fibers of a Lagrangian torus fibration in the SYZ program, to construct the mirrors which capture information about the knot invariants.  In principle their work is coherent with the main theme of our paper, namely we use a general Lagrangian brane instead of a fibration to construct the mirror geometry.  The differences are the following.  First of all we consider a general K\"ahler manifold rather than restricting to the case of the resolved conifold.  Second we use compact Lagrangian immersions rather than non-compact Lagrangian branes to carry out the mirror construction.  Third we use the immersed Lagrangian Floer theory (in particular weakly unobstructed deformations) in our construction rather than a physical theory (and in particular we count discs of Maslov index two rather than zero).  While the detailed techniques and situations are different, both works aim to use general Lagrangian branes rather than tori to generalize the SYZ construction.
\end{remark}

One advantage of our construction is that it avoids complicated scattering and gluing, so the Landau-Ginzburg model $W$, which is roughly speaking the generating function of countings of $J$-holomorphic polygons bounded by the Lagrangian immersion $\BL$, comes out in a natural way.

Our generalized SYZ construction has a much more direct relationship with homological mirror symmetry: we construct a naturally associated $\AI$-functor $\mathcal{LM}^\BL$ from the Fukaya category of $X$ to the category of matrix factorization of $W$:

\begin{theorem}[see Theorem \ref{thm:locmirfun}] \label{thm:locmirfun1}
There exists an $\AI$-functor
$$\mathcal{LM}^\BL:  \mathcal{F}uk\,_\lambda (X) \to \mathcal{MF}(W - \lambda).$$
Here $\mathcal{F}uk\,_\lambda (X)$ is the Fukaya category of $X$ (as an $\AI$-category) whose objects are weakly unobstructed Lagrangians with potential value $\lambda$,
and $\mathcal{MF}(W - \lambda)$ is the dg category of matrix factorizations of $W - \lambda$.
\end{theorem}
%This explains why we call $W$ a localized LG mirror model. So far, there are very few cases where such homological mirror functor is known (to cite who??)

%For instance, consider the elliptic curve quotient $X=E/ (\Z/3) = \bP^1_{3,3,3}$.  We take three Lagrangians in $\mathcal{F}uk\,_0 (X)$, called $\bar{L}, \bar{L}_l, \bar{L}_s$ and transform them using our mirror functor $\mathcal{LM}^\BL$.  We explicitly obtain their corresponding $(4\times 4), (3\times 3), (2\times 2)$ matrix factorizations (see \eqref{MFfinal}, \eqref{longdiagE}, \eqref{MFshort}).

Moreover, our functor can be used to derive homological mirror symmetry under reasonable assumptions.  \emph{This provides a uniform and functorial understanding} to homological mirror symmetry.

\begin{theorem}[see Theorem \ref{thm:criterion_equiv}] \label{thm:intro:equiv}
Suppose that there exists a set of Lagrangians $\{L_i: i \in I\}$ which split-generates 
$D^\pi \Fuk_0 (X)$, and suppose the functor $\functor$ induces an isomorphism on cohomologies
$\HF (L_i,L_j) \overset{\cong}{\rightarrow} \Mor (\functor(L_i),\functor(L_j))$
for all $i,j \in I$.
Then the derived functor
$D^{\pi} (\functor): D^{\pi} \Fuk_0 (X) \to D^{\pi} \MF (W) $
is fully faithful.  Furthermore if $\{\mathcal{LM}^\BL(L_i): i \in I\}$ split-generates $D^\pi \MF (W)$, then $D^{\pi} (\functor)$ is a quasi-equivalence.
\end{theorem}

\begin{remark}
In this paper we work with $\Z/2$ grading: the Lagrangians are $\Z/2$-graded, and the matrix factorizations on the mirror side are also $\Z/2$-graded.  A $\Z$-graded version of our mirror functor can also be natural formulated (when $X$ is Calabi-Yau and $\BL$ is $\Z$-graded), which will appear in our forthcoming paper.
\end{remark}

%In this paper, we first propose another procedure to find localized Landau-Ginzburg mirror, which we call a generalized SYZ construction.  For this, consider a symplectic manifold $X$ and an immersed Lagrangian submanifold $\BL$.  We consider (formal) deformations of immersed Lagrangian submanifolds, and consider  $J$-holomorphic disc potential function $W$ to define a Landau-Ginzburg model. For this, $\AI$-algebra of the immersed Lagrangian $\BL$ has to be weakly unobstructed.  Note that this procedure does {\em not} involve an explicit Lagrangian torus fibration for $T$-duality.  This is called ``localized'', since the constructed disc potential depends on the choice of an immersed Lagrangian submanifold $\BL$ in $X$.

%Next, for such immersed Lagrangian $\BL$, and weak bounding cochains $b$ (we set $\mathbb{L} = (L,b)$), we construct our localized mirror functor $\mathcal{LM}^\BL$.

\emph{Our generalized SYZ construction and mirror functor work in any dimension}.  For instance, the construction can be applied to derive mirror symmetry for \emph{rigid Calabi-Yau manifolds} in any dimensions, \emph{with quantum corrections}, which will be discussed in a future joint work with Amorin.  Even going back to the original setting of SYZ using Lagrangian tori, our mirror functor leads to interesting construction of matrix factorizations in the toric case \cite{CHL-toric}.  Furthermore, in Section \ref{ndim}, we give a conjectural description of our program for Fermat-type Calabi-Yau hypersurfaces.  We will apply our construction to other classes of examples such as toric Calabi-Yau orbifolds, rigid Calabi-Yau manifolds and higher-genus orbifold surfaces in a series of forthcoming works.

In the later part of this paper, we apply our generalized SYZ program to construct the mirror of the orbifold projective line $X=\mathbb{P}^1_{a,b,c}$.  While $X=\mathbb{P}^1_{a,b,c}$ is just one-dimensional, it has very rich geometry due to the presence of the orbifold points and mirror symmetry in this case is very interesting (see Remark \ref{rmk:Ta} and \ref{rmk:Rossi}).  Note that \emph{the original SYZ construction is not available for $\mathbb{P}^1_{a,b,c}$}.  Thus our generalized approach produces \emph{new} mirror pairs which are not reacheable by the original SYZ approach.  

The Landau-Ginzburg mirror of $\mathbb{P}^1_{a,b,c}$ in the existing literature is the polynomial
$$ x^a + y^b + z^c + \sigma xyz. $$
From our point of view it is indeed just an approximation, and quantum corrections are necessary which makes $W$ much more non-trivial.  Namely, the above expression only contains leading order terms, and the actual mirror has higher order terms.  These higher order terms are important to make mirror symmetry works, especially in the hyperbolic case $1/a + 1/b + 1/c < 1$.

We prove homological mirror symmetry for the orbifold projective line $X=\mathbb{P}^1_{a,b,c}$ using our mirror functor.  Namely, 

\begin{theorem} \label{thm:HMS_Pabc}
Let $X=\mathbb{P}^1_{a,b,c}$ and $W$ be its generalized SYZ mirror. Assume $\frac{1}{a} + \frac{1}{b} + \frac{1}{c} \leq 1$.
 The $A_\infty$-functor $\mathcal{LM}^\BL$ in Theorem \ref{thm:locmirfun1} derives homological mirror symmetry, that is,
the split-closed derived functor of
$\mathcal{LM}^\BL:  \mathcal{F}uk (\mathbb{P}^1_{a,b,c}) \to \mathcal{MF}(W)$
is an equivalence of triangulated categories 
$$ D^{\pi} ({\Fuk} (\mathbb{P}^1_{a,b,c}) ) \cong  D^{\pi} ( \mathcal{MF}(W)),$$
for $\frac{1}{a} + \frac{1}{b} + \frac{1}{c} \leq 1 $.
\end{theorem}

For its proof, we verify the conditions in Theorem \ref{thm:intro:equiv}.  Namely, we compute the matrix factorization for $\BL$ under the $\AI$-functor $\mathcal{LM}^\BL$, and show that both $\BL$ and the corresponding matrix factorization split-generate respective derived categories.  We also show that the functor induces an isomorphism on the endomorphisms of $\BL$.

Another purpose of this paper is to explain the enumerative meaning of the mirror map\footnote{The mirror map we mention here is a map from the K\"ahler moduli to the complex moduli of the mirror, which is indeed the ``inverse mirror map'' to match the conventions of existing literatures.  For simplicity we call it the mirror map in the rest of this introduction.}.  It has been expected that coefficients of the mirror map have enumerative meanings in terms of open Gromov-Witten invariants.  A tropical version of such a conjecture was made by Gross-Siebert \cite{GS07} using scattering of tropical discs.  However, it was still an open problem to make a precise statement for compact Calabi-Yau manifolds using open Gromov-Witten invariants. 

Now our generalized SYZ mirror construction produces a map from the K\"ahler moduli of $X$ to the complex moduli of its mirror $W$, which we call the generalized SYZ map\footnote{Such a map can also be produced by Seidel's approach which was studied by Zaslow, namely comparing morphisms between generators of the categories using homological mirror symmetry, see Remark \ref{rmk:Zaslow}.  Thus it can also be called the Seidel's mirror map.  Since our work was motivated by generalizing the SYZ program of constructing the mirror family and the mirror functor, we call it to be the generalized SYZ map.} (see Definition \ref{def:SYZmap}).  We prove that the mirror map equals to the generalized SYZ map for the elliptic curve quotient $E/(\Z/3)$.  Since the generalized SYZ map is written in terms of polygon counting, such an equality establishes a precise enumerative meaning of the mirror map.  Moreover countings in the case of elliptic curves are all integers, and hence it also explains integrality of the mirror map.  We derive this equality for the elliptic curve quotients $E/(\Z/4)$ and $E/(\Z/6)$ in our joint work with Kim \cite{CHKL14}.

\begin{theorem}[Theorem \ref{mir=SYZ}] \label{mir=SYZ1}
The mirror map equals to the generalized SYZ map for the elliptic curve quotient $X = E/(\Z/3) = \bP^1_{3,3,3}$, where $E$ is the elliptic curve with complex multiplication by cube root of unity.
\end{theorem}  

\begin{remark} \label{rmk:Zaslow}
Homological mirror symmetry for elliptic curves was proved by Polishchuk-Zaslow \cite{PZ-E} based on T-duality, which used derived category of coherent sheaves on the B-side.  Moreover the construction of the mirror map by using homological mirror symmetry was investigated by Zaslow \cite{Zaslow} and Aldi-Zaslow \cite{AZ}.  Namely they derive the Seidel's mirror map for elliptic curves and Abelian varieties by matching the morphisms of the generators of the Fukaya category and that of the category of coherent sheaves in the mirror.

In this paper we have different purposes and use different tools.  Namely, we do not start with a given mirror but rather reconstruct the mirror family and the mirror functor, and the mirror map should be produced as a consequence.  We use the machinery of immersed Lagrangian Floer theory on the A-side and Landau-Ginzburg model on the mirror B-side for this purpose.  The work of Orlov \cite{Orlov} is crucial to pass from Landau-Ginzburg model to Calabi-Yau geometry (or more generally geometry of singularities).
\end{remark}

\begin{remark}
There were several interesting physics literatures studying open mirror symmetry for elliptic curve quotients, for instance \cite{B,GJ,HLN,GJLW,KO}.  They identified the effective superpotential coming from summing up disc instantons with the B-model partition function of the mirror via an open-closed mirror map.  

The study of open mirror symmetry in our paper takes a rather different perspective.   Instead of starting with the Landau-Ginzburg superpotential given by the physicists, we construct the mirror Landau-Ginzburg model by a generalized SYZ construction using immersed Lagrangian Floer theory.  The two superpotentials differ by quantum corrections mentioned before.  Moreover we study the enumerative meaning of the mirror map rather than that of the effective superpotential.  The enumerative statement in Theorem \ref{mir=SYZ1} did not appear in existing physics literatures.  

Furthermore the matrix factorizations we studied are obtained by transforming Lagrangian branes using our mirror functor, rather than a B-model construction in physics literatures.  In other words, we focus more on the constructive and functorial aspects of mirror transformations rather than a comparison of partition functions by hand.  We show that mirror symmetry, which matches the disc instantons with the mirror map and Fukaya category with category of matrix factorizations of the mirror, naturally follows from functoriality.
\end{remark}

\begin{remark} \label{rmk:Ta}
Mirror symmetry for $\mathbb{P}^1_{a,b,c}$ is a very interesting subject and has been intensively studied.   For instance Frobenius structures and integrable systems of PDEs in the B-side were studied by Milanov-Tseng \cite{MT} and Rossi \cite{Rossi} (and the general framework were studied in \cite{Saito1989,DZ,Hertling,DS,DS2} and many others).  Explicit expressions of Saito's primitive forms \cite{Saito80,Saito82,Saito1983} were studied and derived by Ishibashi-Shiraishi-Takahashi \cite{IST} and Li-Li-Saito\cite{LLS}.  Global mirror symmetry were investigated by Milanov-Ruan \cite{Milanov-Ruan}, Krawitz-Shen \cite{Krawitz-Shen} and Milanov-Shen \cite{MS}.

Another half of homological mirror symmetry for weighted projective lines $\bP^1_{a,b,c}$ was formulated by Takahashi \cite{Ta} and studied by Ebeling-Takahashi \cite{ET}, Ueda \cite{U} and Keating \cite{Keating}.  Namely derived categories of weighted projective lines should be equivalent to Fukaya-Seidel category of their Landau-Ginzburg mirrors.  Our construction begins with a given symplectic manifold (instead of a Landau-Ginzberg model) and is more adapted to deduce the half of homological mirror symmetry stated in Theorem \ref{thm:HMS_Pabc}.
\end{remark}

\begin{remark} \label{rmk:Rossi}
For $\frac{1}{a} + \frac{1}{b} + \frac{1}{c} > 1$, namely the spherical case, Rossi \cite{Rossi} proved that as Frobenius manifolds, the big K\"ahler moduli of $\bP^1_{a,b,c}$ (whose tangent bundle is given by quantum cohomology $\bP^1_{a,b,c}$) is isomorphic to the space of Laurent polynomials of the specific form $ xyz + P_1(x) + P_2(y) + P_3(z) $,
where $P_1,P_2,P_3$ are polynomials of degrees $p,q,r$ respectively.  Thus from the perspective of Frobenius-manifold mirror symmetry, the mirrors of $\bP^1_{a,b,c}$ in the spherical case are finite Laurent polynomials.

Our paper constructs the Landau-Ginzburg mirror via immersed Lagrangian Floer theory, which is roughly speaking counting holomorphic polygons.  When $\frac{1}{a} + \frac{1}{b} + \frac{1}{c} > 1$, there are only finitely many holomorphic polygons, and hence $W$ has finitely many terms.  This gives a different perspective that mirrors of $\bP^1_{a,b,c}$ in the spherical case are finite Laurent polynomials.
\end{remark}

More introductions to the backgrounds and explanations on constructions and proofs are in order.

\subsection{Generalized SYZ mirror construction}

Strominger-Yau-Zaslow (SYZ) \cite{SYZ96} proposed that mirror symmetry can be understood in terms of duality of special Lagrangian torus fibrations.  Namely, the mirror manifold can be constructed by taking fiberwise torus dual of the original manifold, and Lagrangian branes can be transformed to coherent sheaves on the mirror by a real version of Fourier-Mukai transform.

A lot of efforts have been devoted to the SYZ construction of mirrors, where the main difficulty lies in quantum corrections coming from singular fibers: fiberwise torus duality away from discriminant loci only gives the first order approximation of the mirror, and one needs to capture the additional information of holomorphic discs emanated from singular fibers in order to reconstruct the genuine mirror. 
In toric cases, Lagrangian Floer potential from counting of holomorphic discs defines a Landau-Ginzburg mirror (\cite{CO},\cite{FOOOT}).  

When the discriminant loci of the Lagrangian fibration are relatively simple such as in the case of toric Calabi-Yau manifolds or their conifold transitions, quantum corrections by holomorphic discs have a neat expression and so the SYZ construction can be explicitly worked out \cite{CLL,AAK,L13}.  In general the discriminant locus of a Lagrangian fibration is rather complicated, these holomorphic discs interact with each other and form complicated scattering diagrams studied by Kontsevich-Soibelman \cite{KS} and Gross-Siebert \cite{GS07}.  Deriving closed-string mirror symmetry and homological mirror symmetry from this perspective is a highly-nontrivial open problem.

In this paper we generalize the SYZ approach so that an $\AI$-functor for homological mirror symmetry naturally comes out, and we call this generalized SYZ (see Remark \ref{rmk:genSYZ}).  The original SYZ approach is based on $T$-duality in which tori play a central role.  In our generalized approach Lagrangians we take are not necessarily tori: indeed they are immersed with transverse self intersections.  The notion of \emph{weak unobstructedness} by \cite{FOOO} plays a key role in our construction.  Solving the Maurer-Cartan equation for weakly unobstructed deformations is the key step to apply our construction to actual situations.

Recall that for a torus $T$, the dual torus $T^*$ is given by
$$ T^* = \{\nabla: \nabla \textrm{ is a flat } U(1) \textrm{ connection on } T \} = H^1(T,\real) / H^1(T,\integer)$$
which is the imaginary part of the space of complexified Lagrangian deformations of $T$:
$H^1(T,\real) \oplus \consti (H^1(T,\real) / H^1(T,\integer)).$
This motivates the following construction: for a Lagrangian immersion $L$, we consider the space $V$ of its first-order weakly unobstructed (complexified) deformations, and let $W$ be the generating function on $V$ of $J$-holomorphic polygon countings bounded by $L$. More precisely, we consider
weak bounding cochains, which come from  linear combinations of immersed generators, and 
consider the associated (immersed) Lagrangian Floer potential. Then $(V,W)$ forms a Landau-Ginzburg model, and we call this to be a generalized SYZ mirror of $X$.

\begin{figure}[htp]
\begin{center}
\includegraphics[scale=0.4]{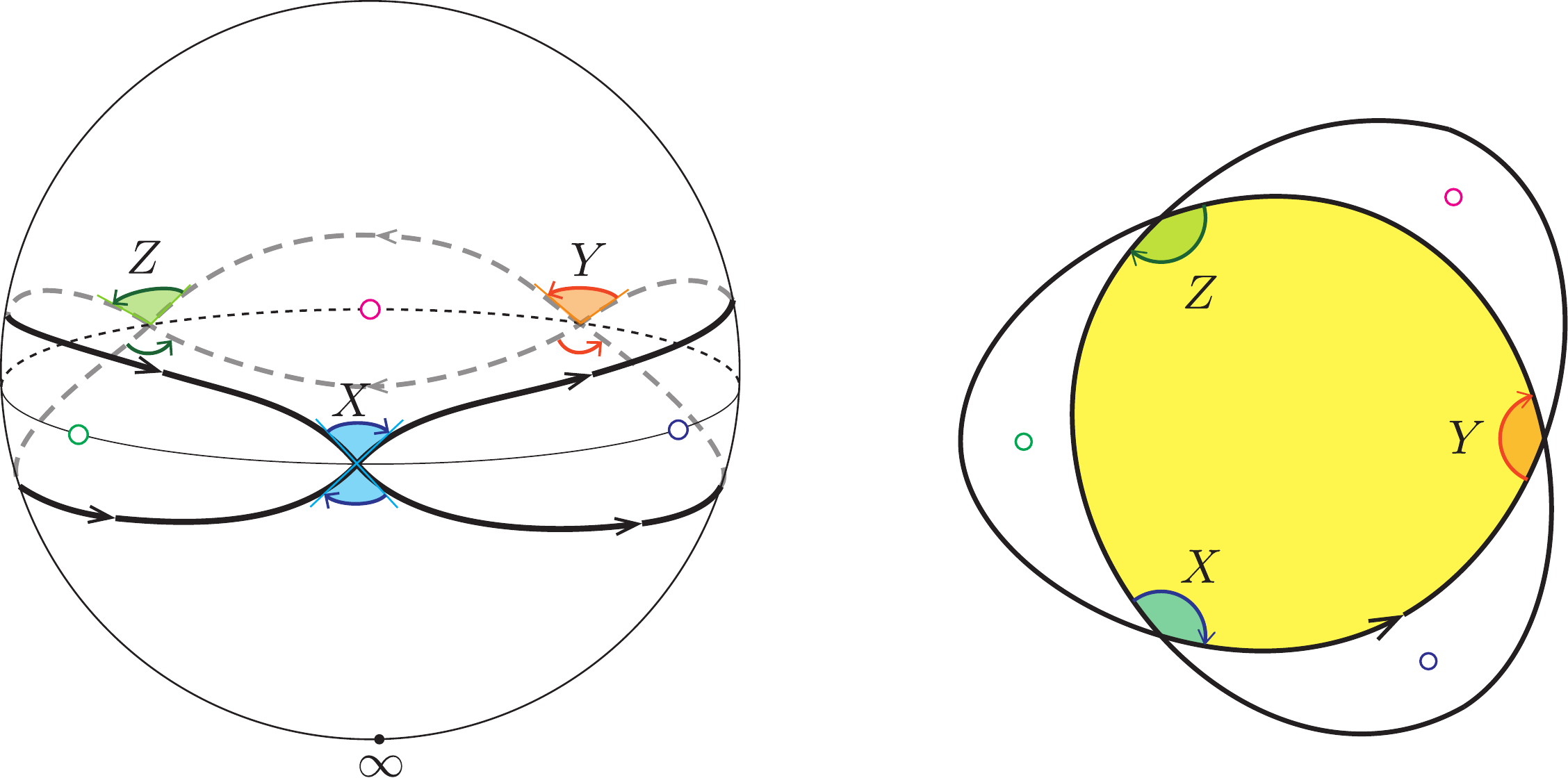}
\caption{The Seidel Lagrangian : Two pictures show the same Lagrangian immersion from different viewpoints.  The three dots on the equator are punctured when $X$ is a pair-of-pants, or they are orbifold points when $X$ is an orbifold projective line.  The shaded triangle on the right contributes to the term $xyz$ of the mirror superpotential.}\label{Seidel_Lag}
\end{center}
\end{figure}

The best example to illustrate the construction is the two-dimensional pair-of-pants $X=\bP^1-\{p_1,p_2,p_3\}$.  Seidel \cite{Se} introduced a specific Lagrangian immersion $L \subset X$ shown in Figure \ref{Seidel_Lag}, and we will call it the Seidel Lagrangian.

$L$ has three immersed points, and they give three directions of (weakly) unobstructed deformations labelled by $x,y,z$.  The only (holomorphic) polygon passing through a generic point of $L$ with $x,y,z$'s as vertices and having Maslov-index two is the triangle shown in Figure \ref{Seidel_Lag}, which corresponds to the monomial $xyz$.  Thus the generalized SYZ mirror of the pair of pants is $W:\C^3 \to \C$ given by
$ W = xyz.$

In \cite{Se},  Seidel used the same type of Lagrangian to prove homological mirror symmetry for genus-two curves (which indeed works for all genus shown by Efimov \cite{Efimov}).  Later Sheridan \cite{Sheridan11,Sh} generalized the construction to higher dimensions and prove homological mirror symmetry for Fermat-type hypersurfaces.  The construction proposed in this paper is largely motivated by their works.

As an application, we apply our construction to the orbifold projective line $\mathbb{P}^1_{a,b,c}$ to produce its Landau-Ginzburg mirror $W$.  The Seidel Lagrangian $L$ plays an essential role in the construction.  As mentioned before, the key step in our construction is solving the Maurer-Cartan equation for weakly unobstructed deformations (Theorem \ref{thm:weak}).  The anti-symplectic involution on $\mathbb{P}^1_{a,b,c}$ is the main geometric input in our method.

While $\mathbb{P}^1_{a,b,c}$ is only one-dimensional, it has a very rich geometry and the quantum-corrected Landau-Ginzburg mirror has an interesting expression.  The leading terms of $W$ are 
$$ - q_\alpha xyz + q_\alpha^{3a} x^a + q_\alpha^{3b} y^b + q_\alpha^{3c} z^c,$$
and there are higher-order terms corresponding to more holomorphic polygons.  The parameter $q_\alpha$ relates with the K\"ahler parameter $q$ (corresponding to the area of $\mathbb{P}^1_{a,b,c}$) by $q = q_\alpha^8$.
When $1/a + 1/b + 1/c \geq 1$, which corresponds to the case that $\Sigma$ has genus less than or equal to one, the superpotential $W$ has finitely many terms; when $1/a + 1/b + 1/c < 1$ (hyperbolic case), which corresponds to the case that $\Sigma$ has genus greater than one, it has infinitely many terms. We determine $W$ explicitly for $\mathbb{P}^1_{3,3,3}$. For general $\mathbb{P}^1_{a,b,c}$, we computed the whole expression of $W$
inductively in a joint work with Kim \cite{CHKL14} using combinatorial techniques.  

Note that $\mathbb{P}^1_{a,b,c}$ can be written as a $G$-quotient of a Riemann surface $\Sigma$.  When $G$ is abelian,  the generalized SYZ mirror of the Riemann surface $\Sigma$ is $(\C^3/\WH{G},W)$ where $\WH{G}$ is the group of characters of $G$ (which is isomorphic to $G$ itself).  Thus it gives a way to construct the quantum-corrected mirrors of higher-genus Riemann surfaces.

As pointed out in Remark \ref{rmk:Rossi}, Rossi \cite{Rossi} proved that the Frobenius structure of quantum cohomology of $\mathbb{P}^1_{a,b,c}$ in the spherical case $1/a + 1/b + 1/c > 1$ is isomorphic to the deformation space of a specific Laurent polynomial (instead of being a series).  This gives another perspective about the Landau-Ginzburg mirror in the spherical case.

% This is in analogous to the toric case, where the disc potential has finitely many terms if and only if the compact toric manifold (or orbifold) is semi-Fano.  In any case, 

\subsection{Localized mirror functor}
Homological mirror symmetry conjecture states that Lagrangian submanifolds (with additional data) correspond to matrix factorizations in the Landau-Ginzburg mirror.  Currently the main approach to prove homological mirror symmetry is to compare generators and their relations (hom spaces) on both sides and show that they are (quasi-)isomorphic.

Constructing a natural functor which transforms Lagrangian submanifolds into matrix factorizations will greatly improve our understanding of (homological) mirror symmetry.  Our generalized SYZ using an immersed Lagrangian $\BL$ provides such a natural $\AI$-functor $\mathcal{LM}^\BL$
from the Fukaya category of $X$ to the category of matrix factorizations of $W$.  
From Theorem \ref{thm:intro:equiv}, under reasonable assumptions the functor derives an equivalence and hence realize homological mirror symmetry.

In general our functor $\functor$ is not an equivalence.  The `mirror' superpotential $W$ that we construct reflects only the local symplectic geometry of $X$ seen by $L$.  Thus we can call it to be a \emph{localized} mirror functor.  While this is weaker than than the mirorr predicted by string theorists, our construction is more flexible and produces generalized mirrors which are not known by string theorists.  The generalized mirrors is particularly useful if we are interested in a particular Lagrangian brane.  See Section \ref{sec:cp1} for an example which shows the local nature of our functor.

In general, one needs to take a finite set of (immersed) Lagrangians in order to generate the Fukaya category, so that the functor $\mathcal{LM}^\BL$ becomes more global and captures more information.  We are developing such a setting in our forthcoming paper, which will involve the use of quivers and non-commutative geometry.

The idea of the construction is the following. It is based on a similarity between the Floer equation and
that of matrix factorization.  As we mentioned earlier, the notion of \emph{weakly ununobstructedness} plays a key role.

One important ingredient is the following observation of Oh \cite{Oh},  which was generalized by Fukaya-Oh-Ohta-Ono \cite{FOOO} to Lagrangian submanifolds deformed by weak bounding cochains.  Suppose that $L_0$ and $L_1$ are two weakly unobstructed Lagrangian submanifolds.  Then the differential $m_1$ of the Floer complex $CF^*(L_0,L_1)$ satisfies
$m_1^2=W_{L_1} - W_{L_0}.$
This follows from $A_\infty$ structure of Lagrangian Floer theory.  More geometrically, it comes from the degenerations of holomorphic strips of Maslov index two with the two boundaries lying in $L_0$ and $L_1$ respectively.
%, which may degenerate into either a broken strip (contributing to $m_1^2$), or a constant strip with a disc bubble. The latter is given by the data of holomorphic discs bounded by either $L_0$ or $L_1$ passing through the intersection point of $L_0$ and $L_1$ where the constant strip is supported.  Thus one obtains the above equality.  %We will consider a more general version when $L_1$ is an immersed Lagrangian that we start with and the compactifications are more complicated in this case.  There could be a non-constant strip of index zero with a disc bubble where the nodal point is immersed, but this additional contributions cancel out when $L_1$ is weakly unobstructed (see Section \ref{ssubsec:magiccancel}).

Recall that the matrix factorization equation of $W$ is given by $\delta^2 = W - \lambda.$
Now if we take $L_1 = \BL$ , the Lagrangian (immersion) that we start with in the beginning with a weak bounding cochain $b$, $L_0 = L$ to be any Lagrangian submanifold with potential $W_{L}=\lambda$, and $\delta = m_1^2$, then the above two equations coincide.  In other words, 
%We take $L_1 = L$ to be the immersed Lagrangian with a weak bounding cochain $b$ (denote $(L,b)$ by $\BL$) such that $m_0^{\BL} = W(b) [L]$.  Under Setup \ref{setup2}, $W(b)$ is the Landau-Ginzburg mirror potential $W$ on a vector space $V$ of $b$'s, and we take another Lagrangian $L_2 $ to be weakly unobstructed with a potential $\lambda$, which is a fixed value.  Then, the Floer equation becomes $(m_1^{\BL, L_2})^2 =W(b) - \lambda$, which is a matrix factorization of $W$.  As a result 
the Floer complex $(CF^*((\BL,b),L),m_1)$ gives a matrix factorization of $W$. %Geometrically $m_1$ is given by counting holomorphic polygons bounded by $\BL$ and $L$. 
Here it is essential that we consider the formal deformation parameter $b$ as variables of $W$ to interpret the Floer complex
as a matrix factorization.

We prove that this definition extends to an $\AI$-functor from the Fukaya $\AI$-category of unobstructed Lagrangians to the dg category of matrix factorizations of $W$.  The proof employs the $\AI$ equations of the Fukaya category.  Our functor is similar to the Hom functor in Yoneda embedding, with the additional input of weak unobstructedness.  %Indeed the statement and the proof generalize naturally to all weakly unobstructed Lagrangian submanifolds, and this gives Theorem \ref{thm:locmirfun1}. 

%We illustrate this localized mirror functor in the case of $\CP^1$ in Section \ref{sec:cp1} by taking a union of two perpendicular equators, and also discuss in detail  the case of an elliptic curve quotient $X=E/(\Z/3)$ in Section \ref{sec:mf333}.

We use this approach to prove that our functor is an equivalence for $X = \bP^1_{a,b,c}$, and thus obtain homological mirror symmetry (Theorem \ref{thm:HMS_Pabc}).  An important step is to compute the matrix factorization mirror to the Seidel Lagrangian.  A priori the mirror matrix factorization takes a rather complicated form.  By doing a non-trivial change of coordinates (Section \ref{sec:transf-Seidel}), we can make it into the following very nice form.
\begin{theorem}[see Corollary \ref{cor:seiext}] 
The Seidel Lagrangian $L$ is transformed to the following matrix factorization under the localized mirror functor $\mathcal{LM}^\BL$:
\begin{equation}
\left(\bigwedge\nolimits^*\langle X,Y,Z \rangle, x X \wedge (\cdot) + y Y \wedge (\cdot) + z Z \wedge (\cdot) + w_x \, \iota_X + w_y \, \iota_Y + w_z \, \iota_Z \right)
\end{equation}
where $w_x, w_y, w_z$ are certain series in $x,y,z$ satisfying $x w_x + y w_y + z w_z = W$. 
\end{theorem}
Employing a result by Dyckerhoff \cite{Dy}, it follows that the above matrix factoriztaion split-generates the derived category of matrix factorizations.

Let us mention some related works. There was also a brilliant approach which directly compares the Fukaya category and its mirror category of matrix factorizations for $\CP^1$ by Chan and Leung \cite{CL}, which were generalized to the case of $\CP^1 \times \CP^1$ in \cite{CHL}. However it is rather difficult to generalize this approach, since Lagrangian Floer differentials may {\em not} have coefficients which are analytic in mirror variables in general.

\begin{remark}
For Lagrangian torus fibrations, Kontsevich-Soibelman \cite{KS-T} has made an attempt to build up a bridge between SYZ and homological mirror symmetry.  There is also an approach using Fourier-Mukai transforms to define mirror functors, which was studied by Tu \cite{Tu} who applied it to toric varieties.  After the publication of this paper as an arXiv preprint, Abouzaid \cite{Abouzaid} used family Floer theory to construct a functor.  However, in the presence of singular fibers in the interior of the base of a Lagrangian torus fibration, one has to carry out order-by-order quantum corrections to the constructions.  How to construct the mirror in the presence of order-by-order corrections is still an open question.
\end{remark}

\subsection{Enumerative meaning of mirror maps}

The mirror map is a central object in mirror symmetry.  It matches the flat coordinates on the K\"ahler moduli and the mirror complex moduli and is essential in the computation of Gromov-Witten invariants using mirror symmetry.  It arises from the classical study of deformation of complex structures and Hodge structures and can be computed explicitly by solving Picard-Fuchs equations.

Integrality of the coefficients of mirror maps has been studied \cite{zudilin02,lian98,krattenthaler10}.  Conjecturally these integers should have enumerative meanings in terms of disc counting.  In the tropical setting of toric degenerations a conjecture of this type was made in the foundational work of Gross-Siebert \cite{GS07}.

In the compact semi-Fano toric case and toric Calabi-Yau case, the mirror map was shown to be equal to the SYZ map \cite{CLT11,CLLT12,CCLT13}.  Since the SYZ map is written in terms of disc invariants, this gives an enumerative meaning of the mirror map.  However for compact Calabi-Yau manifolds the problem is much more difficult: to the authors' knowledge the precise enumerative meaning of mirror map in terms of open Gromov-Witten invariants was not known before even conjecturally.  (Gross-Siebert \cite{GS07} has a precise conjecture in terms of tropical geometry.)

In this paper, we constructs a Landau-Ginzburg model $W$ for a K\"ahler manifold $X$ by the use of a Lagrangian immersion $L \subset X$.  $W$ is counting polygons weighted by their areas and hence depends on the K\"ahler structure $\omega$ of $X$. As a result, we have a map from the K\"ahler moduli to the complex moduli of the Landau-Ginzburg model $W$.  We call this the generalized SYZ map (see Definition \ref{def:SYZmap}).

We prove that the mirror map exactly coincides with our generalized SYZ map in the case of the elliptic curve quotient $\mathbb{P}^1_{3,3,3}= E/(\Z/3)$ (Theorem \ref{mir=SYZ1}), where
$ E = \{x^3 + y^3 + z^3 = 0\}.$
Thus we obtain an enumerative meaning of the mirror map of elliptic curves (which are mirror to the covering elliptic curve $E$).

%The Seidel Lagrangian (Figure \ref{Seidel_Lag}) in $E/(\Z/3)$ lifted to $E$ is a union of three circles, see Figure \ref{qqalpha}.  
In this case the generalized SYZ mirror takes the form (after change of coordainates)
$ W = (x^3 + y^3 + z^3) - \frac{\psi(q)}{\phi(q)} xyz $
where $\phi(q)$ and $\psi(q)$ are generating series counting triangles bounded by the union of three circles with vertices at $x,x,x$ and at $x,y,z$ respectively.
%, see Figure \ref{Deltax3xyz}.  
These generating series are computed explicitly in Section \ref{ex333} (note that their coefficients have signs which requires careful treatments).

On the other side, let $\pi_A (\check{q})$ and $\pi_B (\check{q})$ be the periods of the elliptic curve $E$ which satisfy the Picard-Fuchs equation
$$ u''(\check{q}) + \frac{3\check{q}^2}{\check{q}^3 + 27}u'(\check{q}) + \frac{\check{q}}{\check{q}^3 + 27}u(\check{q}) = 0. $$
The inverse series of $q(\check{q}) = \pi_B (\check{q}) / \pi_A (\check{q})$ is what we refer to as the mirror map, and it can be explicitly written as
$ \check{q} = -3 a(q)$, where
$$ a(q) = 1 + \frac{1}{3}\left( \frac{\eta(q)}{\eta(q^9)} \right)^3 = 1 + \frac{1}{3} q^{-1} \left( \frac{\prod_{k=1}^\infty (1-q^k)}{\prod_{k=1}^\infty (1 - q^{9k})} \right)^3. $$
We verify that the above series equals to $\frac{\psi(q)}{\phi(q)}$.
Thus the mirror map has an enumerative meaning of counting (holomorphic) triangles.
%shown in Figure \ref{Deltax3xyz}.  
We conjecture that the equality between mirror map and generalized SYZ map continues to hold for Fermat hypersurfaces in all dimensions (Conjecture \ref{conj-mm}), and this would give an enumerative meaning of mirror maps of Fermat hypersurfaces.

%When the Lagrangian immersion $L$ we have chosen to start with split generates the Fukaya category, we expec that the Landau-Ginzburg model $W$ we construct is the mirror.  In particular we expect closed mirror symmetry, which can be stated as
%$$ QH^*(X) \stackrel{\cong}{\longrightarrow} \Jac (W) $$
%where $QH^*(X)$ is the (small) quantum cohomology of $X$ defined by counting rational curves, and $\Jac (W)$ is the Jacobian ring of $W$ defined as a quotient of the function ring by the Jacobian ideal generated by the derivatives of $W$.  The map is defined by taking first variations of $W$ under deformations of K\"ahler structures $\omega$ on $X$.  This gives a representation of the quantum cohomology which is in general difficult to compute.

The organization of this paper is as follows.  Section \ref{GMir} algebraically constructs the $\AI$-functor from the Fukaya category to the category of matrix factorizations and hence proves Theorem \ref{thm:locmirfun1} algebraically.  Section \ref{gen_SYZ} and \ref{sect:geom_fctor} formulate the generalized SYZ construction by employing immersed Lagrangian Floer theory of Akaho-Joyce \cite{AJ}, and define the $\AI$-functor more geometrically.  Our construction is applied to a finite-group quotient in Section \ref{sect:quot}.  Section \ref{ex333} and \ref{sect_dim1} apply our general construction to elliptic curves and $\bP^1_{a,b,c}$, in which Theorem \ref{mir=SYZ1} and \ref{thm:HMS_Pabc} are proved.  Section \ref{ndim} is speculative in nature: we apply our construction to Fermat hypersurfaces and brings out more discussions and questions for future research.

\section*{Acknowledgement}
We thank N.-C. Leung for helpful discussions and encouragements. We would also like to thank K. Ono for the kind explanation of the generation criterion for Fukaya categories. The first author thanks Northwestern University for its hospitality during his sabbatical year stay where this work was started. The second author thanks David Favero for his explanation on the induced functor between derived categories. The third author is grateful to S.-T. Yau for his continuous encouragements and useful advice.  He also expresses his gratitude to E. Zaslow for explaining his work with Polishchuk on categorical mirror symmetry on elliptic curves at Northwestern University, D. Auroux for useful discussions on related works on wrapped Fukaya category at UC Berkeley, and thanks to N. Sheridan for explaining his work on homological mirror symmetry during his visit to the Chinese University of Hong Kong.
The work of C.H. Cho and H. Hong was supported by the National Research Foundation of Korea (NRF) grant funded by the Korea Government (MEST) (No. 2012R1A1A2003117).  The work of S.-C. Lau was supported by Harvard University.

\section{Algebraic construction of localized mirror functor}\label{GMir}
The formalism of localized mirror functor can be understood purely algebraically, which will be explained in this section.
We first recall algebraic notions of filtered $\AI$-category, weak bounding cochains, and potentials of
 weakly unobstructed objects. Filtered $\AI$-subcategory of weakly unobstructed objects can be decomposed into $\AI$-subcategories depending on the value of their potentials. Then we propose a setting, in which we construct a filtered $\AI$-functor 
from $\AI$-subcategory of a given value to a dg-category of matrix factorization. 

In geometric situations, Fukaya category is a filtered $\AI$-category, and this algebraic construction gives a localized mirror functor
between a Fukaya sub-category (of a given potential value) to the associated category of matrix factorization.
It is rather amazing that the entire construction of such a mirror functor comes from the structure maps of the filtered $\AI$-category itself.
This provides a geometric (and explicit) construction of (localized) homological mirror functor .

Our localized mirror functor may be understood as certain kind of Hom functor (see  Section \ref{sect_hom}).  Recall that given an $\AI$-category $\CC$, an object $c$ induces a so-called Hom functor from $\CC$ to the dg category of chain complexes, which is used to construct the Yoneda embedding.  Yoneda embedding provides a way to find differential graded (dg) category which is homotopy equivalent to $\CC$.  The construction of the mirror functor is similar to the Hom functor of Yoneda embedding.

\subsection{Filtered $\AI$-category and weak bounding cochains}
Let us recall the definition of filtered $\AI$-algebra and filtered $\AI$-category and relevant notions
from \cite{Fukaya}, \cite{FOOO}.
\begin{defn}
The Novikov ring $\Lambda_0$ is defined as
$$ \Lambda_0= \left. \left\{ \sum_{j=1}^\infty a_j T^{\lambda_j} \right| a_j \in \C, \lambda_j \in \R_{\geq 0}, \lim_{j \to \infty} \lambda_j = \infty \right\}$$
The universal Novikov field $\Lambda$ is defined as a field of fractions $\Lambda = \Lambda_0[T^{-1}]$.  We define a filtration $F^\cdot \Lambda$ on $\Lambda$ such that for each $\lambda \in \R$,
$F^\lambda \Lambda$ consists of elements $\sum_{j=1}^\infty a_j T^{\lambda_j}$ with $\lambda_j \geq \lambda$ for all $j$.
We denote by $F^+\Lambda$ the set of elements $\sum_{j=1}^\infty a_j T^{\lambda_j}$ with $\lambda_j >0$ for all $j$.
\end{defn}
\begin{defn}
A filtered $\AI$-category $\CC$ consists of a collection $Ob(\CC)$ of objects, $\Z$-graded torsion-free filtered
$\Lambda_0$ module $\CC(A_1,A_2)$ for each pair of objects $A_1,A_2 \in Ob(\CC)$, and the operations
\begin{equation}
m_k:\CC[1](A_0,A_1) \otimes \cdots \otimes \CC[1](A_{k-1},A_k) \to \CC[1](A_0,A_k)
\end{equation}
of degree 1 for $k=0,1,\cdots$ and $A_i \in Ob(\CC)$.
Here $m_k$ is assumed to respect the filtration, and they are assumed to satisfy $\AI$-equations:
For $x_i \in \CC[1](A_{i-1},A_i)$ for $i=1,\cdots n$, we have
\begin{equation}\label{eq:ai}
\sum_{n_1+ n_2=n+1} \sum_{i=1}^{n_1} (-1)^{\epsilon_1}m_{n_1}(x_1,\cdots, x_{i-1},m_{n_2}(x_i,\cdots, x_{i+n_2-1}), x_{i+n_2},
\cdots, x_n)=0\end{equation}
where $\epsilon_1 = \sum_{j=1}^{i-1} (|x_j|+1)$.
\end{defn}

\begin{remark}
We may use $\Z_2$-grading instead of $\Z$-grading in the above definition, which is the case when we consider non-Calabi-Yau geometries.
\end{remark}

\begin{defn}
A filtered $\AI$-category with one object is called a filtered $\AI$-algebra.
\end{defn}
We denote by $|x|$ the degree of $x$, and by $|x|'$ the shifted degree of $x$,
with $|x|' = |x|-1$.
\begin{defn}
A filtered differential graded category $\CC$ is a filtered $\AI$-category with vanishing $m_{\geq 3}$ and $m_0$.
\end{defn}

The sign convention for $\AI$-categories is different from the standard one for differential graded category,
and one can define differential $d$ and composition $\circ$ as
\begin{equation}\label{eq:dgsign}
d(x) = (-1)^{|x|} m_1(x),\;\; x_2 \circ x_1 = (-1)^{|x_1|(|x_2|+1)}m_2(x_1,x_2).
\end{equation}
to recover the standard convention of dg-category.

We now explain (weak) bounding cochains (from \cite{FOOO}) defined for a single object in a filtered $\AI$-category.
\begin{defn}
An element $e_A \in \CC^0(A,A)$ for $A \in Ob(\CC)$ is called a {\em unit}
if $m_2(e_A,x_1)=x_1$, $m_2(x_2,e_A) =(-1)^{|x_2|}x_2$  for $x_1 \in \CC(A,A_1), x_2 \in \CC(A_1,A)$,
and
if $m_{k+l+1} (x_1,\cdots, x_l, e_A,y_1,\cdots, y_k) =0$ hold for $k + l \neq 1$.
\end{defn}
For $A \in Ob(\CC)$, we define
$$B\CC[1](A) = \bigoplus_{k=1}^\infty B_k\CC[1](A) = \bigoplus_{k=1}^\infty \CC[1](A,A)^{\otimes k}$$
As usual, $m_k$ defines a coderivation $\WH{d}_k:B\CC[1] \to B\CC[1]$.
Filtration on $\Lambda$ induces a filtration on $\CC[1](A,A)$, and also on $B\CC[1](A)$.
Let us denote its completion by $\WH{B}\CC[1]$, and $\WH{d}_k$ extends to $\WH{B}\CC[1]$.
We denote $\WH{d} = \sum_{k=0}^{\infty} \WH{d}_k.$
\begin{defn}
An element $b \in F^+\CC^1(A,A)$ is called a {\em bounding cochain} (Maurer-Cartan element)
if $\WH{d}(\conste^b) =0$
where  we denote $b^k =  b\otimes \cdots \otimes b \in B_k\CC[1](A)$ and  
\begin{equation}\label{eqeb}
\conste^b := 1 +b + b^2 + b^3 + \cdots \in \WH{B}\CC(A,A).
\end{equation}

We denote by $\WT{\CM}(A)$ the set of all bounding cochains of $A$.
\end{defn}
\begin{remark}
The condition of filtration $b \in F^+\CC^1(A,A)$ is given to ensure the convergence of the infinite sum
\eqref{eqeb}. One can consider $b \in F^0\CC^1(A,A)$ if convergence can be ensured.
One simplest example is $b=0$, where $0$ is a bounding cochain if $m_0=0$ since $\WH{d}(e^0) = m_0 $ for some object $A \in Ob(\CC)$.
%In our case, we will also need $b \in F^0\CC^1(A,A)$. 
In the toric case, the filtration zero part was interpreted as holonomy contribution of the associated complex line bundle
to overcome this issue (\cite{C},\cite{FOOOT}).
\end{remark}
Bounding cochains are introduced to deform the given $\AI$-algebra or $\AI$-category.

\begin{theorem}[Proposition 1.20 \cite{Fukaya}]\label{thm:catbding}
Given a filtered $\AI$-category $\CC$, we can define another filtered $\AI$-category
$\CC'$ as follows.
$$Ob(\CC')  = \bigcup_{A \in Ob(\CC)} \{A\} \times  \WT{\CM}(A),$$
$$\CC'((A_1,b_1), (A_2,b_2))=\CC(A_1,A_2)$$
with the operations $m_k^{b_0,\cdots,b_k}: \CC[1](A_0,A_1) \otimes \cdots \otimes \CC[1](A_{k-1},A_k) \to \CC[1](A_0,A_k)$:
$$m_k^{b_0,\cdots,b_k}:=\sum_{l_0,\cdots,l_k} m_{k+l_0+\cdots+l_k}(b_0^{l_0},x_1,b_1^{l_1},
\cdots,b_{k-1}^{l_{k-1}},x_k,b_k^{l_k}).$$
\end{theorem}

In \cite{FOOO}, Fukaya, Oh, Ohta and Ono constructed filtered $\AI$-algebras from Lagrangian submanifolds, which
has $m_0 \neq 0$ in general.  If an $\AI$-algebra has a bounding cochain $b$, then the deformed $\AI$-algebra $\{m_k^b\}$ has vanishing $m_0^b$, and hence Lagrangian Floer cohomology can be defined since $(m_1^b)^2 = 0$.

More generally we consider the notion of a weak bounding cochain, which can be used for the deformation
as in Theorem \ref{thm:catbding}.
\begin{defn}\label{defn:potential}
An element $b \in F^+C^1(A,A)$ is called a weak bounding cochain
if $\WH{d}(\conste^b)$ is a multiple of unit $e$. i.e. 
$$\WH{d}(\conste^b) = PO(A,b) \cdot e, \;\; \textrm{for some} \;\; PO(A,b) \in \Lambda.$$
We denote by $\WT{\CM}_{wk}(A)$ the set of all weak bounding cochains of $A$.
The coefficient $PO(A,b)$ is a function on the set of all weak bounding cochains $b$ called the potential.  $PO(A,b)$ may also be denoted as $W(b)$ later in this paper.
\end{defn}

In geometric applications, $PO(A,b)$ is the superpotential of a Landau-Ginzburg mirror.
For instance given a compact toric manifold, $A$ is taken to be a Lagrangian torus fiber, and $PO(A,b)$ is written in terms of one point open Gromov-Witten invariants of the fiber $A$.  It has been shown that in the nef case $PO(A,b)$ equals to the Givental-Hori-Vafa superpotential under the mirror map (see \cite{CO} for the Fano case, \cite{FOOO} for defining $PO(A,b)$ for the general case, and \cite{chan-lau,CLLT12} for the nef case).
In our paper, $A$ will be taken to be an immersed Lagrangian and $PO(A,b)$ counts $J$-holomorphic polygons whose corners are specified by $b$'s.

Now we introduce $\AI$-subcategories of a fixed potential value in the Novikov ring.
Consider two objects $(A_0, b_0), (A_1, b_1)$ in the filtered category $\CC'$ of Theorem \ref{thm:catbding}, where $b_0$ and $b_1$ are weak bounding cochains.
One of the associated $\AI$-equations is
$$(m_1^{b_0,b_1})^2(x)  + m_2(m(e^{b_0}), x) + (-1)^{ |x|'}m_2(x, m(e^{b_1})) = 0.$$
Note that if $PO(A_1,b_0) = PO(A_2,b_1)$, then the latter two terms cancel out from the definition of a unit.
This shows that we have $(m_1^{b_0,b_1})^2(x) = 0$ if and only if $PO(A_1,b_0) = PO(A_2,b_1)$.

\begin{defn}
A full filtered $\AI$-subcategory $\CC'_\lambda$ of $\CC'$  for $\lambda \in \Lambda$ is defined by
setting 
$$ Ob(\CC'_\lambda) = \{ (A, b) \in Ob(\CC') \mid b \textrm{ is a weak bounding cochain and } PO(A,b) = \lambda \}.$$
For any two objects in $\CC'_\lambda$, $m_1$ gives a differential on  morphisms between them.
$\AI$-subcategory of unobstructed objects is $\CC'_0$.
\end{defn}

\subsection{Hom functor} \label{sect_hom}
We briefly review the notion of Hom functor in Yoneda embedding for $\AI$-category following \cite{Fukaya}.
In general, for an ordinary category $\mathcal{K}$ and an object $K \in \mathcal{K}$,
$Hom(K,\cdot)$ defines a functor from $\mathcal{K}$ to the category of sets.
For an $\AI$-category $\CC$ and an object $A \in Ob(\CC)$, Fukaya defined a hom functor
$Hom(A,\cdot)$ so that it is $\AI$-functor from the category $\CC$ to
the dg category of chain complexes.

First recall the notion of an $\AI$-functor.
We put 
$$B_k\CC[1](A,B) = \bigoplus_{A=A_0, A_1,\cdots,A_{k-1}, A_k=B} \CC[1](A_0,A_1)\otimes \cdots \otimes \CC[1](A_{k-1},A_k).$$
$$B\CC[1](A,B) = \bigoplus_{k=1}^\infty B_k\CC[1](A,B), B\CC[1] = \bigoplus_{A,B} B\CC[1](A,B).$$
The $\AI$-operation $m_k$ also induces coderivations $\WH{d}_k$ on $B\CC[1]$, and also on its completion $\WH{B}\CC[1]$.
The system of $\AI$-equations \eqref{eq:ai} can be written as a single equation: $\WH{d} \circ \WH{d} =0$.

\begin{defn}
Let $\CC_1, \CC_2$ be $\AI$-categories. An $\AI$-functor $\CF:\CC_1 \to \CC_2$ is a collection of $\CF_k, k\in \Z_{\geq 0}$
such that $\CF_0 : Ob(\CC_1) \to Ob(\CC_2)$ is a map between objects, and
for $A_1,A_2 \in Ob(\CC_1)$, 
$\CF_k(A_1,A_2): B_k\CC_1(A_1,A_2) \to \CC_2[1](\CF_0(A_1),\CF_0(A_2))$ is a homomorphism of
degree 0. 
The induced coalgebra homomorphism $\WH{\CF}_k:B\CC_1[1] \to B\CC_2[1]$ is required to be a chain map
with respect to $\WH{d}$ where  $\WH{\CF}_k(x_1\otimes \cdots \otimes x_k)$ is given by
$$ \sum_m \sum_{0=l_1<l_2< \cdots < l_m=k}
\CF_{l_2-l_1-1}(x_{l_1+1}\otimes \cdots \otimes x_{l_2}) \otimes \cdots \otimes 
\CF_{l_{m}-l_{m-1}-1}(x_{l_{m-1}+1}\otimes \cdots \otimes x_{l_m}).$$
\end{defn}

Let us recall the dg category of chain complexes $\mathcal{CH}$ (see Proposition 7.7 of \cite{Fukaya}).

\begin{defn}
The set of objects $Ob (\mathcal{CH})$ is the set of all chain complexes of free $\Lambda$-modules.
 For $(C,d), (C',d') \in Ob (\mathcal{CH})$,
$$\mathcal{CH}^k ( (C,d), (C',d')) = \bigoplus_l \Hom_\Lambda (C^l, C'^{l+k} ).$$
We define
$$m_1 (x) = d' \circ x - (-1)^{|x|} x \circ d$$
$$m_2 (x,y) =  (-1)^{|x| ( |y| +1 ) } y \circ x$$
and put $m_k = 0$ for $k \geq 3$.  This defines a dg category (as a special case of $\AI$-category).
\end{defn}

We introduce the following notations:
set $\bx = x_1\otimes \cdots\otimes x_k$, and
$$|\bx| = \sum_j |x_j|,\;\;\; |\bx|' = |\bx|-k.$$
Moreover, if $\Delta: B\CC \to B\CC \otimes B\CC$ is
 defined by $$ \Delta(\bx) = \sum_{i=0}^k (x_1 \otimes \cdots \otimes x_i) \otimes(x_{i+1} \otimes \cdots
 \otimes x_k),$$
 then we write $$\Delta^{m-1} (\bx) = \sum_a \bx_a^{(1)} \otimes \cdots  \bx_a^{(m)},$$
 with $\Delta^{m-1}=(\Delta \otimes 1 \otimes \cdots \otimes 1) \circ \cdots \circ (\Delta \otimes 1\otimes 1) \circ (\Delta \otimes 1) \circ \Delta$.
In this notation the $\AI$-equation can be written as
$$\sum_a (-1)^{| \bx^{(1)}_a|'} m(\bx^{(1)}_a, m(\bx^{(2)}_a),\bx^{(3)}_a) =0.$$

% $\epsilon^o(\bx) = \sum_{1 \leq i <j\leq k} (deg(x_i)+1)(deg(x_j)+1)$.
%\begin{defn}
%The opposite $\AI$-category $\CC^o$ of a given $\AI$-category $\CC$ is defined as follows.
%\begin{enumerate}
%\item  $Ob(\CC^o) = Ob(\CC)$.
%\item For $c,c' \in Ob(\CC^o) = Ob(\CC)$, we put $\CC^0(c,c')=\CC(c',c)$.
%\item The operations $m_k^o$ of $\CC^o$ is defined by
%$$m_k^o(x_1,\cdots,x_k)=(-1)^\epsilon m_k(x_k,\cdots,x_1)$$
%where $\epsilon = 1 + \sum_{1 \leq i <j\leq k} (deg(x_i)+1)(deg(x_j)+1)$.
%\end{enumerate}
%\end{defn}
%One can check that this defines an $\AI$-category.

Now we recall the notion of Hom functor for $\AI$-categories.
The following setup is slightly different from that of Definition 7.16 \cite{Fukaya} (see Remark \ref{rmk:hom} below).
\begin{defn}\label{defncaval}
For $A \in Ob(\mathcal{C})$, the $A_\infty$ functor 
$$\mathcal{F}^A = \Hom (A, \cdot) : \mathcal{C} \to \mathcal{CH}$$  is defined as follows.
$\CF_0^A$ is a map between objects: for $B \in Ob(\mathcal{C})$,
$$ \mathcal{F}_0^A (B) = (\CC(A,B), m_1). $$
$\CF_k^A(x_1,\cdots,x_k)$ is a morphism between
chain complexes defined as
$$ \mathcal{F}_k^A (x_1, \cdots, x_k) (y) = (-1)^{\epsilon_2} m_{k+1} (y, x_1, \cdots, x_k )$$
where $y \in \CC(A,B)$, $B_1=B,B_2, \cdots, B_{k+1}  \in Ob (\mathcal{C} )$, $x_i \in  \CC(B_i, B_{i+1})$
and $\epsilon_2 =|y|' |\bx|'$.
\end{defn}
\begin{remark} \label{rmk:hom}
There are similar functors $\CC \to \mathcal{CH}^o, \CC^o \to \mathcal{CH}$ or
$\CC^o \to \mathcal{CH}^o$.  The above definition matches our purpose of defining the mirror functor.  In \cite{Fukaya}, $\mathcal{F}_k^A (\bx) (y) =\pm m_{k+1}(y, \bx^{o})$
for $\bx^o = x_k \otimes \cdots \otimes x_1$ which is different from the above.
\end{remark}
$\CF_0^A$ associates an object $B$ of $\CC$ with a chain complex $(\CC(A,B),m_1)$, and
for $x_1 \in \CC(B_1,B_2)$, $\CF_1^A(x_1)$ defines a homomorphism between two chain complexes
$\CC(A,B_1) \to \CC(A,B_2)$, which is given by multiplication with $x_1$ from the right:
$\pm m_2(\cdot, x_1)$. Now the category $\CC$ can have non-zero $m_k$'s for $k \geq 3$, while for the category $\mathcal{CH}$ $m_k = 0$ for $k \geq 3$. One can check that $\CF^A$ is an $\AI$-functor from the $\AI$-equations of $\CC$.  It is similar to the proof of Theorem \ref{thm_mir_fctor} below that the localized mirror functor is an $\AI$-functor.

\subsection{Localized mirror functor} \label{sect_alg_fctor}
We first recall the notion of a matrix factorization.

\begin{defn}\label{def:mf}
Let $\mathcal{O}$ be the polynomial algebra $\Lambda [x_1,\cdots,x_m]$.\footnote{When $W$ is a series instead of a polynomial, we need to use a suitable completion $\Lambda [[ x_1 ,\cdots, x_m ]]$ of $\Lambda [x_1,\cdots,x_m]$ instead which includes $W$ as an element.  For the moment, let's assume $W$ to be a polynomial over $\Lambda$.  (The definitions can be easily  generalized to the setting when $W$ is a Laurent polynomial or series.)}  The dg category $\mathcal{MF}(W)$ of matrix factorizations of $W$ is defined as follows.  An object of $\mathcal{MF}(W)$ is a $\Z/2$-graded finite-dimensional free $\mathcal{O}$-module $P=P^0 \oplus P^1$
with an odd-degree $\mathcal{O}$-linear endomorphism $d:P \to P$ such that $d^2 = W \cdot Id_P$.
A morphism from $(P,d_P)$ to $(Q,d_Q)$ is given by an $\mathcal{O}$-module homomorphism $f:P \to Q$.
The category $\mathcal{MF}(W)$ is a differential $\Z/2$-graded category 
with a differential defined on the space of morphisms from $(P,d_P)$ to $(Q,d_Q)$, by 
$$d \cdot f = d_Q \circ f + (-1)^{deg(f)} f \circ d_P,$$
and composition between morphisms is defined as usual.
\end{defn}

To define a functor to the above category of matrix factorizations, consider the following algebraic setup.

\begin{setup}\label{setup}
Let $\CC$ be a filtered $\AI$-category.
Assume that we have a weakly unobstructed object $A \in Ob(\CC)$, with finitely many $\{b_1,\cdots,b_m\}$ 
weak bounding cochains such that all linear combinations
$b=\sum_j x_j b_j$ for $(x_1,\cdots,x_m) \in (\Lambda_0)^m$
are weak bounding cochains.
From this, we may regard potential $PO(A,b)$ as a function
$$W = PO_A: (\Lambda_0)^m \to \Lambda_0.$$
\end{setup}

Our localized mirror functor is defined below under Setup \ref{setup}.  Consider $\lambda \in \Lambda$ and full $\AI$-subcategory $\CC'_\lambda$ from Definition \ref{defncaval} and Theorem \ref{thm:catbding}.  For simplicity, we take $\lambda = 0$,  $\CC_0$ is the full $\AI$-subcategory of all unobstructed objects. For notational convenience, we further assume that any unobstructed object $(L,b)$,  has $b=0$.  The general case when $\lambda \neq 0$ is entirely analogous and is stated in Theorem \ref{thm:functor_non0}.

We will construct an $\AI$-functor from the
$\AI$-category $ \CC_0$  to the dg category of matrix factorizations $\mathcal{MF}(W)$
for the potential $W=PO_A$. This mirror functor is given by the hom functor $Hom( (A,b), \cdot)$ of the deformed $\AI$-category $\CC'$.
The deformed $\AI$-operations $\{m_k^{b,0,\cdots,0}\}_{k=0}^\infty$ defined in Theorem \ref{thm:catbding} play a crucial role.

\begin{defn}\label{def:loc_mirr}
Define an $\AI$-functor $\CF^{(A,b)}$ from the
$\AI$-category $ \CC_0$  to the dg category of matrix factorizations $\mathcal{MF}(PO(A,b))$ by the following:
\begin{enumerate}
\item $\CF^{(A,b)}_0$ sends an object $B$ to the matrix factorization $(\CC((A,b), B), m_1^{b,0})$.
\item $\CF^{(A,b)}_1(x_1)$ is defined by
$$(-1)^{\epsilon_2}m_2^{b,0,0}(\cdot,x_1):(\CC((A,b), B_1), m_1) \to (\CC((A,b), B_2), m_1)$$
for $x_1 \in \CC(B_1, B_2)$ (see Definition \ref{defncaval} for the notation $\epsilon_2$).
\item $\CF^{(A,b)}_k(x_1,\cdots,x_k)$ is defined by
\begin{equation}\label{eq:amfunctor}
 \mathcal{F}_k^{(A,b)} (x_1, \cdots, x_k) (y) = (-1)^{\epsilon_2} m_{k+1}^{b,0,\cdots,0}(y, x_1, \cdots, x_k ).
 \end{equation}
  where $y \in \CC((A,b),B)$, $B_1=B,B_2, \cdots, B_{k+1}  \in Ob (\mathcal{C} )$, $x_i \in  \CC(B_i, B_{i+1})$.
\end{enumerate}
\end{defn}

\begin{theorem} \label{thm_mir_fctor}
The collection of maps $\{\CF^{(A,b)}_*\}$ defines an $\AI$-functor.
\end{theorem}

\begin{proof}
We first show that $(\CC((A,b), B), m_1^{b,0})$ is a matrix factorization.
Note that $(A,b)$ is weakly unobstructed with a potential $m^b_{0,A} = PO(A,b)\cdot e$, and $B \in Ob(\CC_{m_0=0})$ is unobstructed with $m^0_{0,B}=0$.
Hence by the filtered $\AI$-(bimodule) equation for 
$ m_1^{b,0}$, we have
\begin{equation}\label{eq:afunctorfirst}
 m_1^{b,0} \circ m_1^{b,0} (x) + m_2^{b,0}(m^b_{0,A}, x) + (-1)^{deg' x} m_2^{b,0}( x, m^0_{0,B})=0.
\end{equation} 
$$m_1^{b,0} \circ m_1^{b,0} (x) = - m_2^{b,0}(m^b_{0,A}, x) = - W \cdot x.$$
Unraveling the sign convention \eqref{eq:dgsign}, if we set $ d(x) = (-1)^{|x|} m_1^{b,0}(x)$, 
then 
$$d^2(x) = d(d(x)) = (-1)^{|x|}(-1)^{|x|+1} m_1^{b,0} \circ m_1^{b,0} (x) = -(-W \cdot x) = W \cdot x$$

(2) is a part of (3), and hence we prove that $\mathcal{F}_k^{(A,b)}$ defined in Equation \eqref{eq:amfunctor} gives
an $\AI$-functor. 
The $\AI$-functor equation can be written as follows since $\mathcal{MF}$ is dg category and hence $m_{\geq 3=0}$:
\begin{equation}\label{eq:functora1}
\sum_a (-1)^{ |\bx_a^{(1)}|'} \mathcal{F}^{(A,b)}(\bx_a^{(1)}, m(\bx_a^{(2)}), \bx_a^{(3)})(y) =
\end{equation}
\begin{equation}\label{eq:functora2}
\sum_c m_2^{b,0,0}(\mathcal{F}^{(A,b)}(\bx_c^{(1)}), \mathcal{F}^{(A,b)}(\bx_c^{(2)}))(y)
\end{equation}
\begin{equation}\label{eq:functora3}
+  m_1^{b,0} \circ \mathcal{F}^{(A,b)}(\bx) (y) + (-1)^{| \bx|'} \mathcal{F}^{(A,b)}(\bx) \circ  m_1^{b,0} (y)
\end{equation}
The first term can be identified with the following from the definition of $\CF$.
\begin{equation}\label{eq:afunctorp1}
\sum_a (-1)^{(| y|')(| \bx|' +1)}  (-1)^{| \bx_a^{(1)}|'} m_{k+1}^{b,0,\cdots,0}(y, \bx_a^{(1)}, m(\bx_a^{(2)}), \bx_a^{(3)}).
\end{equation}
The additional sign appears when we move $y$ to the front.
The second term \eqref{eq:functora2} can be written as
\begin{equation}
\sum_c (-1)^{| \bx_c^{(1)}|(| \bx_c^{(2)}| +1)} \mathcal{F}^{(A,b)}(\bx_c^{(2)}) \circ
\mathcal{F}^{(A,b)}(\bx_c^{(1)}) (y)
\end{equation}
\begin{eqnarray*}
=& \sum_c (-1)^{| \bx_c^{(1)}|(| \bx_c^{(2)}| +1)}   \mathcal{F}^{(A,b)}(\bx_c^{(2)}) (
(-1)^{\epsilon_3}m^{b,0,\cdot,0}(y, \bx_c^{(1)})) \\
=& \sum_c \sum_c (-1)^{| \bx_c^{(1)}|(| \bx_c^{(2)}| +1)}  (-1)^{\epsilon_3}(-1)^{\epsilon_4}
 m^{b,0,\cdots,0}\big(m^{b,0,\cdot,0}(y, \bx_c^{(1)}), \bx_c^{(2)} \big)\\
\end{eqnarray*}
where
$\epsilon_3=| y|'| \bx_c^{(1)}|' $ and
$\epsilon_4 = (|y|' + | \bx_c^{(1)}|' +1)|  \bx_c^{(2)}|$.

Since 
$$\epsilon_3+\epsilon_4 = | y|'| \bx|' + (| \bx_c^{(1)}|' +1)|\bx_c^{(2)}|'
= | y|'| \bx|' + | \bx_c^{(1)}|'(|\bx_c^{(2)}|'+1)$$
the second term \eqref{eq:functora2} is
\begin{equation}\label{eq:afunctorp2}
(-1)^{| y|'| \bx|'} m^{b,0,\cdots,0}\big(m^{b,0,\cdot,0}(y, \bx_c^{(1)}), \bx_c^{(2)} \big)
\end{equation}
The third term  \eqref{eq:functora3} equals
\begin{equation}
(-1)^{| y|'| \bx|'} m_1^{b,0}m_{k+1}^{b,0,\cdots,0}(y,\bx) +
 (-1)^{| \bx|'}(-1)^{(| y|' +1)| \bx|'}m_{k+1}^{b,0,\cdots,0}(m_1^{b,0}(y),\bx))
\end{equation}
Hence we get
\begin{equation}\label{eq:afunctorp3}
(-1)^{| y|'| \bx|'} \big( m_1^{b,0}m_{k+1}^{b,0,\cdots,0}(y,\bx) +
 m_{k+1}^{b,0,\cdots,0}(m_1^{b,0}(y),\bx) \big).
\end{equation}
Sum of the expressions \eqref{eq:afunctorp1},  \eqref{eq:afunctorp2}, \eqref{eq:afunctorp3}
forms the $\AI$-equation for $(y,\bx)$ after we divide by the common sign $(-1)^{|y|' | \bx|'}$.
This proves that $\CF^{(A,b)}$ is an $\AI$-functor. 
\end{proof}

Let's consider the first few equations of an $\AI$-functor.
The first equation is given in Equation \eqref{eq:afunctorfirst}. The second equation is
$$m_1^{b,0}m_2^{b,0,0}(y,x)
+ m_2^{b,0,0}(m_1^{b,0}(y),x) +(-1)^{|y|'} m_2^b(y,m_1(x)) + m_3^{b,b,0,0}(m_0^b,y,m_1)=0$$
and the last term disappears if $A$ is weakly unobstructed, and the third term disappears if $m_1(x)=0$,
in which case $m_2^{b,0,0}( \cdot, x)$ is a chain map of the matrix factorization given by $m_1^{b,0}$.

The general statement for arbitrary $\lambda \in \Lambda$ is stated below and can be proved similarly.
\begin{theorem} \label{thm:functor_non0}
For each $\lambda \in \Lambda$, we have an $\AI$-functor
$$\CF^{(A,b)}:\CC'_\lambda \to \mathcal{MF}(PO(A,b) - \lambda)$$
defined below:
\begin{enumerate}
\item $\CF^{(A,b)}_0$ sends an object $(B,b_1)$ to a matrix factorization $(\CC((A,b), (B,b_1)), m_1^{b,b_1})$.
\item $\CF^{(A,b)}_1(x_1)$ is a morphism of $\mathcal{MF}$-category defined by $(-1)^{\epsilon_2}m_2^{b,b_1,b_2}(\cdot,x_1)$
$$(\CC((A,b), (B_1, b_1)), m_1) \mapsto (\CC((A,b), (B_2, b_2)), m_1),$$
for $x_1 \in \CC((B_1, b_1), (B_2, b_2))$.
\item $\CF^{(A,b)}_k(x_1,\cdots,x_k)$ is  a morphism of $\mathcal{MF}$-category
defined as
\begin{equation}
 \mathcal{F}_k^{(A,b)} (x_1, \cdots, x_k) (y) = (-1)^{\epsilon_2} m_{k+1}^{b,b_1,\cdots,b_{k+1}}(y, x_1, \cdots, x_k ).
 \end{equation}
    where $y \in \CC((A,b),(B,b_1))$, and $(B,b_1)=(B_1,b_1), \cdots, (B_{k+1}, b_{k+1})  \in Ob (\mathcal{C}'_\lambda )$,
     $x_i \in  \CC((B_i, b_i), (B_{i+1}, b_{i+1}))$.
\end{enumerate}
\end{theorem}

\section{Immersed Lagrangian Floer theory and generalized SYZ construction} \label{gen_SYZ}
In this section we explain what we call ``generalized SYZ'' which employs Lagrangian immersions instead of tori to construct the mirror geometry.  The main tool is the Lagrangian Floer theory for immersed Lagrangians developed by Akaho-Joyce \cite{AJ}, which is parallel to the Floer theory for smooth Lagrangians developed by Fukaya-Oh-Ohta-Ono \cite{FOOO}.  Seidel \cite{Se} and Sheridan's work \cite{Sh} can also be used to develop such a theory.  The actual computations in Section \ref{ex333} and \ref{sect_dim1} will be carried out using Seidel's setting.

Let $(M,\omega,J)$ be a K\"ahler manifold with complex structure $J$ and symplectic structure $\omega$, and $(L,\iota)$ be a Lagrangian immersion into $M$, that is an immersion $\iota : L \to M$ with $\iota^\ast \omega =0$, which does not intersect with any orbifold point.  Denote its image by $\bar{L}$, and by abuse of notations sometimes the immersion $\iota: L \to M$ is referred either as $L$ or $\bar{L}$.  We assume that all self-intersections of $\iota$ are transversal double points, which will be enough for our purpose in this paper.  The theory has a natural generalization to Lagrangian immersions with clean self-intersections.  The manifold $L$ is assumed to be oriented and spin.  

\subsection{The space of Lagrangian deformations} \label{sect:space}

For a smooth compact Lagrangian submanifold $L$, its Lagrangian deformations (up to equivalence by Hamiltonian diffeomorphisms) are parametrized by $H^1(L)$.  More generally, cochains of $L$ give formal deformations of $L$ by deforming its associated $A_\infty$-operations.

For $\bar{L}$ being an immersed Lagrangian, we need to consider the extra deformation directions brought by Lagrangian smoothings of the immersed points of $L$.  Note that each immersed point of $L$ give two independent ways of Lagrangian smoothings.  These smoothings is classically known, see for instance Thomas-Yau \cite{TY}.

To be more precise, consider the fiber product
$L \times_\iota L := \{(p,q) \in L \times L: \iota (p) = \iota  (q)\}$
which consists of several connected components, one of them being the diagonal $\{(p,p) \in L \times L\} \cong L$.  Let $R$ be a set of labels of the connected components, and suppose $0 \in R$ labels the diagonal.  Other than $0 \in R$, elements in $R$ can be identified as ordered pairs $(p_-, p_+) \in L \times L$ such that $\iota(p_-) = \iota (p_+) =p$ and $p_-\neq p_+$.  Denote by $\sigma$ the involution on $R$ sending $0$ to $0$, $(p_-, p_+)$ to $(p_+, p_-)$.  We will sometimes denote  $\sigma(r)$ by  $\bar{r}$.

Elements in $R$ other than $0$ are referred as immersed sectors (in analog to twisted sectors in orbifold cohomology theory developed by Chen-Ruan \cite{CR}).

By \cite{FOOO}, cochains of $L \times_\iota L$ give formal boundary deformations of the immersed Lagrangian $\bar{L}$.  To make an explicit choice of independent directions of formal deformations, we can fix a finite set of cycles $\{\Phi_k\}_{k=0}^N$ in $L$ (with pure degrees) which forms a basis of $H^*(L)$ when descending to cohomology.  Let $\Phi_0$ be $L$ itself.  Let $\{\Phi^k\}_{k=0}^N$ be cycles which form the dual basis when descending to cohomology using Poincar\'e duality of $L$.  Then define
$$ H := \mathrm{Span}_{\cpx} \{\Phi_k: k=0, \ldots, N\} \oplus \mathrm{Span}_{\cpx} (R - \{0\}) \cong H^*(L \times_\iota L) $$
which is a finite-dimensional complex vector space.
In Section \ref{sect_obstr} we will study obstructions to the deformations of $\bar{L}$ and restrict to a smaller subspace $V \subset H$.

\subsection{Review on immersed Lagrangian Floer theory} \label{IMLA}
Akaho and Joyce \cite{AJ} defined a filtered $\AI$-algebra of a Lagrangian immersion $\bar{L}$, which is in principle similar to the construction of \cite{FOOO} for Lagrangian submanifolds.  As discusses in the last section, each self intersection point of $\bar{L}$ gives two immersed sectors.
A motivating case is a family of Lagrangian submanifolds $L_1,\cdots, L_k$ 
which intersect each other transversely, and three of them do not have a common intersection.
Then the $\AI$ algebra of the immersion $\iota:\sqcup_{i=1}^k L_i \to M$ can be defined from the Fukaya category with objects $\{L_1,\cdots, L_k\}$ in an obvious way.
Each intersection point $p \in L_i \cap L_j$ contributes two non-zero elements in $CF(L_i,L_j)$ and $CF(L_j, L_i)$ respectively.

\subsubsection{Moduli spaces of holomorphic discs} \label{sect_mod}
Fix $k \in \integer_{\geq 0}$ and $\beta \in H_2(X,\bar{L})$.  We consider the moduli space of stable discs with $k+1$ boundary marked points representing $\beta$.  Since the boundary data $\bar{L}$ is now immersed, there are additional topological data for a disc that we need to fix in order to have a connected moduli space.  Namely, whether a boundary marked point of the disc is a corner, and if it is, which immersed sector it passes through.  Such data is specified by a map $\alpha : \{0, \ldots, k\} \to R.$
The $i$-th marked point is not a corner if and only if $\alpha(i) = 0 \in R$.  Otherwise it is a corner passing through the immersed sector $\alpha(i)$.

Then we have a connected moduli space $\mathcal{M}_{k+1} (\alpha, \beta)$ of stable discs associated to each $(k, \alpha, \beta)$.\footnote{This is denoted as $\overline{\mathcal{M}}_{k+1}^{main} ( \alpha, \beta, J)$ in the paper of Akaho-Joyce, where $J$ is the almost complex structure we have fixed from the very beginning.}  An element in $\mathcal{M}_{k+1} (\alpha, \beta)$ is a stable map $u : (\Sigma,\partial \Sigma)  \to (M, \iota(L))$ (more precisely an equivalence class of stable maps) from a genus $0$ (prestable) bordered Riemann surface $\Sigma$ with mutually distinct marked points $z_0, \cdots, z_k \in \partial \Sigma$ that are not nodes.  Write $\vec{z} = (z_0, \cdots, z_k)$.  This is an straightforward generalization from the case when $\iota$ is an embedding, with the following additional condition coming from $\alpha$: for all $i$ with $\alpha(i) \not= 0 \in R$,
$$\left( \lim_{z \uparrow z_i} u(z), \lim_{z \downarrow z_i} u(z) \right) = \alpha(i) \in R$$
where the limit is taken over $z \in \partial \Sigma$ where $z \uparrow z_i$ means $z$ approaches $z_i$ in the positive orientation of $\partial \Sigma$.  Since $z \not= z_i$ is sufficiently closed to $z_i$, $u(z)$ is not an immersed point of $\bar{L}$ and can be identified as an element in $L$.  Thus the above limits are taken in $L$.

%To define the notion of stable map from $\Sigma$, we additionally choose a continuous and orientation-preserving map $l : S^1 \to \partial \Sigma$ (up to reparametrization) such that
%\begin{equation}\label{l}
%\left\{
%\begin{array}{l}
%|l^{-1} (x)|=1 \,\, \mbox{if}  \,\, x \in \partial \Sigma \,\, \mbox{is smooth (i.e. not a node);}\\
%|l^{-1} (x)|=2 \,\, \mbox{if} \,\, x \in \partial \Sigma \,\, \mbox{is singular}.
%\end{array} \right.
%\end{equation}
%$u \circ l$ is a continuous map $S^1 \to \iota(L)$ (which may not be able to be lifted to $S^1 \to L$).

%\begin{defn} \label{disc moduli}
%An element of $\mathcal{M}_{k+1} (\alpha, \beta)$ is given by an equivalence class of quintuples $(\Sigma, \vec{z}, u, l, \bar{u})$ such that
%\begin{itemize}
%\item $(\Sigma, \vec{z}, u)$ is a stable $J$-holomorphic map bounding $\iota(L)$ with $u_\ast \left( [\Sigma] \right) =\beta$,
%\item $l : S^1 \to \partial \Sigma$ satisfies \eqref{l},
%\item $\bar{u}: S^1 \setminus \{ \zeta_i : \alpha(i) \not= 0\} \to L$ is a continuous lifting of $u \circ l$ (i.e. $\iota \circ \bar{u} = u \circ l$) where $\zeta_i := l^{-1} (z_i)$,
%\item For all $i$ with $\alpha(i) \not= 0 \in R$,
%$$\left( \lim_{\theta \uparrow 0} \bar{u} \left( e^{i \theta} \zeta_i\right), \lim_{\theta \downarrow 0} \bar{u} \left( e^{i \theta} \zeta_i\right)  \right) = \alpha(i) \in R.$$
%\end{itemize}
%\end{defn}

Define the evaluation maps $ev_i : \mathcal{M}_{k+1} ( \alpha, \beta) \to L \times_\iota L$ by
\begin{equation*}
ev_i \left( [u, \vec{z}] \right) = \left\{
\begin{array}{ll}
(u(z_i),u(z_i)) \in L, & \alpha(i) = 0 \in R\\
\alpha(i) \in R, & \alpha(i) \not= 0 \in R
\end{array} \right.
\end{equation*}
for $i=0, \cdots, k$.  For $i=0$ (the output marked point), we also consider the `twisted' evaluation map $ ev = \sigma \circ ev_0:  \mathcal{M}_{k+1} ( \alpha, \beta) \to L \times_\iota L$, where $\sigma:L \times_\iota L \to L \times_\iota L$ sends $(p_-,p_+)$ to $(p_+,p_-)$:
\begin{equation*}
ev \left( [u, \vec{z}] \right) = \left\{
\begin{array}{ll}
(u(z_0),u(z_0)) \in L, & \alpha(0) = 0 \in R\\
\sigma \circ \alpha(0) \in R, & \alpha(0) \not= 0 \in R.
\end{array} \right.
\end{equation*}

\subsubsection{Maslov index and dimension formula}\label{sec:mas}
It is well-known in Lagrangian Floer theory that the virtual dimensions of moduli spaces are written in terms of Maslov indices of disc classes.  Here we recall the explicit formula in our setting.  For the moment we allow the self-intersections to be clean, which is a little bit more general than requiring them to be transverse.

First recall the Maslov index of a $J$-holomorphic polygon.  Let $L_0, L_1$ be two oriented Lagrangian submanifolds which intersect cleanly.
For $p \in L_0 \cap L_1$, we consider a smooth Lagrangian path $\lambda_{p_0,p_1}$
such that $\lambda_{p_0,p_1}(0) = T_pL_0, \lambda_{p_0,p_1}(1) = T_pL_1$. Without loss of generality we may assume that $\lambda_{p_0,p_1}(t)$ always contain
$T_p(L_0 \cap L_1)$ for all $t$.

For such a path, one can associate an index $\mu(\lambda_{p_0,p_1})$
as follows. Concatenate $\lambda_{p_0,p_1}$ with a positive-definite path from $L_1$ to $L_0$
(this path do not concern orientations of $L_i$'s) to obtain a loop of Lagrangian subspaces in $T_pX$, and
denote the standard Maslov index of this loop by $\mu(\lambda_{p_0,p_1})$. In \cite{FOOO}, it was
shown that this index is related to the Fredholm index of the Cauchy-Riemann operator
on a half-infinite strip with Lagrangian boundary condition given by $\mu(\lambda_{p_0,p_1})$.
Note that the parity of  $\mu(\lambda_{p_0,p_1})$ is independent of the chosen path
$\lambda_{p_0,p_1}$: a different choice of $\lambda_{p_0,p_1}$ results in an even-integer
difference in $\mu(\lambda_{p_0,p_1})$ due to orientation.
One can invert the path $\lambda_{p_0,p_1}$ to obtain $\lambda_{p_1,p_0}$ from
$T_pL_1$ to $T_pL_0$, and the Maslov indices are related by $\mu(\lambda_{p_0,p_1}) + \mu(\lambda_{p_1,p_0}) = n,$
as the concatenation of a positive-definite path contributes $n$ to the Maslov index.  In particular,

\begin{lemma}
If the Lagrangian path $\lambda_{p_0,p_1}$ is homotopic to the positive-definite path
from $L_0$ to $L_1$, then $\mu(\lambda_{p_0,p_1}) =n$. 
\end{lemma}

%\begin{remark}
%Hence, we may want to modify the Definition \ref{positivity}. We need
%a condition that $\mu(\lambda_{p_0,p_1}) =1$.
%\end{remark}

If $L_i$'s for $i=0,1$ are graded, then there is a canonical choice of a path $\lambda_{p_0,p_1}$ (up to homotopy) at an intersection point $p$.  The Maslov index of this path is independent of the choice of $p$ within the same connected component.  Thus each component of intersections is associated with a Maslov index.

For a Lagrangian immersion $\bar{L}$, the above notion of Maslov index applies by taking $L_0,L_1$ to be the two local branches at a self-intersection point $p$ of $\bar{L}$.  Now consider a holomorphic disc $u:(\Delta, \partial\Delta) \to (X, \bar{L})$ bounded by $\bar{L}$ with $k+1$ boundary marked points as in Section \ref{sect_mod}. 
Let $\beta \in H_2(X, \bar{L})$ be the relative homology class of $u$, and $\alpha: \{0,\ldots,k\} \to R$ be a specification of turnings of the boundary marked points.  Consider the pull-back bundle $u^*TX$ and its trivialization over $\Delta$.  Lagrangian sub-bundles pulled back along $\partial \Delta$ can be connected via choices of paths $\lambda_{\alpha(i)}$
for $i=0,\cdots, k$ ($\lambda_{\alpha(i)}$ can be chosen to be a constant path when $\alpha(i) = 0$). Let us denote the resulting Lagrangian loop by $\lambda_{\beta}$, which also depends on choices $\lambda_{L_iL_{i+1}}$.  The difference $\mu(\lambda_{\beta}) -\sum_{i=0}^k \mu(\lambda_{\alpha(i)})$ depends only on $\beta$ and $\alpha$ (rather than the choice of paths $\lambda_{\alpha(i)}$), and is denoted as $\mu(\alpha,\beta)$.

\begin{lemma}
If all self-intersections of $\bar{L}$ are transversal, the virtual dimension of $\moduli_{k+1}(\alpha,\beta)$ is
given by
$$ \mu(\alpha,\beta) + n + \sum_{i=0}^k (1 - \mathrm{codim}(\alpha(i))) - 3$$
where $\mathrm{codim} (r) = 0$ if $r = 0 \in R$, and $\mathrm{codim} (r) = n$ if $r \not= 0 \in R$.

The formula reads the same in the case of clean self-intersections, where $\mathrm{codim}(r)$ is the codimension of the self-intersection component labelled by $r \in R$.
\end{lemma}

%Hence, one may regard $\mu(\lambda_{L_iL_{i+1}})$ as the degree in the normal direction of clean intersection $L_i \cap L_{i+1}$.

%In our transverse situation the dimension formula reads as
%$$ \dim \moduli_{k+1}(\alpha,\beta) = \mu^{\mathrm{poly}}(\beta) + n + \sum_{i=0}^k (1 - \mathrm{codim}(\alpha(i))) - 3$$
%where $\mathrm{codim} (r) = 0$ if $r = 0 \in R$, and $\mathrm{codim} (r) = 1$ if $r \not= 0 \in R$.

%\begin{defn}
%The minimal Maslov index of $(\bar{F}_r,\alpha)$ is defined as
%$$\min \big\{\mu^{\mathrm{poly}}(\beta): [\beta] \in \pi_2(\bar{X},(\bar{F}_r,\alpha)) - \{0\}. \big\}$$
%$(\bar{F}_r,\alpha)$ is said to be admissible if the minimal Maslov index of $(\bar{F}_r,\alpha)$ is at least $2 + \sum_i (n-\dim \alpha|_{\Phi_i}-1)$. 
%\end{defn}

%\begin{prop}
%If $(\bar{F}_r,\alpha)$ is admissible, then disc bubbling does not occur.  Moreover, for every $\beta \in \pi_2(\bar{F}_r,\alpha)$, $[\moduli_1(\beta;\alpha)]_{\mathrm{virt}} = c [F_r]$ for some $c \in \rat$.
%\end{prop}

\subsection{Weakly unobstructedness and the disc potential} \label{sect_obstr}
Given an oriented spin immersed Lagrangian $\bar{L}$, the moduli space $\moduli_{k+1}(\alpha,\beta)$ associated to each disc class $\beta$ and specification of turnings $\alpha: \{0,\ldots,k\} \to R$  at $k+1$ marked points has a compact oriented Kuranishi structure by the work of Akaho-Joyce \cite{AJ}. They used these to define a filtered $\AI$-algebra $A = C^*(L \times_{\iota} L)$ of $\bar{L}$, and we refer the readers to \cite{AJ} for more details.

The $A_\infty$ structure
$\{m_k: A^{\otimes k} \to A: k \geq 0\}$
is, roughly speaking, defined by counting holomorphic polygons with one output marked point and $k$ input marked points.  Recall that in Section \ref{sect:space} we have chosen a finite-dimensional subspace $H \subset A$ spanned by a finite set of cocycles which descend to be a basis in the cohomology level.

The fundamental class $e_L$ of $L$ (where $\iota:L\to\bar{L}$ is the Lagrangian immersion) serves as a unit of the $\AI$-algebra.  As in the algebraic setting of Section \ref{GMir}, we consider the set $\WT{M}_{wk}$ of weak bounding cocycles $b\in H$.

\begin{defn}
$b \in H$ is called to be a weak bounding cocyle if it satisfies
\begin{equation} \label{eqn:weak}
m(\conste^b) = \sum_{k=0}^\infty m_k(b,\ldots,b) = W(b) \cdot e_L.
\end{equation}
The above equation is called the Maurer-Cartan equation for weak bounding cocycles.
\end{defn}

We make the following setup (see Setup \ref{setup} in the algebraic setting).

\begin{setup}\label{setup2}
Take a finite dimensional vector space $V \subset \WT{M}_{wk}$ over the Novikov field $\Lambda_0$.
Namely, consider finitely many linearly independent weak bounding cocycles $b_1,\cdots,b_m$ 
such that all the linear combinations
$b=\sum_j x_j b_j$ for $(x_1,\cdots,x_m) \in (\Lambda_0)^m$
are weak bounding cocycles.  $V$ is the subspace spanned by $b_i$'s.  Then for each $b \in V$ we have $W(b) \in \Lambda_0$ defined by Equation \ref{eqn:weak}.

$(V,W)$ forms a Landau-Ginzburg model, and we call this a \emph{generalized SYZ mirror} of $X$.
%From this, we may regard potential $PO(A,b)$ as a function
%$$W: (\Lambda_0)^m \to \Lambda_0.$$
%{\color{red} We assume that $W$ when restricted to $\C^m$, becomes a holomorphic function
%after replacing the formal parameter $T$ to a suitable complex number.
%Convergence issue as we take $x_i$ from $\Lambda_0$ or $\C$}
\end{setup}

Up to this point we fix the K\"ahler structure $\omega$ of $X$ and define the generalized SYZ mirror $W$, which is in general a formal power series whose coefficients are elements in the Novikov ring $\Lambda_0$.  Now we consider the K\"ahler moduli and define the generalized SYZ map as follows under reasonable assumptions.

\begin{defn}[Generalized SYZ map] \label{def:SYZmap}
Suppose $X$ is compact K\"ahler.  Let $\cM_K(X) = K_C(X) \oplus \consti (H^{1,1}(X,\C) \cap H^2(X,\R)) \subset H^{1,1}(X,\C)$ be the complexified K\"ahler cone of $X$, where $K_C(X)$ denotes the K\"ahler cone.  Fix a smooth lift of $\cM_K(X)$ (or if necessary its universal cover) to the space of K\"ahler forms of $X$, such that $\bar{L}$ is Lagrangian and every $b \in V$ is weakly unobstructed with respect to each complexified K\"ahler form $t \in \cM_K(X)$.  Then we have a formal power series $W_t(x_1,\ldots,x_m)$ which depends on $t$.  

Now assume further that there exists $R \gg 0$ such that for every $t \in \cM_K(X)$ with $\|t\| > R$ and $\|(x_1,\ldots,x_m)\| < R$  (with respect to fixed linear metrics on $(H^{1,1}(X,\C))^2$ and $V$), $W_t(x_1,\ldots,x_m)$ is convergent over $\C$ (by substituting the formal parameter $T$ in the Novikov ring to be a fixed positive real number) and has isolated singularities in $\{\|(x_1,\ldots,x_m)\| < R\}$.  Let $\cM_C(W)$ be the base space of the universal unfolding of $W$ (see \cite{Saito82,DS} for the deformation theory of a function).  Then $t \mapsto W_t$ defines a map
$\cM_K(X) \to \cM_C(W)$, which we call to be the generalized SYZ map.
\end{defn}

The generalized SYZ map is expected to coincide with the (inverse) mirror map, which is a central object in mirror symmetry.  In particular it should satisfy certain Picard-Fuchs equations.  This is proved in \cite{CLLT12,CCLT13} for compact semi-Fano manifolds and toric Calabi-Yau orbifolds respectively.

The superpotential takes the form
\begin{align*}
W(b) %&= \mathop{\sum_{\beta}}_{k \geq 0} q^\beta n(\beta; \underbrace{b, \ldots, b}_{k}) \\
&= \sum_{(i_1, \ldots, i_m) \in \Z_{\geq 0}^m} \sum_{\beta} \sum_{P} q^\beta n(\beta; b_{P(1)}, \ldots, b_{P(k)}) x_1^{i_1} \ldots x_m^{i_m}
\end{align*}
where $b = \sum_{i=1}^m x_i b_i$, $k = i_1 + \ldots + i_m$, and we are summing over all permutations $P$ which are maps $P: \{1,\ldots,k\} \to \{1, \ldots, m\}$ such that $|P^{-1} \{j\}| = i_j$ for all $j=1,\ldots,m$.  Roughly speaking $n(\beta; b_{P(1)}, \ldots, b_{P(k)})$ is counting the number of discs representing the class $\beta$ whose boundary passes through a generic point of $L$ and $b_{j_1}, \ldots, b_{j_k}$, which in general depends on the choice of perturbations in Kuranishi structures of disc moduli.  The weakly unobstructedness of $b$ ensures that while the individual numbers $n(\beta; b_{P(1)}, \ldots, b_{P(k)})$ depend on perturbations, $W$ is well-defined as a whole expression.

\subsection{Fukaya category for surfaces}\label{sec:fuksur}
The theory is already very interesting for $\dim_\R X = 2$. 
For the explicit computation in Section \ref{ex333} and \ref{sect_dim1}, we follow Seidel 
\cite{Se} for the Morse-Bott construction of the Fukaya category of surfaces, which is briefly explained here. 

Let $M$ be a surface, and consider a countable set $\mathcal{L}$ of
closed embedded curves in $M$ as a set of Lagrangian submanifolds.
We assume that  any two curves of $\mathcal{L}$
intersect transversally, and any three curves do not have a common intersection, and also that  curves in $\mathcal{L}$  do not bound a holomorphic disc in $M$. 

\begin{figure}[h]
\begin{center}
\includegraphics[height=1.2in]{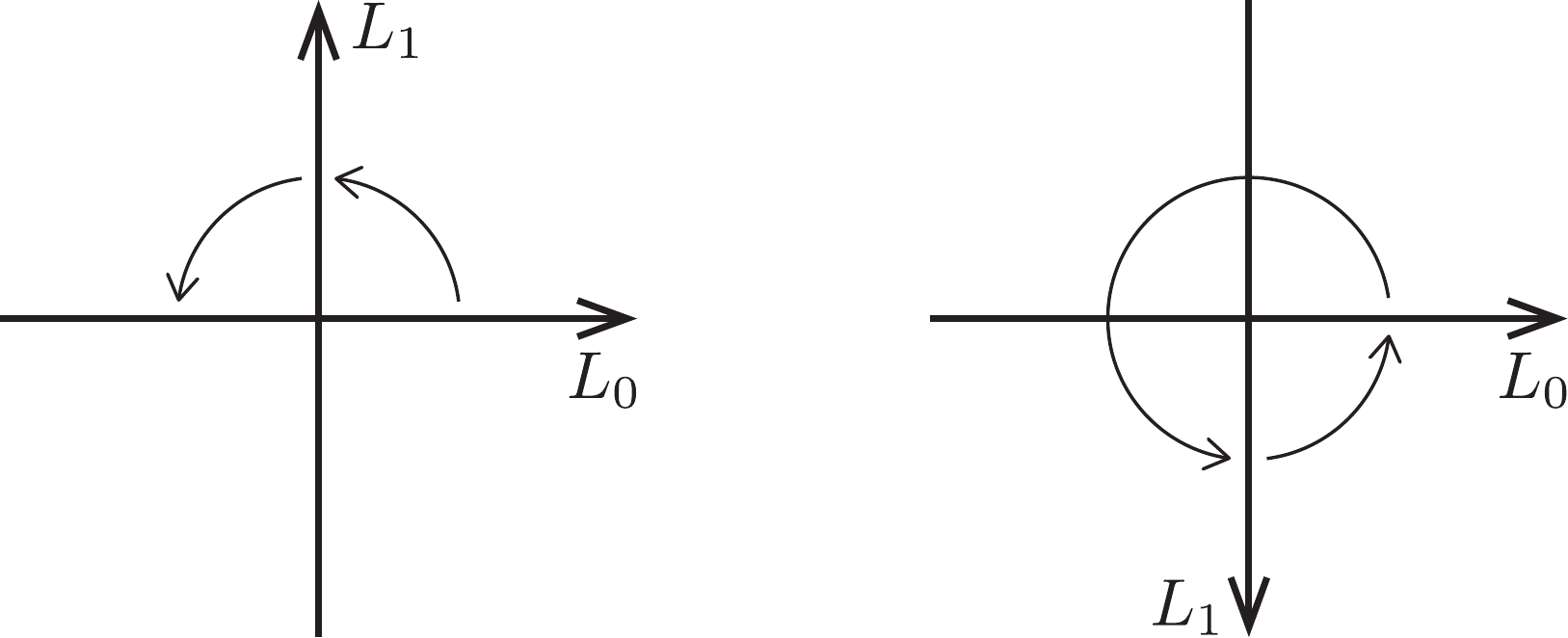}
\caption{odd degree (left) and even degree (right) intersections}\label{figdeg}
\end{center}
\end{figure}

Given two distinct objects $L_0, L_1$ in $\mathcal{L}$, the morphism space between them is generated by intersection points
$$CF^*(L_0,L_1)= \bigoplus_{x\in L_0 \cap L_1} \Lambda x.$$
We will think of them as being $\Z/2$-graded, and degree of a generator is even or odd  according to the rule explained in Section \ref{sec:mas}. The  Figure \ref{figdeg} illustrates the odd and even degree intersection points, which are measured by the Maslov indices of the Lagrangian paths which are described by arrows in the figure.

Boundary operator $\partial $ of $CF^*(L_0,L_1)$ is defined by the counting of 
strips $u: \R \times [0,1] \to M$ which satisfy the $J$-holomorphic equation
$\frac{\partial u}{\partial \tau} + J \frac{\partial u}{\partial t} =0$ and
the boundary condition $u(\tau, 0) \in L_1$, and $u(\tau,1) \in L_0$. 
(We follow the convention of \cite {FOOOanchor}).

To define $\AI$-operations, first consider the case that the collection of objects $(L_0, L_1,\cdots, L_d)$ are pairwise distinct.
In this case the $A_\infty$-operations are defined as
$$m_k:CF(L_{i_0},L_{i_1}) \times \cdots 
\times CF(L_{i_{k-1}},L_{i_k}) \to CF(L_{i_0},L_{i_k}).$$
$$m_k (x_1, \cdots, x_k) := \sum_y \left( \# \mathcal{M} (x_1, \cdots, x_k, y) \right) y,$$
where $\mathcal{M} (x_1, \cdots, x_k, y)$ is the moduli space of holomorphic polygons which turn at $x_1, \cdots, x_k, y$ in counter-clockwise order
and $\#$ means the signed count of such holomorphic polygons whose symplectic areas are recorded in the exponents of $q$'s.

The sign is determined as follows. First, each Lagrangian $L$ may be equipped with a non-trivial spin structure, in which case
we fix a  point ``$\circ$" in $L$ which is regarded as a point where this nontrivial spin bundle on $\bar{L}$ is twisted.
For a holomorphic disc $P \in \mathcal{M}_{k+1} (x_0, \cdots, x_{k})$, if there is no ``$\circ$" on the boundary of $P$ and the orientation of the boundary $P$ agrees with that of the Lagrangian, then the contribution of $P$ to $m_k(x_0, \cdots,x_{k-1})$ is $+q^{\omega(P)} x_k$.
If orientations do not coincide, the resulting sign is as follows. Disagreement of the orientations on $\arc{x_0 x_1}$ is irrelevant to the sign. Disagreement of the orientations on $\arc{x_i x_{i+1}}$ affects the sign by $(-1)^{|x_i|}$. If two orientations are opposite on $\arc{x_k x_0}$ then we multiply by $(-1)^{|x_0| + |x_k|}$. Finally, we multiply by $(-1)^{l}$ where $l$ is the number of times $\partial P$ passes through ``$\circ$". 

Now consider two objects $L_0, L_1$ such that the underlying two curves are the same.  We choose a metric and a Morse function $f_{01}$ on this curve, with a unique maximum $p$ and a unique minimum $e$.  Then $$CF^*(L_0,L_1)= CM^*(f_{01}) =   \Lambda e \oplus \Lambda p.$$
If $L_0, L_1$ have the same orientation and spin structure, the differential on $CF^*(L_0,L_1)$ are Morse differential and hence vanishes.

A general $\AI$-operation, when some of the $L_i$'s coincide, is defined from
the counting of holomorphic pearly trees (using the combination of
Morse flows and $J$-holomorphic polygons). This construction was sketched in \cite{Se} and
details has been given in Section 4 of \cite{Sheridan11}. Such a construction uses
modulus- and domain-dependent perturbations to achieve transversality, and
we refer readers to the above references for more details.
We remark that disc or sphere bubbling do not appear in defining $\AI$-operations
for every surfaces $X$ except a sphere. When $X$ is a sphere, and $L$ is an embedded
curve, it may bound Maslov index two discs and Chern number two spheres, but it can
be handled in the same way with minor modifications.
We remark that our weakly unobstructedness and  computations in Section \ref{ex333} and \ref{sect_dim1} will not   involve complicated holomorphic pearly trees, but rather a counting of  polygons or the following elementary situations.

 Consider the case $(L_0,\cdots, L_k)$ such that $L_i = L_{i+1}$ for one $i \in (\Z/k)$.  We may choose a Morse function on $L_i$ which has only one minimum $e$ and one maximum $p$.  The $\AI$-operation $m_k(\cdots,e,\cdots)$ (or $m_k(\cdots,p,\cdots)$), where $e$ (or $p$) is put in the $(i+1)$-th place for $i=0,\ldots,k-1$, is defined by counting holomorphic polygons whose $i$-th marked points lie in the unstable submanifold of $e$ (or $p$).  Similarly when $i=k$, the coefficient of $e$ (or $p$) in the $\AI$-operation $m_k(\cdot)$ is defined by counting holomorphic polygons whose $k$-th marked points (the output points) lie in the stable submanifold of $e$ (or $p$).

Here are the situation that only constant polygons can contribute.
When we consider $m_k(\cdots,e,\cdots)$ where $e$ is put in the $i$-th place, we count discs which intersect the unstable submanifold $W^u (e)$ at the $i$-th marked point. Since $W^u (e)$ is $L_i \setminus \{e \}$, the corresponding moduli space is discrete only when the disc is constant. 
Similarly, concerning the coefficient of $p$ of $m_k(\cdots)$, we consider discs which intersect $W^s (p)=L_i \setminus \{ p \}$ and hence only constant discs may contribute.
In this setting, we have
$$m_2 (x,e) = x = (-1)^{|x|} m_2 (e,x)$$
which implies that $e$ serves as the unit of $L_i$ in the $A_\infty$-category.

\subsection{Fukaya category for orbi-surfaces}\label{sec:orbisurface}
The above construction of Seidel on the Fukaya category of surfaces can be
extended to orbifold surfaces (orbi-surface for short) following 
Section 9 of \cite{Se}, which we now explain briefly.

For simplicity we restrict to the case of a compact orbi-surface $\chi$ obtained 
as a global quotient of a compact surface $(\Sigma,\omega)$ by an effective action of a finite group $G$.
%(we will see later that  $\mathbb{P}^1_{a,b,c}$ can be defined in this way).
We denote the projection map by $\pi:\Sigma \to \chi$.

Embedded curves in $\chi$ which avoid orbifold points may be regarded as objects of
the Fukaya category. In fact, we also allow our objects to be immersed Lagrangian $\bar{L} \subset \chi$ which does not pass through orbifold points, whose pre-image $\pi^{-1}(\bar{L})$ is a union of embedded curves in $\Sigma$ (as a union of
$G$-action images).  The condition to avoid orbifold points can be achieved using an Hamiltonian isotopy in $\Sigma$.
%We expect that the derived Fukaya category do not depend on the choice of a cover, but we do not know how to prove it. 
Note that such a definition depends on the choice of a cover $\Sigma$.  For the orbifold $\bP^1_{a,b,c}$ we consider in this paper, we will fix a smooth cover of $\bP^1_{a,b,c}$ in which the Seidel Lagrangian (see Figure \ref{Seidel_Lag}) lifts as a submanifold.  Since we prove homological mirror symmetry for $\bP^1_{a,b,c}$, the derived Fukaya category is indeed independent of the choice of such a cover (Corollary \ref{cor:indep-cover}).

Hence we take the objects of $\mathcal{F}uk(\chi)$ to be 
$G$-equivariant family of Lagrangian objects in $\Sigma$, given by smooth connected non-contractible curves $L_1,\cdots, L_k$ in $\Sigma$ such that $\bigcup_{i=1}^k L_i$ is a $G$-orbit of $L_1$.
If there is no $g \in G$ such that $g(L_1)=L_1$, then we have $k=|G|$.
We denote by $G_{L_1}$ the subgroup of elements $g\in G$, with $g(L_1)=L_1$.
We assume that $G_{L_1}$ preserves an orientation of $L_1$ (otherwise, it is a non-orientable object).

%
%Fukaya category of $\chi$, $\mathcal{F}uk(\chi)$, can be defined following the construction for the surface case(\cite{Se}, \cite{Sh}) by considering  $G$-equivariant Lagrangian objects in $\Sigma$. Namely, an object of $\mathcal{F}uk(\chi)$  is given by a  $G$-equivariant family of Lagrangian submanifolds in $\Sigma$:
%consider . Such a collection $\{L_i\}$ is called a $G$-equivariant family of Lagrangian submanifolds.

We remark that $G$-equivariant family of Lagrangians in $\Sigma$ which passes through an orbifold point has the structure of an orbifold embedding, which we refer readers to \cite{CHS} for more details.
But we do not consider them in this paper.

Let us now explain the construction of  $\mathcal{F}uk(\chi)$ analogous to those in Seidel \cite{Se} and Sheridan \cite{Sheridan11}.
First, by choosing the curves transversal to each other, we may set the
Floer datum between two different Lagrangians to be trivial, and Floer datum between
the same Lagrangian $L$ to be  a Morse function on $L \cong S^1$.
One have to consider the perturbation datum (modulus and domain dependent) which is
used to perturbed the holomorphic pearly tree equation. We require that the space of Hamiltonian
functions $\mathcal{H}$ are smooth functions of $\chi$ (whose local lifts are smooth,
and associated Hamiltonian vector fields at orbifold points vanish),  almost
complex structures $J$ to be the standard complex structure$J_0$ near orbifold points of $\chi$.

In this paper, we will not consider orbifold holomorphic polygons, or the case that the domain
has  orbifold singularities (such polygons will lead to bulk deformations by twisted sectors
which will be considered in another paper). Hence the domain of holomorphic  polygons and pearly trees are the same as in manifold cases (we refer readers to \cite{CR} for more
details on (orbifold) holomorphic curves on orbifolds).
Note also that polygons without orbifold singularity can pass through orbifold points if it locally
 maps to the uniformizing cover of an orbifold point.

Therefore, the moduli space of perturbed holomorphic polygons  and pearly trees can be defined in the same way as in the manifold cases, and as explained in  \cite{Se}, a generic perturbation datum may be used to achieve transversality of the relevant moduli spaces which contribute to $A_\infty$-operations on the Fukaya category $\mathcal{F}uk(\chi)$

In fact, one can pull-back the Floer datum and perturbation datum of $\mathcal{F}uk(\chi)$
to construct the Fukaya category $\mathcal{F}uk(\Sigma)$. Note that there is an one-to-one correspondence between holomorphic polygons(and pearly trees) bounded by a $G$-equivariant family of Lagrangians and those in downstairs bounded by the corresponding immersed Lagrangian. This construction of $\mathcal{F}uk(\Sigma)$ implies that there exists a strict $G$-action on this $\AI$-category, whose
$G$-invariant part is the orbifold Fukaya category  $\mathcal{F}uk(\chi)$.
Namely, the morphisms in $\mathcal{F}uk(\chi)$ are the $G$-invariant part of the corresponding total morphisms between their $G$-equivariant families.

%
%Now, $\mathcal{F}uk(\chi)$ can be defined as a kind of a $G$-invariant of $\mathcal{F}uk(\Sigma)$, where  we enlarge the objects of $\mathcal{F}uk(\Sigma)$ by including their $g$-action images.
%As we can assume that the objects of $\mathcal{F}uk(\Sigma)$ intersect transversely to each other, we may set the Floer datum between two different objects $L$ and $L'$ of $\mathcal{F}uk(\Sigma)$  to 
%be trivial, and set the Floer datum between  $L$ and itself to be a $G_{L}$-invariant Morse-function $f:L \to \R$.
%In this way, we can make the whole Floer data to be $G$-equivariant (this step may not be possible in general dimensions).

%Following the construction of Seidel \cite{Se} and Sheridan \cite{Sh}, we can equip $\mathcal{F}uk(\Sigma)$
%with a strict $G$-action because Floer data are $G$-equivariant. Thus, we can define the morphisms in
%$\mathcal{F}uk(\chi)$ to be the $G$-invariant part of the corresponding total morphisms between their $G$-equivariant families.

\section{Geometric construction of localized mirror functor} \label{sect:geom_fctor}
In this section we apply the algebraic construction in Section \ref{GMir} to the setting of generalized SYZ in Section \ref{gen_SYZ} and obtain a functor for the purpose of homological mirror symmetry.  Moreover we exhibit locality of the functor by working out the simplest possible example, namely the complex projective line $\bP^1$.

Recall that in Section \ref{gen_SYZ} we fix an immersed Lagrangian $L$ and use its local deformations (see Setup \ref{setup2}) to construct a localized mirror $(V,W)$.  Applying Definition \ref{def:loc_mirr}, Theorem \ref{thm_mir_fctor} to the filtered $\AI$-category $\mathcal{C} = \mathcal{F}uk\,_\lambda (X)$ and the object $A = \bar{L} \in \mathcal{C}$, we obtain a functor $\mathcal{LM}^\BL$:

%satisfying Setup \ref{setup2}.
%We consider the weak bounding cochain $b=\sum_j x_j b_j$, where $x_j$'s are
%regarded as variables over the Novikov field $\Lambda_0$. We write $\BL=(L, b)$ for simplicity.
%We may regard $T^{\omega(\beta)}$ as  K\"{a}hler parameters with values in $\C$.
%\begin{defn}
%We define the {\em localized mirror of the immersed Lagrangian $L$ } to be 
%the Landau-Ginzburg model
%$$W: V \to \C$$ 
%\end{defn}
%This is called a  localized mirror as it is constructed from the infinitesimal deformation space
%of a given immersed Lagrangian. But the following localized mirror functor is global in the
%sense that we have a map from the entire Fukaya category (not just from the
%tubular neighborhood of the immersed Lagrangian) to the dg category of matrix factorization of $W$.

%As explained in the introduction, one major difference from the original SYZ is that we need to take the `dual' of a Lagrangian which is no longer a torus.  The main idea is that the set of complexified deformations gives such a `dual'.  For a torus $T$, the dual torus $T^*$ is given by
%$$ T^* = \{\nabla: \nabla \textrm{ is a flat } U(1) \textrm{ connection on } T \}.$$
%On the other hand, the complexified Lagrangian deformations of $T$ are parametrized by
%$$H^1(T,\real) \oplus \consti (H^1(T,\real) / H^1(T,\integer))$$
%which contains $T^*$ as the imaginary part.  Thus for a general Lagrangian $L$, it is natural to define its `dual' to be the complexified infinitesimal Lagrangian deformations of $L$.

\begin{theorem}[Geometric version of Theorem \ref{thm_mir_fctor}] \label{thm:locmirfun}
Let $X$ be a symplectic manifold and $(V,W)$ be its localized mirror as in Setup \ref{setup2}.  We have an $\AI$-functor
\begin{equation}\label{eq:LMAIfunc}
\mathcal{LM}^\BL:  \mathcal{F}uk\,_\lambda (X) \to \mathcal{MF}(W(b)- \lambda)
\end{equation}
where $\mathcal{F}uk\,_\lambda (X)$ is an $\AI$-category and $\mathcal{MF}(W(b)- \lambda)$ is a dg category.
\end{theorem}

Let us now consider the derived version of the functor \eqref{eq:LMAIfunc}. Any $A_\infty$-functor $\mathcal{C} \to \mathcal{D}$ admits a canonical extension $D^\pi (\mathcal{C}) \to D^\pi (\mathcal{D})$ where $D^\pi$ means taking the split-closure of a derived category (see the proof of Lemma 2.4 of \cite{Se2}). Thus, the $A_\infty$-functor $\mathcal{LM}^\BL$ induces a functor of triangulated categories 
$$D^\pi \mathcal{LM}^\BL : D^{\pi} ({\Fuk}\,_\lambda (X) ) \to  D^{\pi} ({\mathcal{MF}} (W - \lambda)).$$

Homological mirror symmetry occurs as special instances where $D^\pi \mathcal{LM}^\BL$ is an equivalence.  We have the following theorem which basically compares generators on the two sides using our functor.  The idea of matching split-generators and their morphism spaces on the two sides was exploited in literatures such as \cite{Se} and \cite{Sh} to prove homological mirror symmetry.  Theorem \ref{thm:criterion_equiv} and its proof employ such an idea.  The upshot here is that we have a geometric construction of the mirror $W$, and the functor is naturally constructed from the geometric setup (rather than matching the split generators by hand).

\begin{theorem}\label{thm:criterion_equiv}
%Suppose that the immersed Lagrangian $\BL$ split-generates 
Suppose that there exists a set of Lagrangians $\{L_i: i \in I\}$ which split-generates 
$D^\pi \Fuk_0 (X)$, and suppose the functor $\functor$ induces an isomorphism on cohomologies
$HF (L_i,L_j) \overset{\cong}{\rightarrow} \Mor (\functor(L_i),\functor(L_j))$
for all $i,j \in I$.
Then the derived functor
$$D^{\pi} (\functor): D^{\pi} \Fuk_0 (X) \to D^{\pi} \MF (W) $$
is fully faithful.  Furthermore if $\{\mathcal{LM}^\BL(L_i): i \in I\}$ split-generates $D^\pi \MF (W)$, then $D^{\pi} (\functor)$ is a quasi-equivalence.
\end{theorem}

\begin{proof}
The proof employs standard techniques in homological algebra.  By assumption any object in $D^\pi \Fuk_0 (X)$ can be obtained from taking direct sums and cones of $\{L_i: i \in I\}$ in finite steps.  First we show that the derived functor is fully faithful.  For any two objects $A, B$ of $D^\pi \Fuk_0 (X)$, we consider the induced map $h: HF (A,B) \to \Mor (\functor(A),\functor(B))$ from the functor.  We already know that this map is an isomorphism when $A,B$ are taken from the set $\{L_i: i \in I\}$.  Now if $A = A_1 \oplus A_2$ and $HF (A_i,B) \to \Mor (\functor(A_i),\functor(B))$ for $i=1,2$ are isomorphisms, then $h$ is just the direct sum and hence is also an isomorphism.  Similarly it is true for $B = B_1 \oplus B_2$.

Then we consider the case that $HF (A,B_i) \to \Mor (\functor(A),\functor(B_i))$ for $i=1,2$ are isomorphisms, and $B$ is obtained as a cone of $B_1$ and $B_2$, namely we have the exact triangle $B_1 \to B_2 \to B \overset{[1]}{\rightarrow}$.  The derived functor $D^\pi \mathcal{LM}^\BL$ sends exact triangles to exact triangles, and hence we have the exact triangle $D^\pi \mathcal{LM}^\BL(B_1) \to D^\pi \mathcal{LM}^\BL(B_2) \to D^\pi \mathcal{LM}^\BL(B) \overset{[1]}{\rightarrow}$.  By applying five lemma to $D^\pi \mathcal{LM}^\BL$ acting the exact sequence
$$ \ldots \to HF^k(A,B_1) \to HF^k(A,B_2) \to HF^k(A,B) \to HF^{k+1}(A,B_1) \to HF^{k+1} (A,B_2) \to \ldots $$
we conclude that $h: HF (A,B) \to \Mor (\functor(A),\functor(B))$ is also an isomorphism.  The proof for the case that $A$ is a cone of $A_1$ and $A_2$ is similar.  This proves that the functor is fully faithful.

If further $\{\mathcal{LM}^\BL(L_i): i \in I\}$ split-generates $D^\pi \MF (W)$, then the functor sends generators to generators, and induces isomorphisms on the morphism spaces between each pair of generators.  Thus the two derived categories are equivalent.
\end{proof}

%Idea similar to this theorem was used in several literatures (for instance, \cite{Se} and \cite{Sh}) to prove homological mirror symmetry.

Thus we can interpret homological mirror symmetry as a special case where our functor satisfies the conditions in Theorem \ref{thm:criterion_equiv}.  We will study concrete examples, namely the orbifold $\proj^1_{a,b,c}$ in Section \ref{ex333} and \ref{sect_dim1}.

\subsection{Geometric transform of Lagrangian branes} \label{geom_transf}
Consider the localized mirror functor $\mathcal{LM}^\BL$ in this geometric setup.
On the level of objects, it sends Lagrangian submanifolds or immersions to matrix factorizations. 
More precisely, we need to consider the Fukaya category $\mathcal{F}uk_\lambda$
consists of objects which are weakly unobstructed ($\Z/2$ graded) spin Lagrangian submanifolds with potential $\lambda \in \Lambda_0$.

Such a category arises naturally since Floer cohomology between two weakly unobstructed Lagrangian submanifolds are defined only if their potentials are the same (\cite{FOOO}).
Namely, for two Lagrangian submanifolds $L_i$ for $i=0,1$ with bounding cochains $b_i$ and associated potentials $\lambda_i$ (for $i=0,1$),
the differential $m_1$ for $CF(L_0,L_1)$ satisfies the equation
%\begin{equation}\label{eq:flm1}
$m_1^2 = \lambda_1 - \lambda_0.$
%\end{equation}
Hence $m_1^2=0$ if $\lambda_1 = \lambda_0 = \lambda$.  Now taking $(L_0,b_0) = \BL$ which is the reference Lagrangian immersion, the equation becomes
%\begin{equation}\label{eq:flm2}
$m_1^2 = W - \lambda_0.$
%\end{equation}

On the other hand, recall from Definition \ref{def:mf} that a matrix factorization $(P, d)$ for a critical value $\lambda$ of $W$ satisfies
$d^2 = (W -\lambda) Id.$
Thus Lagrangian Floer theory naturally gives rise to matrix factorizations.  Our localized mirror functor is a categorical formulation of this observation.

Given a weakly unobstructed ($\Z/2$ graded) spin Lagrangian immersion $L'$ with potential $\lambda \in \Lambda_0$ and a weak-bounding cochain $b'$, we associate to it a matrix factorization $(P=P_0 \oplus P_1)$ which is the Floer complex $CF(\BL, L')$ with differential $d = m_1^{b,b'}$ (where $P_0$ and $P_1$ are the even and odd parts of the complex respectively).  Recall that $b$ varies in $V$, which probes the mirror matrix factorization.  For an unobstructed Lagrangian $L_i$ with $b_i=0$ and $\lambda_i=0$,  the corresponding matrix factorization under the localized mirror functor $\mathcal{LM}^\BL$ is given by 
the Floer complex $(CF(\BL,L_i), m_1^{b,0})$.

On the level of morphisms, an element $\alpha \in CF(L_1,L_2)$ provides a morphism between two
matrix factorizations $(CF(\BL,L_1), m_1^{b,0})$, $(CF(\BL,L_2), m_1^{b,0})$
which is induced from the product $m_2: CF(\BL, L_1) \times CF(L_1,L_2) \mapsto CF(\BL, L_2)$.
See Section \ref{sect_alg_fctor} for the details.

The similarity of the above two equations was first used by Chan and Leung in \cite{CL} to compute mirror matrix factorizations of Lagrangians in $\CP^1$, and was further explored in the work of Cho-Hong-Lee \cite{CHL} for the case of $\CP^1 \times \CP^1$ and weighted projective lines $\CP^1_{n,m}$.  In this paper we construct the whole functor in full generality rather than just on the object levels.
\subsection{An example: the projective line}\label{sec:cp1}
In this section we exhibit locality of our functor by studying the example $X = \bP^1$.  We will construct two localized mirrors $W^{\pm}$ and functors $\LM^{\BL^{\pm}}:\Fuk(X) \to \MF(W^\pm)$ using two different Lagrangian branes $\BL^+$ and $\BL^-$.  $W^{\pm}$ will be \emph{different} from the Hori-Vafa mirror of $\bP^1$, and they can be regarded as the Hori-Vafa mirror localized at each of the two critical points.  Moreover we will show that each of the functors $\LM^{\BL^{\pm}}$ only sees a part of the Fukaya category.

In this sense the mirrors constructed in our approach are more general and flexible than the mirrors provided by the string theorists.  Our mirror functor in general reflects a subcategory of the Fukaya category and comes in a much more direct way.  Our forthcoming paper will develop a mirror construction using more than one Lagrangian immersions.  This will capture more information about the global symplectic geometry.

%It is well-known that the Fukaya category of $\bP^1$ has two objects, one is an equator $L$ with flat line bundle of trivial holonomy, the other is also $L$, but with flat line bundle of holonomy $(-1)$.

$\BL^{\pm}$ is taken to be the union of two
equators intersecting with each other transversely endowed with flat $U(1)$ connections.  Identify $\bP^1$ as the unit sphere of $\R^3$.
 Let $L_V$ be the vertical equator (with the standard orientation) and $r(L_V)$
the rotation of $L_V$ by $\pi/2$ about $z$-axis.  $L$ is taken to be the immersion $L_V \cup r(L_V) \subset \bP^1$.
Both $L_V$ and $r(L_V)$ are equipped with the trivial spin structure.  Then $\BL^{\pm}$ are defined as the immersion $L$ equipped with the flat $U(1)$ connections $\nabla^+ := 1 \cup 1$ and $\nabla^-:=-1 \cup -1$ respectively.  $\nabla^+ := 1 \cup 1$ denotes the trivial connection, and $\nabla^-:=-1 \cup -1$ denotes the flat connection with holonomy $-1$ on both circles $L_V$ and $r(L_V)$.

 $L$ has two self intersection points, namely the north and south poles.  Each self-intersection point correspond to two immersed generators, one has degree odd and one has degree even.  Denote by $X, Y$ the degree odd generators corresponding to north and south poles respectively.  Let $\bar{X}, \bar{Y}$ be the corresponding degree even generators.  Now we consider the formal deformations $b = xX + yY$ for $x, y \in \C$.  First of all we need to solve the Maurer-Cartan equation for weak bounding cocycles.

%We will show that $b:= xX + yY$ are weak Maurer-Cartan solutions.  We denote the collection $(L, \nabla^\pm, b)$ by  $\mathbb{L}^{\pm}$.  We will compute the associated potential $W^{\pm}$, and consider  localized mirror functors $\mathcal{LM}^{\mathbb{L}^\pm}: \mathcal{F}uk (S^1) \to \mathcal{MF}(W^\pm)$.  Namely, we will find the corresponding matrix factorization of a great circle $L_0$ with holonomy $\pm 1$ under $\mathcal{LM}^{\mathbb{L}^\pm}$.
  
\begin{figure}[h]
\begin{center}
\includegraphics[height=2in]{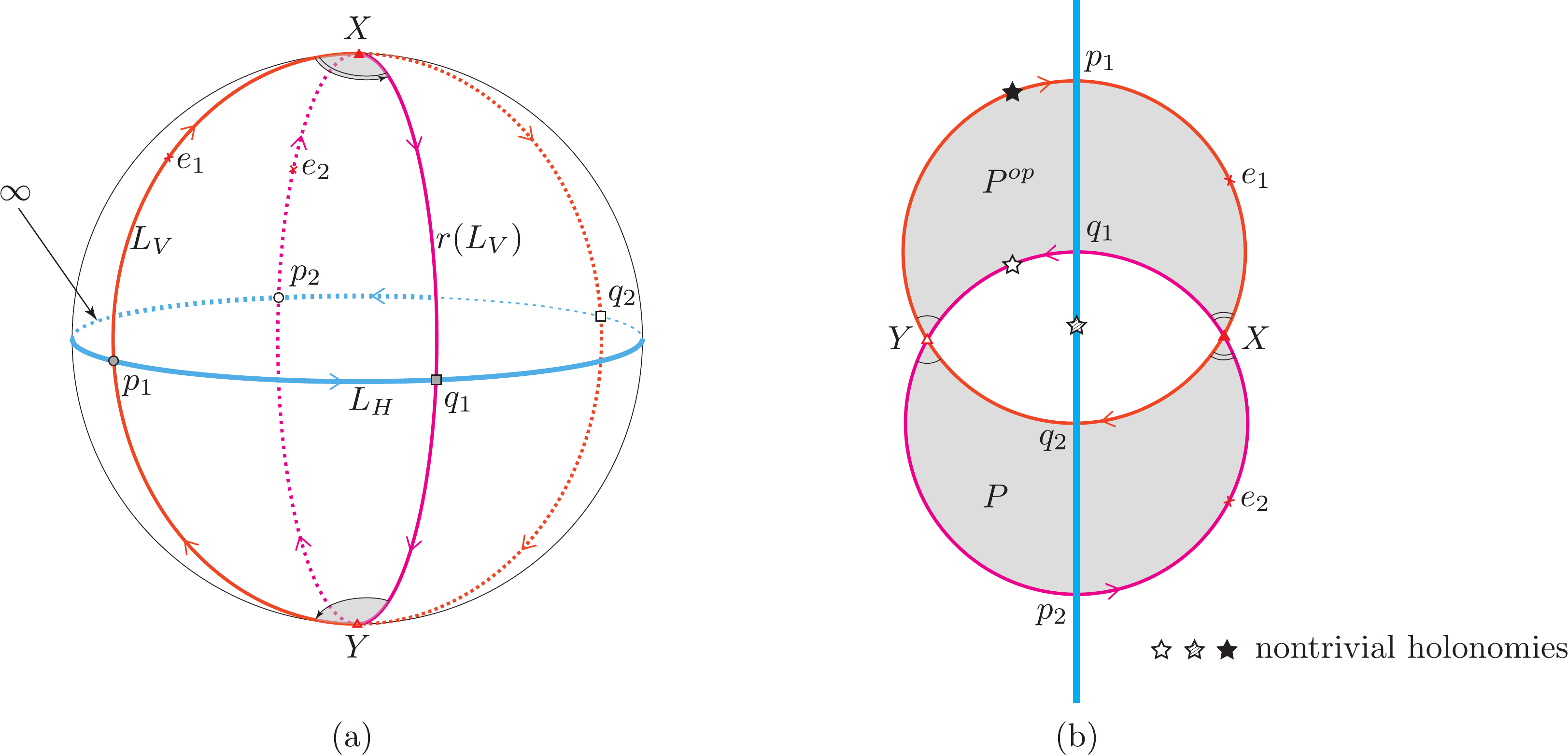}
\caption{(a) $\mathbb{L}$ and $L_H$ (b) holomorphic discs contributing to the potential of $\mathbb{L}$}\label{twoeqCP1}
\end{center}
\end{figure}

\begin{lemma} \label{lem:P1-MC}
For both $\BL^+$ and $\BL^-$, $b:= xX + yY$ is a weak bounding cocycle for all $x, y \in \C$. 
\end{lemma}
\begin{proof}
We will use the Bott-Morse model of Seidel \cite{Se}.  Take the standard height function $f$ on $L_V$ and $r^\ast f$ on $r(L_V)$.  Denote the minima to be $e_1$ on $L_V$ and $e_2$ on $r (L_V)$.  Then $e= e_1 \oplus e_2$ serves as a unit of the $A_\infty$-algebra $CF(\mathbb{L}^{\pm}, \mathbb{L}^{\pm})$.

First consider $m_0$, which is contributed from the holomorphic discs bounded by $L_V$ and that bounded by $r(L_V)$.  Since they bound the same holomorphic discs, $m_0$ for both $L_V$ and $r(L_V)$ are proportional to the unit, and $m_0$ for $\BL^\pm$ is also proportional to the unit.  We will show that $m_{k}(x_1,\cdots, x_k) =0$ for $k=1$ or $k \geq 3$ and $x_i$ to be either $X$ or $Y$ for every $i=1,\cdots, k$.  Moreover $m_2(X,X) = m_2(Y,Y)=0$, and $m_2(X,Y) + m_2(Y,X)$ is proportional to the unit $e$.  Putting all these together, we obtain that $b= xX + yY$ is a weak bounding cocycle.

Second $m_{k}(x_1,\cdots, x_k) \not= 0 $ for $x_i$ to be either $X$ or $Y$ only when $\mu(\beta)=2$.  There are exactly two holomorphic polygons of Maslov index two, namely $P$ and $P^{op}$ shown in Figure \ref{twoeqCP1}.  They do not contribute to $m_2(X,X)$, $m_2(Y,Y)$ nor $m_k$ for $k \geq 3$.  Hence $m_2(X,X) = m_2(Y,Y)=0$, and $m_k(x_1,\cdots, x_k) =0$ for $k \geq 3$.  $P$ and $P^{op}$ contribute to $m_2(X,Y) + m_2(Y,X)$, and it follows from direct computation that $m_2(X,Y) + m_2(Y,X)$ is proportional to $e$.

We now show that the contributions from $P$ and $P^{op}$ in $m_1 (X)$ cancel each other and hence $m_1(X)=0$.  Similarly $m_2(Y)=0$ for the same reason.  The computations for $\mathbb{L}^+$ and $\mathbb{L}^-$ are similar and so we will only do it for $\mathbb{L}^+$.  The contribution of $P$ to $m_1 (X)$ is $q^{\omega(P)} \bar{Y}$, since the orientation induced by $P$ on $\partial P$ agrees with the orientation of the Lagrangian.  On the other hand $P^{op}$ contributes by $(-1)^{|X| + |\bar{Y}|} q^{\omega(P^{op})} \bar{Y} = - q^{\omega(P)} \bar{Y}$.  Thus $m_1(X) = q^{\omega(P)} \bar{Y} - q^{\omega(P)} \bar{Y} = 0$.
\end{proof}
% For simplicity, we write by $\mathbb{L}^+$ the collection $(L, \nabla^+, b = xX + yY)$
%  We regard $\mathbb{L}$ as an immersed Lagrangian in $\C P^1$ and compute the local mirror functors induced by $(\mathbb{L}, \nabla^{\pm})$. For simplicity, we write $\mathbb{L}^+$ for $(\mathbb{L}, \nabla^{+})$ and $\mathbb{L}^-$ for $(\mathbb{L}, \nabla^{-})$
%
%$\mathbb{L}$ has two self-intersection points which give rise to four immersed generators of $CF(\mathbb{L}^{\pm}, \mathbb{L}^{\pm})$ denoted by $X$, $\bar{X}$, $Y$, $\bar{Y}$.
%We firstly show thatfor any $x,y \in \C$.
%
%We set $V = \C^2$  and compute the superpotential $ W^{\pm} : V \to \C$.

Now we compute the superpotential $W^\pm: \C^2 \to \C$ by counting holomorphic polygons bounded by $\BL^\pm$ respectively.

\begin{prop}
The superpotential constructed from each of $\bL^\pm$ is $W^{\pm}= q^2 xy \pm q^4$ respectively, where $q^8 = \exp \left( - \int_{\bP^1} \omega \right)$ is the K\"ahler parameter of $\bP^1$.
\end{prop}
\begin{proof}
From the proof of Lemma \ref{lem:P1-MC}, weh have $m^{\mathbb{L}^{\pm}} ( \conste^b) = W^{\pm}(b) e$ for $b=xY+yY$ and $(x,y) \in \C^2$.
As higher order terms vanish and $m_1=0$, we have
$m ( \conste^b) = m_0  + m_2 (b,b).$ 
The term $m_0$ is given by Maslov index two holomorphic discs bounding $L_V$ and $r (L_V)$, which is $m_0 = \pm q^4 e_1 \pm q^4 e_2 = \pm q^4 e$ for $\mathbb{L}^\pm$. $P$ and $P^{op}$ contribute to $m_2^{\mathbb{L}^+} (b, b)$ by $q^2 xy e_1$ and $q^2 xy e_2$. So, $m_2 (b, b) = q^2 xy e$. 
Therefore, $W^+ = q^2 xy + q^4$. Similarly for $\mathbb{L}^-$, two nontrivial holonomies cancel each other along $\partial P^{op}$ to give the same $xy$ term in $W^-$.  Thus $W^- = q^2 xy - q^4$.
\end{proof}

Next we consider our functor $\LM^{\BL^\pm}$ and use it to transform the two split-generators of the derived Fukaya category, namely $L_H^+$ which is the horizontal equator equipped with the trivial flat $U(1)$ connection, and $L_H^-$ which is the horizontal equator equipped with the flat $U(1)$ connection with holonomy $(-1)$ along the equator.  It is well-known that the potential value of $L_H^\pm$ is $\pm q^4$ (see \cite{C}).

There are four intersection points between $\mathbb{L}$ and $L_H$ : $p_1$ and $p_2$ of odd degree, $q_1$ and $q_2$ of even degree, see Figure \ref{twoeqCP1}.  By a direct count of holomorphic strips bounded by $\BL^\pm$ and $L_H$, we obtain the following proposition, which implies that the functor $\LM^{\BL^\pm}$ reflects the subcategory generated by $L_H^\pm$.

\begin{prop}
$\LM^{\BL^+}(L_H^+)$ is the matrix factorization
\begin{equation}\label{PLH}
H = \,\,\,
\bordermatrix{ & p_1 & p_2 \cr
q_1 & -xq & 0 \cr
q_2 & 0 & -yq \cr              } \qquad
K = \,\,\, \bordermatrix{ & q_1 & q_2 \cr
p_1 & -yq & 0 \cr
p_2 & 0 & -xq \cr              }
\end{equation}
of $W^+ - q^4$ (namely $HK = KH = (W^+ - q^4) \Id$), and $\LM^{\BL^+}(L_H^-)$ is a matrix factorization
%$$ $$
of $W^+ + q^4$, which is a trivial object since $-q^4$ is a regular value of $W^+$.
Similarly $\LM^{\BL^-}(L_H^+)$ is a trivial object, while $\LM^{\BL^-}(L_H^-)$ takes the same expression \eqref{PLH} which serves as a matrix factorization of $W^- + q^4$ (which equals to $W^+ - q^4$).
\end{prop}

\section{Generalized SYZ for a finite-group quotient} \label{sect:quot}
In the early development of mirror symmetry, physicists employ finite group actions to construct mirrors.  The symmetry provided by a finite group simplifies the geometry and allows one to study a manifold via its quotient. 

In this section, we take a Lagrangian immersion into a global quotient orbifold $X=[\tilde{X}/G]$ and carry out the generalized SYZ mirror construction.  We will consider $G$-equivariant Lagrangian immersions in $\tilde{X}$
(satisfying Assumption \ref{assump1}), which provides a Lagrangian immersion to $X$.
Floer theory of the Lagrangian immersion in $X$ can be obtained by taking the
$G$-invariant part of the Floer theory of $G$-equivariant Lagrangian immersion in $\tilde{X}$.

We will apply our construction of localized mirror functor in the previous sections \emph{$G$-equivariantly} on $\tilde{X}$.  The purpose is two-fold: On one hand we study mirror symmetry for the symplectic orbifold $X$; on the other hand, we obtain the mirror of $\tilde{X}$ as a quotient of the mirror of $X$ by the dual group of $G$.

\subsection{Generalized SYZ construction for a global quotient}

Let $\tilde{X}$ be a K\"ahler manifold equipped with an effective action by a finite group $G = \{g_1, \ldots, g_{|G|}\}$.  Then $X = \tilde{X}/G$ is a K\"ahler orbifold.  Denote the canonical quotient map by $\pi: \tilde{X} \to X$.

\begin{assumption}\label{assump1}
Let $L$ be a compact oriented connected spin Lagrangian submanifold of $\tilde{X}$.
We assume that for each $g \in G$ with $g \neq 1$, the $g$-action image of $L$, $g \cdot L$, intersects 
transversely with $L$. 
\end{assumption}

 Then the image $\bar{L} = \pi (L) \subset X$ is a Lagrangian immersion whose self-intersections are transverse.  We denote the normalization map $\iota := \pi|_L: L \to \bar{L}$ and $H^*(L \times_\iota L)$ gives the deformation space of $\bar{L}$ (before taking obstructions into account).

With Assumption \ref{assump1}, it is relatively easy to define the immersed Lagrangian Floer theory of
$\bar{L}$. We consider the Floer theory of  the family of Lagrangians $\pi^{-1}(\bar{L}) = \bigcup_{g \in G} g \cdot L \subset \tilde{X}$, which is also an immersed Lagrangian with transverse self-intersections, and take a $G$-invariant part. As $G$ acts freely on the set of objects $\{gL \mid g \in G\}$, it is easy to make the associated $\AI$-structure to have a strict $G$-action.
 Without the assumption \ref{assump1}, the construction becomes much
more involving and we refer readers to \cite{CH} for more details.

%Since $X$ is a global quotient, the immersed Lagrangian Floer theory of $\bar{L} \subset X$ is reflected from ($G$-invariant part of) that of $\pi^{-1}(\bar{L}) = \bigcup_{g \in G} g \cdot L \subset \tilde{X}$, which is also an immersed Lagrangian with transverse self-intersections.  

Write $\coprod_g L_g = \coprod_{g \in G} g \cdot L$, which gives a normalization $\tilde{\iota}: \coprod_g L_g \to \pi^{-1}(\bar{L})$
of $\pi^{-1}(\bar{L})$.  Let $\tilde{R}$ be the set of components of $(\coprod L_g) \times_{\tilde{\iota}} (\coprod L_g)$.

\begin{prop} \label{prop:up_down}
$G$ acts freely on $\tilde{R}$, and the quotient has a natural identification with $R$, the set of components of $L \times_\iota L$. $C^*(L \times_\iota L)$ (the space of Floer cochains downstairs) has a natural identification with $(C^*((\coprod L_g) \times_{\tilde{\iota}} (\coprod L_g)))^G$ (the space of $G$-invariant Floer cochains upstairs).
\end{prop}

\begin{proof}
$(\coprod L_g) \times_{\tilde{\iota}} (\coprod L_g)$ consists of $|G|$-copies of $L$ and ordered pairs $(p_+,p_-)$, where $p_+,p_-$ belong to $L_g, L_h$ for $g \not= h$ respectively and $\tilde{\iota} (p_+) = \tilde{\iota} (p_-)$.  In other words,
$$ (\coprod L_g) \times_{\tilde{\iota}} (\coprod L_g) = \coprod_{g \in G} \left( L_g \amalg \left(\coprod_{h \not= 1} L_g \times_{\tilde{\iota}} L_{gh}  \right) \right). $$
$a \in G$ acts by sending
$ L_g \amalg \left(\coprod_{h \not= 1} L_g \times_{\tilde{\iota}} L_{gh}  \right) \mapsto L_{ag} \amalg \left(\coprod_{h \not= 1} L_{ag} \times_{\tilde{\iota}} L_{agh}  \right)$
and hence it induces a free action on the set of components $\tilde{R}$.  The quotient map induces a one-one correspondence between
$ L \amalg \left(\coprod_{h \not= 1} L \times_{\tilde{\iota}} L_{h}  \right) $
and $L \times_\iota L$.
\end{proof} 

\begin{lemma} \label{g_X}
Each immersed sector $X \in R - \{0\}$ of $\bar{L}$ is canonically associated with an element $g_X$ in the finite group $G$.
\end{lemma}

\begin{proof}
$X = (p_-,p_+)$ can be identified with an immersed sector $\tilde{X} = (\tilde{p}_-,\tilde{p}_+)$ of $L\subset \tilde{X}$ upstairs.  Then $\tilde{p}_- \in L_{g_-}$ and $\tilde{p}_+ \in L_{g_+}$ for certain $g_-, g_+ \in G$.  $g_-, g_+ \in G$ themselves depend on the choice of lifting, but the quotient $g_X = g_-^{-1} g_+$ does not.  Thus $X$ is associated $g_X \in G$.
\end{proof}

Explicitly let $X_1, \ldots, X_{|R|-1} \in R - \{0\}$ be all the immersed sectors of $\bar{L}$.  $|R|-1$ is an even integer since the elements in $R - \{0\}$ naturally form pairs: $\bar{X} = (p_+,p_-)$ if $X = (p_-,p_+)$.  For each ${X_i} \in R - \{0\}$ we can choose a representative in $\tilde{R}-\{0\}$ upstairs, which is also denoted by $ X_i$ by abuse of notations.  Upstairs $ X_i$ is of the form $(p_g, p_h)$ for $p_g \in  g \cdot L$ and $p_h \in  h \cdot L$ such that $\tilde{\iota}(p_g) = \tilde{\iota}(p_h)$.  Then $ \tilde{R} = G \amalg \{g \cdot X_i: g \in G, i =1, \ldots, |R|-1 \}. $

%\begin{figure}[htp]
%\begin{center}
%\includegraphics[scale=0.5]{angles}
%\caption{Angles indicating $(a)$ $(X_g, X_h)$ and $(b)$ $(X_h, X_g)$}\label{angles}
%\end{center}
%\end{figure}
%
%

We may also choose cycles $\Phi_j$'s which forms a basis of $H^*(L)$.  Then the deformation space $H \subset C^*(L \times_\iota L)$ is defined as the span of $ X_1, \ldots,  X_{|R|-1}$ and $\Phi_1, \ldots, \Phi_{h^*(L)}$, and we denote the corresponding coordinates as $x_1,\ldots,x_{|R|-1},z_1, \ldots, z_{h^*(L)}$.

Upstairs we consider the span of $\left(\sum_g g\cdot\Phi_{j} \right)$ and $\left( \sum_g g \cdot X_i \right)$, which is a subspace of $(C^*((\coprod L_g) \times_{\tilde{\iota}} (\coprod L_g)))^G$.  Elements in the vector space $H^*((\coprod L_g) \times_{\tilde{\iota}} (\coprod L_g))$ are written uniquely as
$ \sum_{j,g} z_{j,g} [g \cdot \Phi_j] + \sum_{i,g} x_{i,g} [g \cdot X_i] $
and $G$-invariant elements are written as
$$ b = \sum_{j} z_{j} \left(\sum_g g \cdot \Phi_j \right) + \sum_{i} x_{i} \left( \sum_g g \cdot X_i \right) \in (H^*((\coprod L_g) \times_{\tilde{\iota}} (\coprod L_g)))^G $$
which are identified with the deformation element
$ \sum_{j=1}^{h^*(L)} z_{j} \Phi_{j} + \sum_{i=1}^{|R|-1} x_{i} X_i \in H^*(\bar{L} \times_\iota \bar{L}) $
downstairs.

\begin{prop} \label{prop:equiv mk}
For $\phi_1, \ldots, \phi_k \in (C^*((\coprod L_g) \times_{\tilde{\iota}} (\coprod L_g)))^G$,
$m^{(\tilde{X},L)}_k (\phi_1,\ldots,\phi_k)$ still belongs to $(C^*((\coprod L_g) \times_{\tilde{\iota}} (\coprod L_g)))^G$.
\end{prop}

\begin{proof}
$G$-action on $\tilde{X}$ induces a free $G$-action on the moduli space of stable discs $\bigcup_{\beta} \CM_{k+1} (\beta;\phi_1,\ldots,\phi_k)$ such that the evaluation map $\ev_0$ is $G$-equivariant.  Thus $m_k$, which is the image chain under $\ev_0$, is $G$-invariant.
\end{proof}

By the above proposition and the identification between cochains upstairs and downstairs in Proposition \ref{prop:up_down}, we can define
$ m_k^{X,\bar{L}}: C^*(L \times_\iota L)^{\otimes k} \to C^*(L \times_\iota L) $
to be $m_k^{(\tilde{X},L)}$ restricted to $(C^*((\coprod L_g) \times_{\tilde{\iota}} (\coprod L_g)))^G$.

The unit $1_L$ downstairs is identified with the unit $\sum_{g \in G} 1_{L_g}$ upstairs.  Then we consider weakly-unobstructed cocycles, which are elements $b \in H$ such that
$$m_0^b = \sum_{k \geq 0} m_k (b,\ldots,b) = W \sum_{g \in G} 1_{g \cdot L}$$
for some coefficient $W$.  Let $V$ be the span of $\{\Phi_{j}\}_{j=1}^{l}$ and $\{X_i\}_{i=1}^{N-l}$ for some $j$ and $l$, and assume that all elements in $V$ are weakly unobstructed cocycles.  Take $V$ to be the subspace of all these elements.  Then $(V,W)$ forms a Landau-Ginzburg model.  For simplicity of notations we take $l = 0$.  Write
$b =  \sum_{i=1}^{N} x_i X_i $,
then in terms of the coordinates $x_i$'s of $V$,
$$W = \sum_{(i_1, \ldots, i_N) \in \Z_{\geq 0}^N} \sum_{\beta \in H_2(\tilde{X},L)} \sum_{\substack{P\\{g_1,\ldots,g_k \in G}}} q^\beta n(\beta;\, g_1 \cdot X_{P(1)}, \ldots, g_k \cdot X_{P(k)}) x_1^{i_1} \ldots x_N^{i_N}$$
where $k = i_1 + \ldots + i_N$, and we are summing over all permutations $P$ which are maps $P: \{1,\ldots,k\} \to \{1, \ldots, N\}$ such that $|P^{-1} \{j\}| = i_j$ for all $j=1,\ldots,N$.

\subsection{Dual group action on the mirror and $\WH{G}$-invariance}\label{sec:dualgponmir}
We show that there exists a canonical action of the dual group $\WH{G}$ on the mirror, which leaves $W$ invariant.  Then $(V/\WH{G},W)$ will be taken to be a generalized SYZ mirror of $\tilde{X}$.

 We first recall the notion of the dual group (character group) $\WH{G}$
of $G$. A character of a finite group $G$ is a group homomorphism $\chi:G \to S^1=U(1)$,
where the group structure of $U(1)$ is the multiplication. The set of all characters of  $G$ form
a group $\WH{G}$, with the group law that two homomorphisms are multiplied
pointwise: for two characters, $\chi, \psi:G \to U(1)$, we have $(\chi \psi)(g) = \chi(g)\psi(g)$.
The trivial character $\chi_0$ is given by the constant map to $1$, and is the identity of 
the group $\WH{G}$. It is well-known that $\WH{G}$ is isomorphic to $G$, but not canonically.

We work in the setting of the previous section, with Assumption \ref{assump1}.
Fix a representation of $G$, denoted it as $\phi:G \to End(V)$.
Given an equivariant family $\sqcup_{g \in G} g \cdot L$, we consider the following trivial bundle $\CL_\phi = \sqcup_{g \in G} g \cdot L \times V.$
For its $G$-action, given $h \in G, x \in g \cdot L$, we define
$h \cdot (x, v) := ( hx, \phi(h)v),$
where $hx \in hg \cdot L$.  Note that this $G$-action is compatible with the projection $\pi_1$ to the first component.

The Floer complex with the above bundle data will lead us to the dual group $\WH{G}$ action in the mirror construction.
Consider  $\chi \in \WH{G}$ as an one-dimensional representation $\phi$.
When Lagrangians are equipped with a trivial bundle (with the trivial holonomy),
each intersection point $p_i \in L \cap g \cdot L$ (for the mirror variable $x_i$) corresponds to the generator of
$Hom(\CL_\chi|_{L, p_i}, \CL_\chi|_{g L, p_i}).$
Hence for the mirror variable $x_i$ corresponding  to $p_i$, it is natural to define an action of $\chi \in \WH{G}$ to be
$\chi \cdot x_i = \chi(g) x_i.$

One way to interpret the above action is from the comparison of the generators of
the trivial character $\chi_0$ and that of  $\chi$.
We define $\WH{G}$-action on the $z_j$ variables to be trivial.

\begin{prop}\label{identifypo}
The above $\WH{G}$-action leaves the superpotential $W$ invariant.
\end{prop}
\begin{proof}
Recall from Lemma \ref{g_X} that each immersed sector is associated with an element in $G$.  At each $i$-th turn of a disc corresponding to a group element $g_i\in G$, we get an additional effect of 
multiplication of $\chi(g_i)$. For any holomorphic polygon, which contributes to the
potential $W_\chi$, we should have $g_1 \cdots, g_k =1$ since the fiber should come back to
the original one after the whole set of turns.  It follows that $\chi(g_1) \cdots \chi(g_k) =1$.
\end{proof}
Note that the group action obtained as in the above is always a diagonal action.  We may summarize the above construction into the following definition.
\begin{defn}
The dual group $\WH{G}$  acts on  mirror space, leaves $W$ to be invariant,  whose action of $\chi \in \WH{G}$ on the variables corresponding
to the intersection $L\cap g L$ is given by the scalar multiplication of $\chi(g) \in U(1)$.
\end{defn}

\subsubsection{Non-abelian case}
When $G$ is non-abelian, it has higher-dimensional irreducible unitary representations $U$.  Correspondingly we can consider the non-abelian infinitesimal deformation space of the immersed Lagrangian $\bar{L}$, namely, $H \otimes \mathfrak{gl}(U)$.  Conceptually we are deforming the pair $(\bar{L},\underline{U})$, where $\underline{U}$ is the trivial bundle over $\bar{L}$ with trivial connection.  A formal deformation of an immersed sector is specified by a matrix in $\mathfrak{gl}(U)$.  When $U$ has dimension one, $\mathfrak{gl}(U) = \cpx$ and this reduces to the previous considerations.

Suppose that $X_i$'s are immersed sectors spanning $V$ and $b = \sum_i x_i X_i \in V$.  The coordinates $x_i$'s are now regarded as $\mathfrak{gl}(U)$-valued variables and hence \emph{non-commutative}.  We assume that $m^b_0 = W(b) 1_L$.  Then $W$ is a non-commutative series in $X_i$'s, and it is independent of the choice of representations $U$ of $G$ (as long as we regard $X_i$'s as non-commutative variables).  As before $W$ is given by disc counting.  The key point is that in such a non-commutative version of $W$, the disc counting not just take care which immersed sectors each disc pass through, but also the \emph{order} of the sectors.  Thus we obtain a \emph{non-commutative Landau-Ginzburg model} $W$.

Now consider the character group $\WH{G}$.
For the same reason as in Proposition \ref{identifypo}, each term in $W$ is invariant under the action of $G$: a term in $W$ is contributed by a stable disc which turns at some immersed sectors $X_{i_1}, \ldots, X_{i_k}$.  Each immersed sector $X_{i_j}$ corresponds to a group element $g_j$, and $g_1 \ldots g_k = 1$ because the Lagrangian goes back to itself after turning around the disc once.  In particular, $W$ is invariant under the action of $\WH{G}$.

To conclude, when $G$ is non-abelian, we have to consider the non-commutative version of $W$, which is disc-counting taking care of the order of immersed sectors on the boundary of the disc.  $W$ is invariant under the action of the dual group $\WH{G}$.  Since all the examples we consider in this paper involve only abelian groups, it suffices for us to consider commutative Landau-Ginzburg model.

%The first approach may be just to consider only one dimensional representations. Any homomorphism $G \to U(1)$ in fact factors through the abelianization $G/[G,G]$, and hence the dual group will be much smaller in a sense, losing certain informations.  The induced action of $\WT{G}$ will be a diagonal action as in the previous section, with Lemma analogous to \ref{identifypo}.

%If we consider higher dimensional representations, say $r$-dimensional irreducible representation $\chi$, then, each geometric generator of the Floer complex can give rise to several variables in the potential.  Namely, as we will consider Hom spaces between the fibers of orbi-bundles at the intersection point, we have $r^2$ worth of generators. Suppose we use all these generators, and define the associated potential function $W_\chi$.  Given another representation $\psi$, we will get $W_\psi$.The relationship of $W_\chi$ and $W_\psi$ is also given by $\psi(g) \chi(g^{-1})$, which is not just scalar multiplication, but could involve linear change of variables.  There should be an analogue of Lemma \ref{identifypo}, which probably means that these potentials can be identified up to linear change of variables.

%I guess that to consider an automorphism of the potential (to consider the dual group in the mirror),  instead of considering the whole structure, one could restrict to only the diagonal symmetries, which fixes the potential $ W$ but I could be wrong. And this seem to come from one dimensional representations.

\subsection{Mirror functors and $\WH{G}$-equivariance}\label{subsec:gequivfuct}

With an immersed Lagrangian $\bar{L}$ chosen in the global quotient $X = \tilde{X}/G$ and a choice of subspace $V$ of weak bounding cochains, we obtained a Landau-Ginzburg model $W$ which is invariant under $\WH{G}$ by Proposition \ref{identifypo}. 
We show in this section that our localized mirror functor can be lifted to give a map from the Fukaya category of $\tilde{X}$ to
the $\WH{G}$-equivariant matrix factorizations category of $W$. 

 Since we have the (ramified) covering map $\tilde{X} \to X$, we may use the terms upstairs and downstairs to refer to $\tilde{X}$ and $X$ respectively. We still make Assumption \ref{assump1} and that $G$ is Abelian ($\WH{G} \cong G$). For simplicity, we assume that the weak Maurer-Cartan space $V$ is generated by immersed sector (only $x_i$-variables).
 
Downstairs, for each compact oriented spin Lagrangian immersion $L' \subset X$ in the global quotient
(with a potential value $\lambda$, which we assume to be $0$ for simplicity), we obtain a corresponding matrix factorization of $W$
as explained in Section \ref{GMir},  \ref{geom_transf}, and  an $\AI$-functor $\Fuk_0(X) \to \mathcal{MF}(W)$.
% Namely, we perturb $L'$ by a Hamiltonian isotopy such that it intersects transversely with $\bar{L}$ and consider the vector space $F = CF^*(\bar{L},L')$ spanned by the intersection points $L \cap L'$ together with the endomorphism $m_1^{b,0}$, which can be written as a matrix $M$.  It gives a matrix factorization of $W$:
%$$M^2 = W.$$
%We say that $(F,M)$ is a matrix factorization mirror to $L'$.
%Conceptually $F$ can be identified as a trivial bundle $\uF$ over the mirror space $V$, and $M$ is identified as an endomorphism of this bundle.  Thus by our previous construction we obtain an $\AI$-functor downstairs:
%$$ \Fuk_0(X) \to \mathcal{MF}(W).$$

Upstairs, we would like to construct a mirror functor from Fukaya category of $\tilde{X}$ to the category of $\WH{G}$-equivariant matrix factorizations of the superpotential $\tilde{W}$. Recall that we have a family of Lagrangian $\pi^{-1}(\bar{L}) = \bigcup_{g \in G} g \cdot L$, with weak bounding cochain $b = \sum_{i=1}^{N} x_i X_i$ for $X$, or 
$\tilde{b} = \sum_{i=1}^{N} x_i \left(\sum_{g \in G} g \cdot \tilde{X}_i \right)$ for $\tilde{X}$.
By abuse of notation we may denote $\tilde{b}$ and $\tilde{X}_i$ simply by $b$ and $X_i$, when it is clear from the context that we are talking the geometry upstairs.

%
%
% In $\tilde{X}$ we have the preimage $\pi^{-1}(\bar{L})$ which is also an immersed Lagrangian.  Throughout this section, we assume that
%$$
%where $g \cdot L$'s are Lagrangian submanifolds of $\tilde{X}$ intersecting each other transversely.
%
%Each weak bounding cochain $b = \sum_{i=1}^{N} x_i X_i$
%downstairs induces a $G$-invariant weak bounding cochain
%$$\tilde{b} = \sum_{i=1}^{N} x_i \left(\sum_{g \in G} g \cdot \tilde{X}_i \right)$$
%upstairs, where $\tilde{X}_i$ is a certain chosen lift of the immersed sector $X_i$ downstairs.  By abuse of notation we may denote $\tilde{b}$ and $\tilde{X}_i$ simply by $b$ and $X_i$, when it is clear from the context that we are talking the geometry upstairs.

The same construction upstairs using $\tilde{\BL} = (\pi^{-1}(\bar{L}), \tilde{b})$ gives a generalized SYZ mirror $(V,\tilde{W})$ together with a mirror functor $\Fuk (\tilde{X}) \to \MF(\tilde{W})$. 

\begin{remark}
The counting of polygons  in  upstairs and downstairs are the same but their respective weights given by areas are related by a change of K\"ahler coordinates, since the K\"ahler parameters of $\tilde{X}$ and $X$ are different.
\end{remark}
 
%Every class $(\beta,\alpha)$ downstairs corresponds to $|G|$ disc classes $(\beta_g,\alpha_g)$ upstairs: we can choose a lift of $(\beta,\alpha)$ to $(\beta_1,\alpha_1)$ bounded by $\pi^{-1}(\bar{L})$, and the $G$ action on $(\beta_1,\alpha_1)$ produces $|G|$ disc classes.  These classes are pairwise distinct because $G$ acts freely on $\pi^{-1}(L)$, and hence $g \cdot \alpha_1 \not= \alpha_1$ whenever $g \not= 1$.  By definition each of the classes $(\beta_g,\alpha_g)$ have the same disc invariant as $(\beta,\alpha)$.  If the area of $\beta_g$ upstairs has a clear relationship with that of $\beta$ downstairs for every $\beta$, then we can express $\tilde{W}$ in terms of $W$ by multiplying by a factor of $|G|$ and change of K\"ahler parameters.  For instance for $X = \bP^1_{(3,3,3)}$, we have
%$$ \tilde{W} = |G| \cdot W(q_{\alpha}^3,q_{\beta}). $$

Now we would like to take the action of $G$ on $\tilde{X}$ into account and make the above to be a functor to the category $\MF_{\WH{G}}(\tilde{W})$ of $\WH{G}$-equivariant matrix factorizations.  To do this, we use the canonical equivalence of dg categories $ \MF_{\WH{G}}(\tilde{W}) \cong \Tw(B_W\#\WH{G})$
which is proved by Tu \cite[Section 6]{Tu2},
where $B_W$ is the coordinate ring $\C[x_1,\ldots,x_N]$, and $B_W\#\WH{G}$ is a dg category whose objects are characters of $\WH{G}$ (which is $G$) defined by Caldararu-Tu \cite[Section 2.15]{CT}.  (In \cite{Tu2} $B_W\#\WH{G}$ was taken to be an algebra instead of a category, but we will use the setting of \cite{CT} because it is more canonical.) 

Thus for each $g \in G$ we have a corresponding object $A_g$ of $B_W\#\WH{G}$.  The morphism space
$\Hom_{\WH{G}} (A_{g_1},A_{g_2}) $
is defined as the space of polynomials in $\C[x_1,\ldots,x_N]$ which is invariant under the $(g^{-1}_1,g_2)$-twisted action of $\WH{G}$.  

Recall from Lemma \ref{g_X} that each intersection point $X_i$ corresponds to a group element $g_{X_i} \in G$, and the dual
group action of  $\chi \in \WH{G}$ on the corresponding variable $x_i$ is given by $ \chi \cdot x_i = \chi(g_{X_i}) x_i.$
 Then the $(g^{-1}_1,g_2)$-twisted action of $\chi \in \WH{G}$ is defined as
$$ \chi \cdot_{(g^{-1}_1,g_2)} x_i := \chi(g^{-1}_1) \big(\chi \cdot (\chi(g_2) x_i)\big) = \chi(g^{-1}_1 g_2) \big( \chi \cdot x_i \big) . $$
``$\Tw$" denotes the twisted-complex construction where differentials of a complex squares to $W$ instead of $0$.

Given a Lagrangian $L' \subset \tilde{X}$ upstairs (which can be assumed to be transversal to $\pi^{-1}(\bar{L})$ by applying Hamiltonian perturbation if necessary), the corresponding matrix factorization was defined to be $(CF^*(\pi^{-1}(\bar{L}),L'),m_1^{(b,0)}).$
To make it $\WH{G}$-equivariant, we `categorize' it to record the relations between different branches of $\pi^{-1}(\bar{L})$ by $G$-action.  Namely:
\begin{enumerate}
\item For each intersection point between $L'$ and $g \cdot L$ for some $g \in G$, we make a copy of the object $A_g$ of $B_W\#\WH{G}$.  
\item Sum over all intersection points $p$ to obtain an object of $\Sigma(B_W\#\WH{G})$
$$F = \bigoplus_{p \in L' \cap \pi^{-1}(\bar{L})} A_{g(p)}$$
where for an intersection point $p$, it is associated to a group element $g(p) \in G$ such that $p \in L' \cap (g(p) \cdot L)$.
\item The differential $\delta$ is defined by $m_1^{\tilde{\BL},L'}$ as before, which has the property $\delta^2 = \tilde{W}$.   The higher terms of the $\AI$-functor are defined by $m_k^{(L,b),\cdot,\ldots,\cdot}$ as before.
\end{enumerate}

In summary, the above construction `categorizes' $(CF^*(\pi^{-1}(\bar{L}),L'),m_1^{(\BL,L')})$ by the partition
$$ (\pi^{-1}(\bar{L})) \cap L' = \bigcup_{g \in G} (g\cdot L) \cap L'.$$

\begin{prop}\label{prop:gendualact}
The above defines an $A_\infty$ functor
$$\mathcal{LM}^\BL: \Fuk_0(\tilde{X}) \to \mathcal{MF}_{\WH{G}}(\tilde{W}) $$
where $\Fuk_0(\tilde{X})$ denotes the subcategory of Lagrangians which are unobstructed and $\mathcal{MF}_{\WH{G}}(\tilde{W})$ denotes the category of $\WH{G}$-equivariant matrix factorizations.
\end{prop}

\begin{proof}
We need to check that the matrix factorization $(F,\delta)$ defined above is an object of $$\Tw(B\#\WH{G}) \cong \MF_{\WH{G}}(\tilde{W}).$$
In other words, $\delta$ belongs to $\Hom_{\WH{G}}(\bigoplus_{p \in I} A_{g(p)}, \bigoplus_{p \in I} A_{g(p)})$, namely, polynomials involved in $\delta$ are twisted $\WH{G}$-invariant in the sense defined above.

Consider a holomorphic polygon bounded by $\BL$ and $L'$ contributing to $m_1^{(\BL,L')}$, see the left hand side of Figure \ref{equivariance}.  Let $g_1 \cdot L, \ldots, g_k \cdot L$ be the branches of $\BL$ bounding the polygon labelled counterclockwisely, and let $X_1,\ldots,X_k$ be the involved intersection points.  The input point is at $g_k \cdot L$ and the output point is at $g_1 \cdot L$.  The monomial contributed is $x_1\cdots x_{k-1}$.  We want to deduce that $x_1\cdots x_{k-1}$ is invariant under the $(g_k^{-1},g_1)$-twisted action of $\WH{G}$.

\begin{figure}[htp]
\includegraphics[scale=0.4]{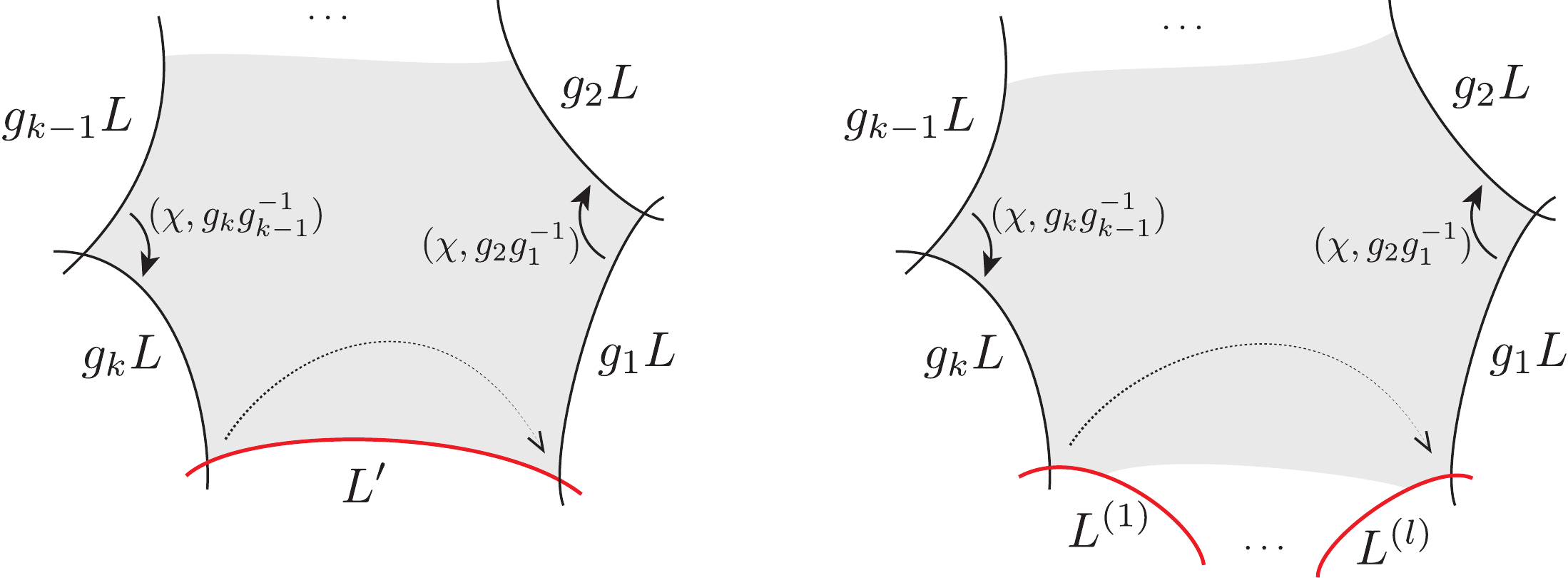}
\caption{$\WH{G}$-equivariance.}
\label{equivariance}
\end{figure}

The action of $\chi \in \WH{G}$ on $x_i$ is
$$ \chi \cdot x_i =\chi(g_{i+1}g_i^{-1}) x_i $$
and so
$$ \chi \cdot (x_1 \ldots x_{k-1}) = \chi(g_{k}g_1^{-1}) (x_1 \cdots x_{k-1}).$$
Then the twisted action is
$$ \chi \cdot_{(g_k^{-1},g_1)} (x_1 \cdots x_{k-1}) = \chi(g_k^{-1}) \chi(g_{k}g_1^{-1}) \chi(g_1)(x_1 \cdots x_{k-1})= x_1 \cdots x_{k-1}.$$

The same argument, by replacing $L'$ with an ordered collection of Lagrangians $(L^{(1)},\cdots,L^{(l)})$ (see the right hand side of Figure \ref{equivariance}), shows that the higher terms $m_k^{\BL,L^{(1)},\ldots,L^{(l)}}$ of the $\AI$-functor have their targets to be the $\WH{G}$-equivariant morphism space $$\Hom_{\WH{G}}\left(\bigoplus_{p \in L^{(1)} \cap \pi^{-1}(\bar{L})} A_{g(p)}, \bigoplus_{q \in L^{(l)} \cap \pi^{-1}(\bar{L})} A_{g(q)}\right).$$
\end{proof}

\section{The elliptic curve and its $\Z/3$-quotient $\mathbb{P}^1_{3,3,3}$}\label{ex333}
In the rest of this paper, we apply our general theory to a series of interesting examples.  In this section we consider the elliptic curve quotient $\mathbb{P}^1_{3,3,3}$.  This produces the mirror superpotential $W: \C^3 \to \C$.  We will construct the mirror of $\mathbb{P}^1_{a,b,c}$ for general $a,b,c$ and prove homological mirror symmetry in Section \ref{sect_dim1}.  Since $\mathbb{P}^1_{3,3,3}$ serves as a nice illustration and its mirror map is particularly interesting, we separate it into a single section.

Consider the elliptic curve $E = \C / (\Z + \Z \tau)$ with $\tau = e^{2\pi \sqrt{-1}/3}$. The curve $E$ admits a $\Z/3$ action, which
is generated by the complex multiplication of $\tau$ on $\C$. The quotient is the orbisphere $\mathbb{P}^1_{3,3,3}$ with three $\Z/3$ orbifold points,
and the quotient map is denoted as $\pi: E \to \mathbb{P}^1_{3,3,3}$.

We take the $\Z/3$-equivariant Lagrangian immersion $\WT{\iota}: \WT{L} \times \Z/3  \to \C$
where $\WT{L} = \R$, and $$\WT{\iota}(x, 1) = \frac{1+\tau}{2} + \sqrt{-1}x, \quad \WT{\iota}(x, \tau) = \tau \WT{\iota}(x,1),\quad  \WT{\iota}(x, \tau^2) = \tau^2 \WT{\iota}(x, 1).$$
Let $pr : \C \to E$ be the quotient map, and set $L := \WT{L} / Im(\tau) \Z$. Then  
$$\iota := pr \circ \WT{\iota} : L \times \Z/3 \to E$$
defines a $\Z/3$-equivariant  Lagrangian immersion. By taking $\Z/3$-quotient, we obtain a Lagrangian immersion (also denoted as $\bar{L}$) $\OL{\iota}:= \pi \circ \iota: S^1 \to \mathbb{P}^1_{3,3,3}.$  %We will equip it with a non-trivial spin structure.
Figure \ref{qqalpha}  illustrates the universal cover $\C$ of the elliptic curve $E$ and the equivariant immersion $\iota$. See Figure \ref{SeiLagr} for its image in the orbi-sphere $\mathbb{P}^1_{3,3,3}$.

There is another $\Z/3$-equivariant Lagrangian immersion
given by the imaginary axis of $\C$, and its $\Z/3$-images. 
More precisely, denote
 $\iota_{L_l} : L \times \Z/3  \to \C$ with  $\iota_{L_l} (x, 1) =  \sqrt{-1}x$.
 We will use $\bar{L}_l$ to denote this Lagrangian immersion, and call it long diagonal. (In fact,  each branch  $L_l$ of $\iota_{L_l}$ was called long diagonal in \cite{B}). 
Note that $L$ and $L_l$  are related by translation in the  cover $\C$, and they are different objects in Fukaya category, since translation in $E$ is not a Hamiltonian isotopy. 
 
Another $\Z/3$-equivariant Lagrangian immersion given by the real axis of $\C$
will be called short diagonal, and denoted as $\bar{L}_s$.

First, we will give explicit computation of the Lagrangian Floer potential of $\bar{L}$ to find the LG mirror
$W_{3,3,3}$, by counting suitable holomorphic triangles in $\C$.

This enables us to show that the generalized SYZ map equals to the mirror map for $X= \proj^1_{3,3,3}$ (\cite{Saito},\cite{Milanov-Ruan,ST}), which means that our work provides an enumerative meaning of the mirror map, and also explains the integrality of the mirror map

We will also study the localized mirror functor for $\bar{L}$ in detail in this section.
Alternatively one could study the localized mirror functor for $\bar{L}_l$.
The localized mirror functor for  $\bar{L}$ is indeed a genuine homological mirror functor
for $\mathbb{P}^1_{3,3,3}$ as we will show in the next section that $\bar{L}$ split-generates Fukaya category of $\mathbb{P}^1_{3,3,3}$,

We show that under this mirror functor, the long diagonal $\bar{L}_l$ goes to a $3 \times 3$ matrix factorization. 
and the short diagonal $\bar{L}_s$ to a $2 \times 2$ matrix factorization. And $\bar{L}$ itself
goes to a $4 \times 4$ matrix factorization, which is computed explicitly in the next section.

\subsection{Quantum-corrected superpotential}
\begin{figure}[htp]
\begin{center}
\includegraphics[scale=0.45]{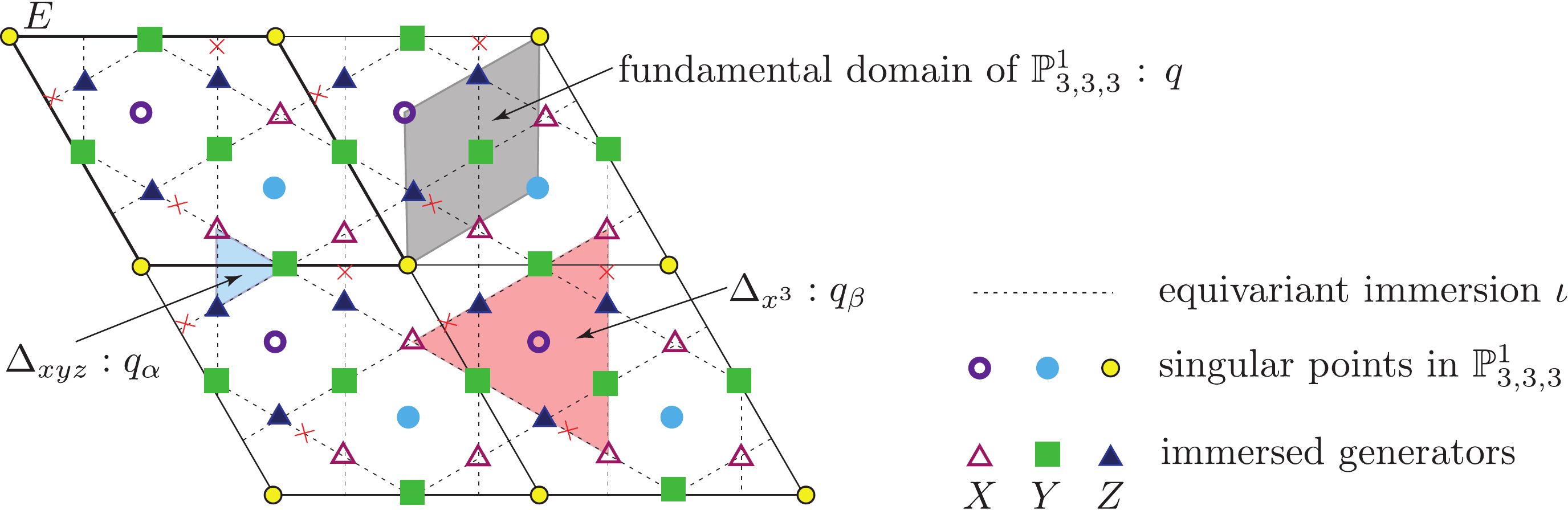}
\caption{Universal cover of $\mathbb{P}^1_{3,3,3}$ and the equivariant immersion $\iota$}\label{qqalpha} 
\end{center}
\end{figure}

The $\AI$-algebra for the Lagrangian immersion $\oi$ into $\mathbb{P}^1_{3,3,3}$ is defined as follows. We first consider three Lagrangians in $E$
given by components of the $\Z/3$-equivariant immersion $\iota$, and construct Fukaya sub-category of these three Lagrangian submanifolds,
as explained in Section \ref{sec:fuksur}. This subcategory has a natural $\Z/3$-action, and
 the $\Z/3$-invariant part provides a required $\AI$-algebra structure on the Floer complex of $\bar{\iota}$.

We have three  immersed points of $\bar{\iota}$  in  $\mathbb{P}^1_{3,3,3}$.
Each immersed point gives rise to two generators of the Floer complex, namely $X,Y,Z$, the odd-degree generators, and $\OL{X}, \OL{Y}, \OL{Z}$, the even-degree generators of $HF(\oi, \oi)$.
Also, we choose a Morse-function $f$ on the domain of $\bar{\iota}$ ($\cong S^1$), with one minimum $e$ and one maximum $p$. We will think of $e$ as
a unit and $p$ as a point class of the Lagrangian. The Floer cohomology $HF(\oi, \oi)$ is generated by these 8 elements.

\begin{lemma}
An element $b=xX + yY +zZ$ for any $x,y,z \in \Lambda$ is a weak bounding cochain.
\end{lemma}
This will be proved in more general setting in Theorem \ref{thm:weak}.
Now we compute the mirror superpotential. 
For this we have to compute $m_k$-operations in $E$ modulo $\Z /3$-action with inputs
given by lifts of $X,Y$ and $Z$. 
\begin{lemma}
Let $X_i \in \{X,Y,Z\}$ for $i=1,\cdots, k$. 
Then $m_k(X_1, X_2,\cdots, X_k)=0$ if $k\geq 4$.
\end{lemma}
\begin{proof}
As inputs are immersed generators, $m_k$ counts holomorphic $k$ or $k+1$-gons mapping into $\mathbb{P}^1_{3,3,3}$ with
appropriate boundary conditions.
As the boundary of the holomorphic polygon is trivial in the orbifold fundamental group of $\mathbb{P}^1_{3,3,3}$, we can take a lift of this polygon to the universal cover $\C$. Now let us consider the index. If we denote by $\beta$ the homotopy class of a holomorphic polygon, 
then the expected dimension of $m_{k,\beta}(X_1,\cdots, X_k)$ is $\mu(\beta) + 1 -3+1$. Since the Lagrangian immersions are straight lines, and only turns at the immersed generators,  we have
$\mu(\beta) \geq \frac{k(\pi - \pi/3)}{\pi}$. This implies that $m_k(X_1,\cdots, X_k)$ are nontrivial only if $k \leq 3$.
\end{proof}

Since the potential $PO$ is defined by $ m(\conste^b) = PO(\oi,b) e$ for weak bounding cochains $b$,
it is enough to determine the coefficients of $e$ in $m_k(X_1,\cdots, X_k)$ for $k \leq 3$.
The coefficients of $e$ in $m_1(X_1)$ and $m_2(X_1,X_2)$ are obviously zero for any degree one immersed generators $X_1, X_2$. Thus, it suffices to compute the contribution of $m_3(X_1,X_2,X_3)$ to $e$. 
Such contributions come from triangles (with $\mu(\beta)=2$) in $\C$ passing through one of
pre-images of $e$ under the covering map.

%\subsubsection{Minimal Potential}
Suppose that the (image under $\bar{\iota}$ of) minimum $e$ of the chosen Morse function on $S^1$ lies on the arc $\arc{XZ}$ and the maximum $p$ is the reflection of $e$ across the long diagonal (equator of $\mathbb{P}_{(3,3,3)}$).
In Figure \ref{qqalpha}, $e$ is represented by {\color{red} $\times $}. 
 Recall that the Lagrangian $\bar{\iota}$ is equipped with a non-trivial spin structure. We can enforce this by picking a generic point (for convenience ) right next to $e$ in the
arc $\arc{XZ}$. Taking the effect of the non-trivial spin structure into account, we reverses the sign of  moduli spaces if the boundary of
holomorphic polygons pass through this generic point. Since the generic point is right next to $e$, we may regard $e$ as the generic point for this purpose. 

Now we consider smallest triangles contributing to the potential.
For this,  pick a lift $\WT{e}$ of $e$ in $\C$ and find triangles with turns given by $X,Y,Z$,
and passing through $\WT{e}$.
It is east to
see from Figure \ref{qqalpha} that there exist 4 types of triangles passing through $\WT{e}$,
which contribute to 
\begin{equation}\label{trim3}
m_3(X,X,X),\quad m_3(Y,Y,Y),\quad m_3(Z,Z,Z),\quad m_3(X,Y,Z).
\end{equation}
Note that  triangles with corners given by $X,Z,Y$ in a counter-clockwise order do not
meet any lift of $e$ at its boundary, hence is not counted.

The smallest triangle with corners in the order $X,Z,Y$ has a unique lift of $e$ on its boundary, and hence contributes to $m_3$ exactly once. However, the smallest triangle with corners $X,X,X$ has three pre-images of $e$ on its boundary. Nevertheless, as a disc on the orbifold $\mathbb{P}^1_{3,3,3}$, not on the cover, this triangle contributes
also once due to its symmetry.

To record symplectic areas, we denote by $q$ the inverse of the exponential of the area of $\bP^1_{3,3,3}$, and by $q_\alpha$ that of the minimal $XYZ$ triangle $\Delta_{xyz}$,
and by $q_\beta$ that of minimal $X^3$ triangle $\Delta_{X^3}$.
We have relations $q_\alpha^8 = q$ and $q_\beta = q^9$.
See Figure \ref{qqalpha}. 

Finally, we assign the signs to these polygons according to the rule given in Section \ref{sec:fuksur}. The full superpotential takes the form
%Then, the minimal area triangles contribute to $PO(\oi, b)$
%as
\begin{equation}\label{eq:W333o25}
 W_{3,3,3}=- q_\beta x^3 + q_\beta y^3 - q_\beta z^3 - q_\alpha xyz + O(q_\alpha^{25}).
\end{equation}
%Roughly speaking, the negative sign for $x^3, y^3, xyz$ comes from the non-trivial spin structure, and
%the positive sign for $y^3$ comes from the additional sign from the mismatch of the orientation of
%the Lagrangian immersion with that of 
%the boundary of the holomorphic polygon. We will examine the precise sign rule in the next subsection.
In addition to these, there are infinite sequences of triangles for each cubic term in Equation \eqref{eq:W333o25}. By counting these triangles weighted by areas and signs, we obtain the following theorem.

\begin{theorem}\label{thm:pot333formula}
The Lagrangian Floer potential $PO(\oi, b)$ is given by 
\begin{equation}\label{superpotential}
W_{3,3,3} (x,y,z) = \phi (q_\alpha^3) (x^3 - y^3 + z^3)  + \psi (q_\alpha) xyz,
 \end{equation}
where $b= xX + yY + z Z$ and 
$$\phi (q_\alpha) = \sum_{k=0}^\infty  (-1)^{k+1} (2k+1) q_\alpha^{(6k+3)^2},$$
$$\psi(q_\alpha) = -q_{\alpha} + \sum_{k=1}^\infty  (-1)^{k+1} \left( (6k+1) q_\alpha^{(6k+1)^2} - (6k-1) q_\alpha^{(6k-1)^2} \right).$$
\end{theorem}

The precise form of the potential above has not been computed previously, and is different from the folklore mirror potential $x^3 + y^3 +z^3 - \sigma xyz$.
One can regard $W_{3,3,3}$ as a quantum corrected mirror potential. The benefit
is that we can compute the mirror map explicitly from this construction. We will show in Section \ref{sec:syzmirror} that the ratio of $\phi (q_\alpha)$ and $\psi(q_\alpha)$ gives rise to the mirror map in classical mirror symmetry.

On the other hand, closed string mirror symmetry can be explained via an closed to open map. For toric manifolds, Fukaya-Oh-Ohta-Ono \cite{FOOO_MS} used holomorphic discs
with one interior insertion from quantum cohomology cycle to define a map from the quantum cohomology
to the Jacobian ring.  The same line of proof should work in this case, or in general  $\bP^1_{a,b,c}$
%(which will be discussed in the next section)
to derive the isomorphism
\begin{equation}\label{eq:ks1}
\mathfrak{ks}: QH^*(\bP^1_{a,b,c}) \cong \Jac (W_{a,b,c}). 
\end{equation}
We hope to discuss this in more detail in the future.
Given the explicit computations of $W_{a,b,c}$ in \cite{CHKL14}, the isomorphism
leads to an explicit representation of the quantum cohomology ring.

\subsection{Generalized SYZ map equals to the mirror map}\label{sec:syzmirror}

In this subsection, we show that the generalized SYZ map equals to the mirror map for $X= \proj^1_{3,3,3}$.  The mirror map for $\proj^1_{3,3,3}$ is well-known from Saito's theory \cite{Saito}, and readers are referred to \cite{Milanov-Ruan,ST} for more recent treatments.  Our work here gives an enumerative meaning of the mirror map.  Moreover, it explains the integrality of the mirror map. i.e. its coefficients when expanded in the K\"ahler parameter $q$ are all integers.

The mirror map is given by the quotient $q(\check{q})=\exp\left( \frac{2\pi\consti}{3} \cdot \pi_B (\check{q}) / \pi_A (\check{q})\right)$, where $\pi_A (\check{q})$ and $\pi_B (\check{q})$ are functions satisfying the Picard-Fuchs equation
$$ u''(\check{q}) + \frac{3\check{q}^2}{\check{q}^3 + 27}u'(\check{q}) + \frac{\check{q}}{\check{q}^3 + 27}u(\check{q}) = 0. $$
See Section 6.2 of \cite{Milanov-Ruan} for the explicit expressions of $\pi_A (\check{q})$ and $\pi_B (\check{q})$.  The mirror map is of the form
$$ q = - \check{q}^{-1} \left( 1 + \sum_{k=1}^\infty c_k (-\check{q})^{-3k} \right). $$
By inverting the above series, one obtains $ \check{q} = -3 a(q)$, where
$$ a(q) = 1 + \frac{1}{3}\left( \frac{\eta(q)}{\eta(q^9)} \right)^3 = 1 + \frac{1}{3} q^{-1} \left( \frac{\prod_{k=1}^\infty (1-q^k)}{\prod_{k=1}^\infty (1 - q^{9k})} \right)^3. $$
Thus the Landau-Ginzburg mirror to $(\proj^1_{3,3,3}, \omega_q)$ is $$ W^{\mathrm{LG}} = (x^3 + y^3 + z^3) + \check{q}(q) xyz. $$

\begin{theorem}\label{mir=SYZ}
The generalized SYZ mirror $W^{\mathrm{SYZ}}(x,y,z):=W_{3,3,3} (x,y,z)$ equals to the Landau-Ginzburg mirror $W^{\mathrm{LG}}(x,y,z)$ for $\proj^1_{3,3,3}$, up to a coordinate change in $(x,y,z)$.
\end{theorem}

\begin{proof}
Recall that the generalized SYZ mirror that we have constructed for $\proj^1_{3,3,3}$ is
\begin{equation}\label{eq:yminusy}
 \phi (x^3 + y^3 + z^3) - \psi xyz 
\end{equation}
where
$$ \phi = \sum_{k=0}^{\infty} (-1)^{3k+1} (2k+1) q_\alpha^{3(12 k^2 + 12 k + 3)} $$
and
$$ \psi = -q_\alpha + \sum_{k=1}^{\infty} \left( (-1)^{3k+1} (6k+1) q_\alpha^{(6k+1)^2} + (-1)^{3k}(6k-1)q_\alpha^{(6k-1)^2} \right). $$
Here we made a coordinate change $y \mapsto -y$. Compare Equation \eqref{eq:yminusy} with Equation \eqref{superpotential}.

By a change of coordinates $(x,y,z) \mapsto \phi^{-1/3}(x,y,z)$, the generalized SYZ mirror can be written as
$$ W^{\mathrm{SYZ}} = (x^3 + y^3 + z^3) - \frac{\psi}{\phi} xyz. $$
Recall that $q = q_\alpha^8$.
 
Now to prove $W^{\mathrm{SYZ}} = W^{\mathrm{LG}}$, it suffices to check the equality
$$ \frac{\psi}{\phi} = - \check{q}(q_\alpha^8) = 3 + q_\alpha^{-8} \prod_{k=1}^\infty \frac{(1-q_\alpha^{8k})^3}{(1 - q_\alpha^{72k})^3}.$$
%The equality can be rearranged as
%$$\phi(q_\alpha) \prod_{k=1}^\infty (1-q_\alpha^{8k})^3 = q_\alpha^{8} \big(\psi(q_\alpha) - 3 \phi(q_\alpha)\big) \prod_{k=1}^\infty (1 - q_\alpha^{72k})^3. $$
This can be derived directly by using the following identity
$$ \prod_{k=1}^\infty (1 - q_\alpha^{k})^3 = \sum_{k=0}^\infty (-1)^k (2k+1) q^{\frac{k(k+1)}{2}}.$$
%Then the left hand side becomes
%$$ \left( \sum_{k=0}^\infty (-1)^k (2k+1) q^{4k(k+1)} \right) \left(\sum_{l=0}^{\infty} (-1)^{l+1} (2l+1) q_\alpha^{9(2l+1)^2} \right) $$
%and the right hand side is
%\begin{align*}
%& q_\alpha^{8} \left( \sum_{l=0}^\infty (-1)^l (2l+1) q^{36l(l+1)} \right)\\
%& \cdot \left( \sum_{k=0}^{\infty} (-1)^{3k+1} (6k+1) q_\alpha^{(6k+1)^2} + \sum_{k=0}^{\infty} (-1)^{3k+1}(6k+5)q_\alpha^{(6k+5)^2} - \sum_{k=0}^{\infty} (-1)^{3k+1} (6k+3) q_\alpha^{(6k+3)^2} \right)\\
%=& \left(\sum_{l=0}^{\infty} (-1)^{l} (2l+1) q_\alpha^{9(2l+1)^2} \right) \left( \sum_{k=0}^\infty (-1)^{k+1} (2k+1) q^{4k(k+1)} \right).
%\end{align*}
%Therefore, the equality follows.
\end{proof}

Note that the coefficients of $\phi$ and $\psi$ are all integers.  Moreover the leading coefficient of $\phi$ is $-1$.  Hence $-\psi / \phi$ have integer coefficients.  By the above theorem, this implies that the mirror map $\check{q}(q)$ also has integer coefficients.  Thus,
\begin{corollary}
The mirror map $\check{q}(q)$ as a series in $q$ has integer coefficients.
\end{corollary}

\subsection{Mirror functor to matrix factorizations}\label{sec:mf333}
Now we explain our localized mirror functor for $\mathbb{P}^1_{3,3,3}$. We will transform several
Lagrangian branes into the corresponding matrix factorizations using our mirror functor.

The Lagrangian immersion $\oi $ (or $\bar{L}$) together with the weak bounding cochains $b = xX + yY + zZ$
gives rise to an $\AI$ functor from the Fukaya category $\mathcal{F}uk_0(\mathbb{P}^1_{3,3,3})$ to the category of matrix factorizations of $W_{3,3,3}$, 
$\mathcal{LM}^{(\oi, b)} : \mathcal{F}uk_0(\mathbb{P}^1_{3,3,3})  \to \mathcal{MF}(W_{3,3,3}).$

There are three  notable Lagrangian objects in $\mathcal{F}uk_0(\mathbb{P}^1_{3,3,3})$:
$\oi$ itself, the long diagonal $\bar{L}_l$, and  the short diagonal $\bar{L}_s$, which
are defined in the beginning of this section.

These are mapped under the functor $\CF^{(\oi, b)}$ to matrix factorizations.
We will compute the matrix factorizations of $\bar{L}_l$ and $\bar{L}_s$ in the following.  See Proposition \ref{prop:MF333} for the matrix factorization transformed from the Seidel Lagrangian.

\subsubsection{Mirror matrix factorization of a long brane $\bar{L}_l$}
\begin{figure}[htp]
\begin{center}
\includegraphics[scale=0.4]{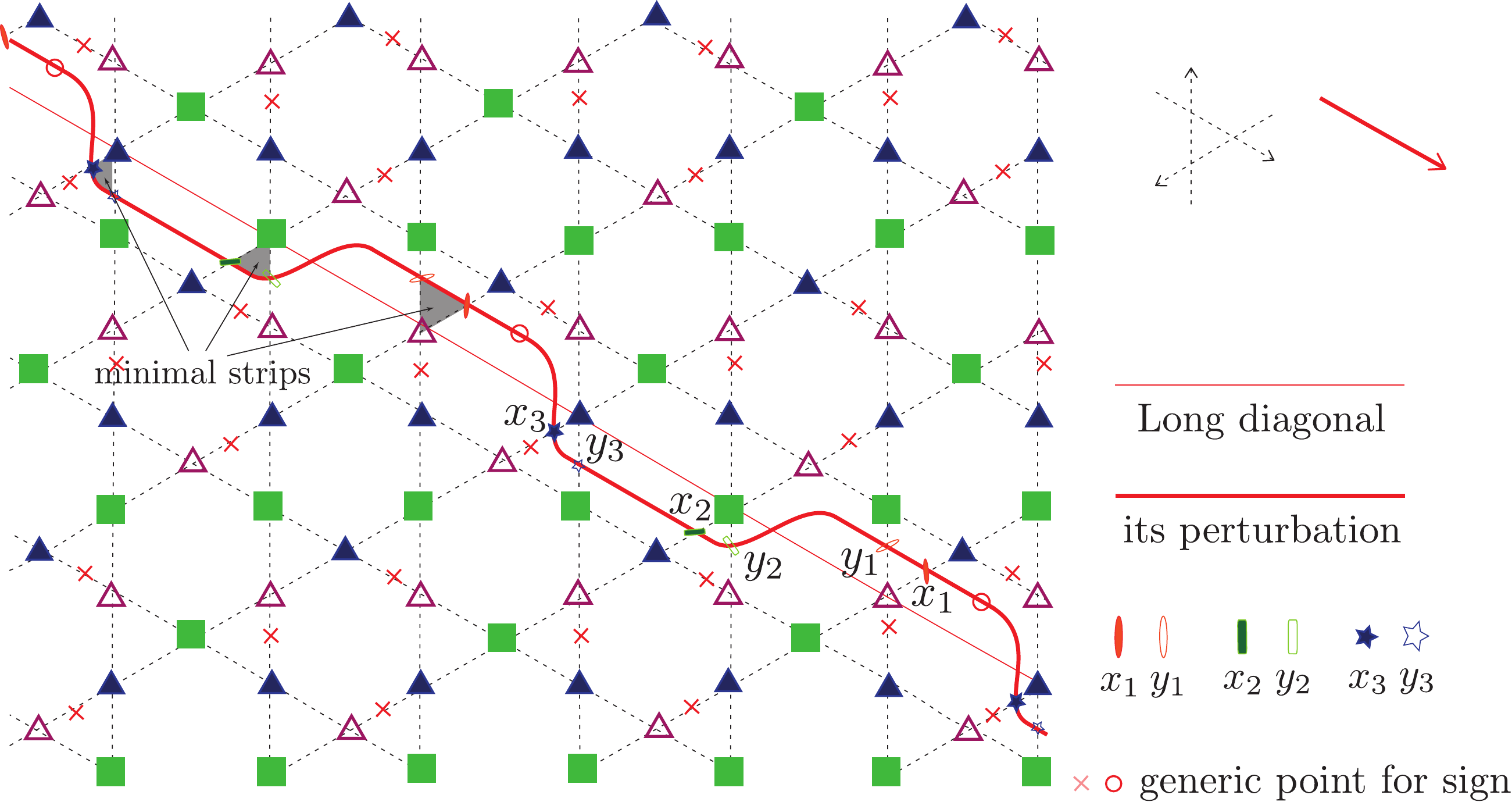}
\caption{Long diagonal and its perturbation and their orientations}\label{longdiag} 
\end{center}
\end{figure}

Let us first consider the long diagonal $L_l$ and find its mirror matrix factorization.
We will equip $\bar{L}_l$ with a non-trivial spin structure as we did for $\bar{L}$ which makes the signs of counting more regular. (See {\color{red}$\circ$} on $\bar{L}_l$ in Figure \ref{longdiag}.) 

%Since $\bar{L}_l$ is an equator in the orbi-sphere, it is easy to see that  $\bar{L}_l$ and the Lagrangian $\bar{L}$ have three triple intersections (see Figure \ref{SeiLagr}). 
It is possible to define $CF((\bar{L},b), \bar{L}_l)$ directly using the Morse-Bott model, but we work with the perturbation $\bar{L}_l^\epsilon$ of $\bar{L}_l$ by a Hamiltonian isotopy as in Figure \ref{longdiag} and take the limit $\epsilon \to 0$. 

The new Lagrangian $\bar{L}_l^\epsilon$ intersects $\bar{L}$ transversely at 6 points as illustrated in the Figure  \ref{longdiag}.
More precisely, if we perturb $\bar{L}_l$ to $\bar{L}_l^\epsilon$ around a neighborhood of each triple intersection of $\bar{L}_l$ and $\bar{L}$, 
we get a pair of  odd and even intersection points, denoted by $(x_i,y_i)$ for $i=1,2,3$.

The matrix factorization obtained from $CF( (\bar{L}, b), \bar{L}_l)$ (where $b$ varies in $\mathcal{M}_{weak}^1 (\bar{L})$) is of the following form after taking the limit $\epsilon \to 0$.

\begin{equation}\label{longdiagE}
E=
\bordermatrix{
       &  x_1   		 &  x_2             &  x_3               \cr
y_1 &  \alpha_1 x    &  \alpha_2 z  &  -\alpha_3 y    \cr
y_2 &  \alpha_3 z    &  -\alpha_1 y  &  \alpha_2 x    \cr
y_3 &  -\alpha_2 y    &  \alpha_3 x  &  \alpha_1 z    \cr}
\end{equation}
\begin{equation}\label{longdiagJ}
J=
\bordermatrix{
       &  y_1   	&  y_2  & y_3  \cr
x_1 &\beta_1 x^2-\gamma_1 yz &\beta_3 z^2-\gamma_3 xy &\beta_2 y^2+\gamma_2 zx  \cr
x_2 &\beta_2 z^2 -\gamma_2 xy &\beta_1 y^2+\gamma_1 zx &\beta_3 x^2-\gamma_3 yz   \cr
x_3 &\beta_3 y^2+\gamma_3 zx &\beta_2 x^2-\gamma_2 yz &\beta_1 z^2-\gamma_1 xy   \cr}
\end{equation}
where the $(i,j)$-th entry of $E$ is given by $\langle m_1^{b,0}(x_j), y_i \rangle$, and
the $(i,j)$-th entry of $J$ is given by $\langle m_1^{b,0}(y_j), x_j \rangle$.
Here $y_i$-coefficients of $m_1^{b,0}(x_j)$ is denoted as $\langle m_1^{b,0}(x_j), y_i \rangle$.

\begin{figure}[htp]
\begin{center}
\includegraphics[scale=0.35]{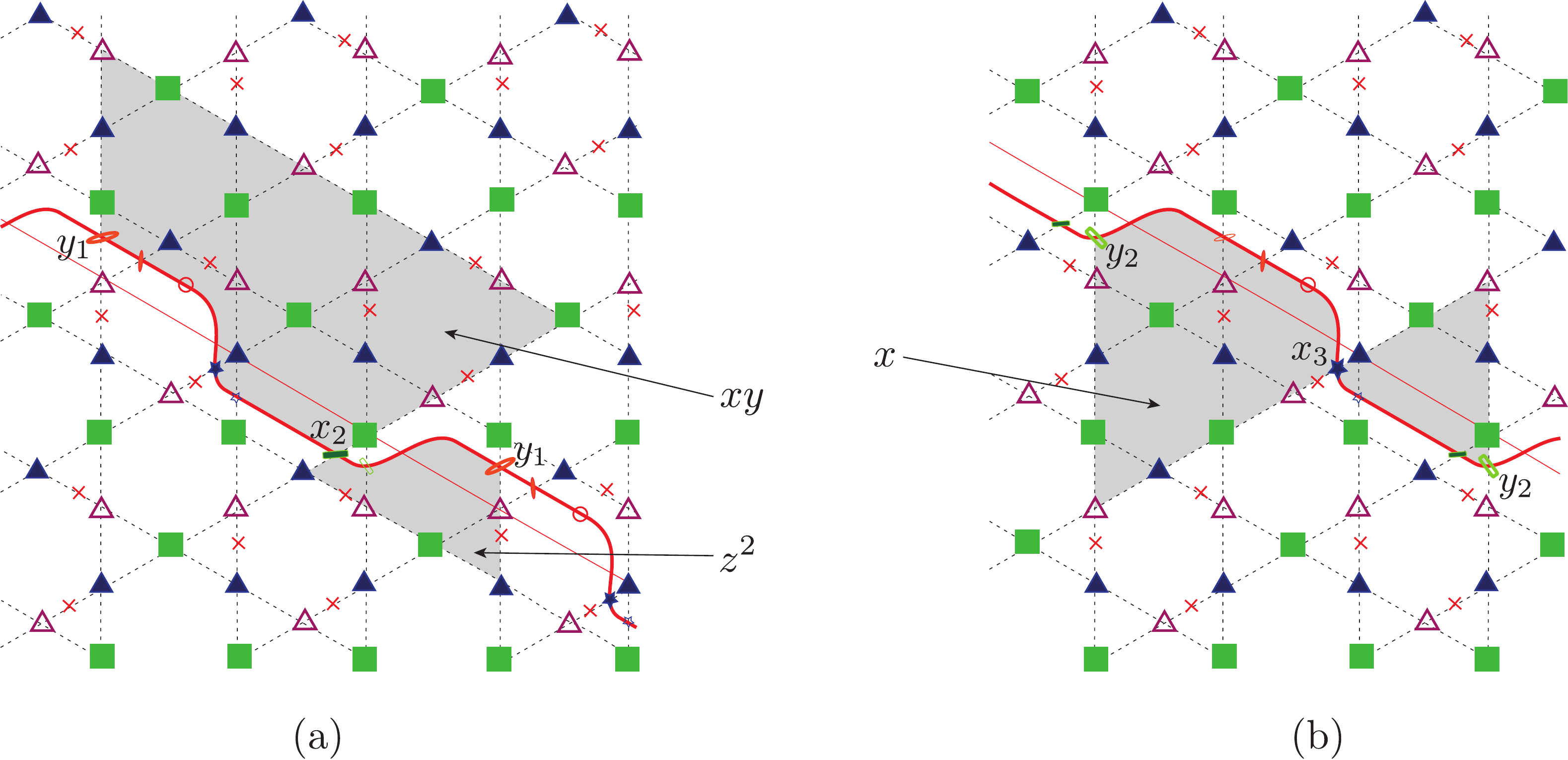}
\caption{Some strips from (a) $y_1$ to $x_2$ and (b) from  $x_3$ to $y_2$}\label{fig:ldx2y1} 
\end{center}
\end{figure}

For example, the $(2,1)$-th entry of $J$ is obtained by counting strips from $y_1$ to $x_2$, two of which are shown in Figure \ref{fig:ldx2y1} (a), and the $(2,3)$-th entry of $E$ by counting strips from $x_3$ to $y_2$ as shown in Figure \ref{fig:ldx2y1} (b).

The coefficients $\alpha_i$, $\beta_i$ and $\gamma_i$ ($i=1,2,3$) are given by
countable sums from sequences of strips, which are explicitly calculable if we take the limit $\epsilon \to 0$. 
As an example, let us explain the computation of  $\langle m_1^{b,0}(x_1), y_1\rangle$ or $\alpha_1x$.
\begin{figure}[htp]
\begin{center}
\includegraphics[scale=0.4]{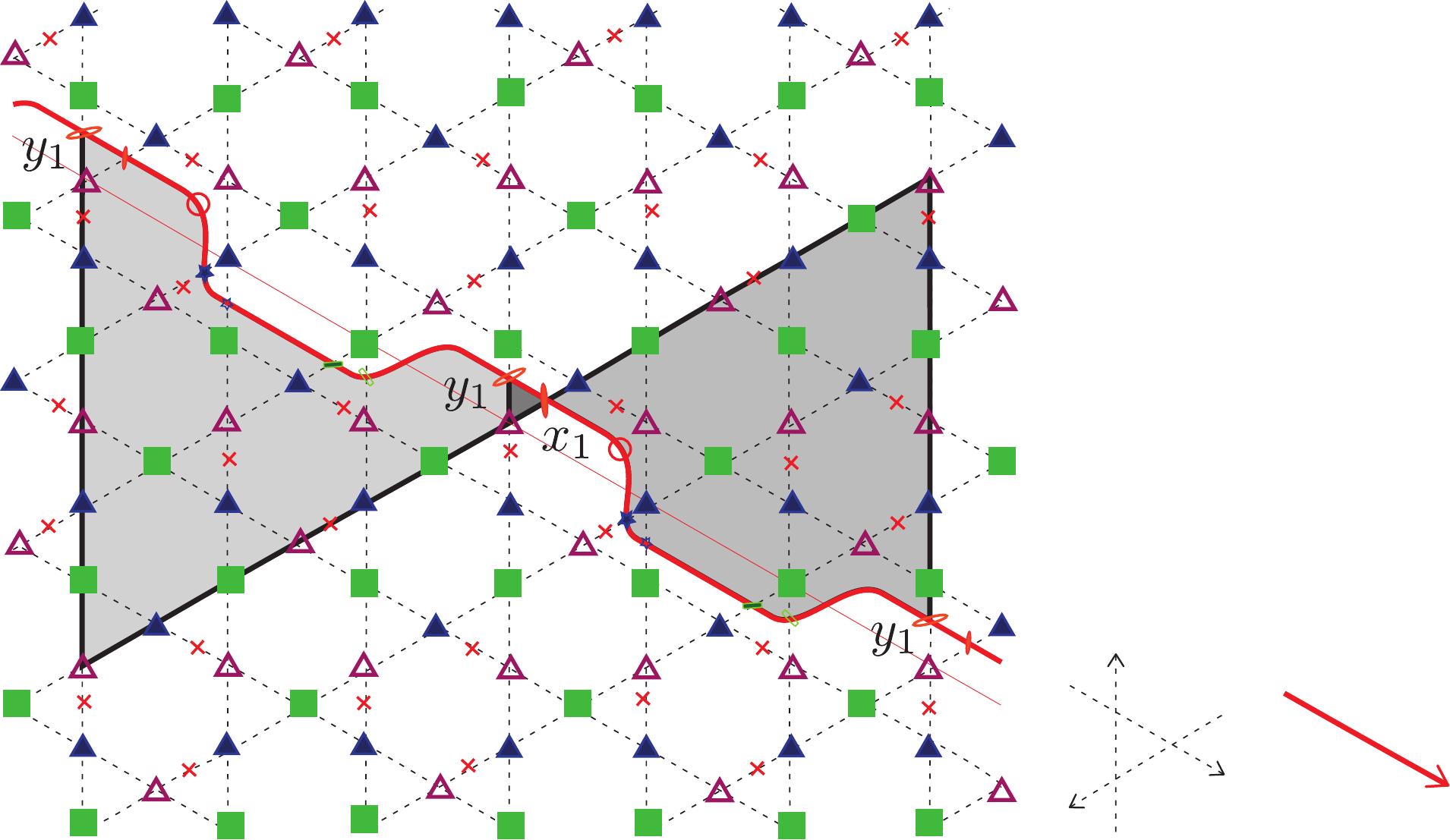}
\caption{Computation of $\alpha_1$}\label{fig:alpha1cal} 
\end{center}
\end{figure}
One type of strips  are  given by triangles, whose vertices
are $x_1, y_1, X$ in a counter-clockwise order (see  the left side of Figure \ref{fig:alpha1cal}). The weighted signed sum of these strips
(in the limit  $\epsilon \to 0$) is given as
$1 - q_\alpha^{6^2} + q_\alpha^{9^2} + \cdots.$
%$$1 - q_\alpha^{so(11)} + q_\alpha^{so(17)} + \cdots.$$
The other type of strips are given by triangles, whose vertices
are $x_1, X, y_1$ (the right side of Figure \ref{fig:alpha1cal}). Their sum gives 
$ -q_\alpha^{6^2} + q_\alpha^{9^2} + \cdots.$
%$$ -q_\alpha^{so(11)} + q_\alpha^{so(17)} + \cdots.$$
Summing up contributions from two types of triangles, we have
$\alpha_1 = 1 - 2q_\alpha^{6^2} + 2 q_\alpha^{9^2} + \cdots. $
%$$\alpha_1 = 1 - 2q_\alpha^{so(11)} + 2 q_\alpha^{so(17)} + \cdots. $$
The other terms can be also computed in a similar way, and we omit the details.

%{\color{red}
%\begin{remark}
%The long diagonal is equipped with the non-trivial spin structure. I'm not sure $\alpha_i$ is computed with the correct sign. (Please correct it)
%\end{remark}
%}

\begin{figure}[htp]
\begin{center}
\includegraphics[height=2in]{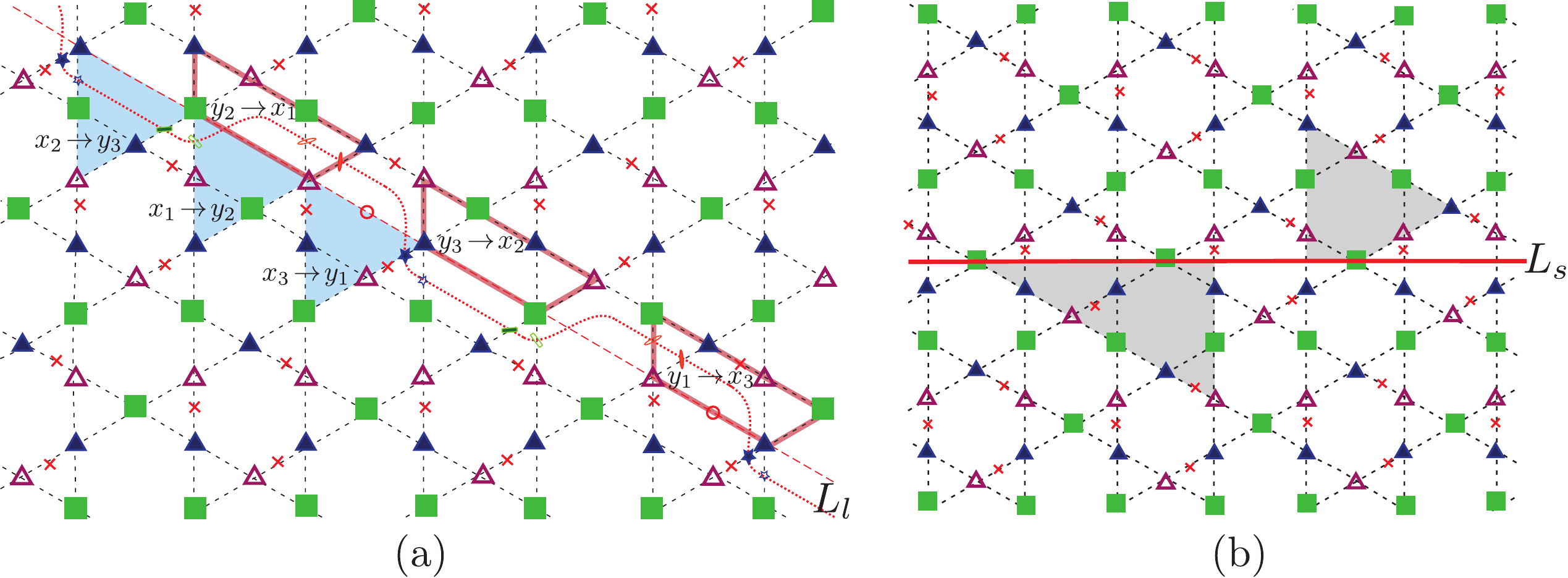}
\caption{Long and short diagonal branes and holomorphic strips between $(\bar{L},b)$ and diagonal branes}
%Translation symmetry among strips (in the limit $\epsilon \to 0$)}
\label{symlongshortdiag} 
\end{center}
\end{figure}

The symmetry appearing in coefficients $\alpha_i$, $\beta_i$, $\gamma_i$ for both matrices \eqref{longdiagJ}, \eqref{longdiagE} can be understood as a translation symmetry of $\bar{L}$ together with the fact that the image of long diagonal in $\mathbb{P}^1_{3,3,3}$ is the horizontal equator. See (a) of Figure \ref{symlongshortdiag}. Three trapezoids above the long diagonal represent (limits of) strips from $y_2$ to $x_1$, from $y_3$ to $x_2$, from $y_1$ to $x_3$, respectively. Coefficients of corresponding entries in the expression \eqref{longdiagJ} all agree up to sign. Indeed, one can check that $(-)$-sign appears whenever there are odd number of $y$ in the monomial for both matrices $J$ and $E$. (We use the same sign rule as in Subsection \ref{FloerSeidel} by fixing a generic point ${\color{red} \circ}$ on $\bar{L}_l$ for its
non-trivial spin structure.)

Similarly, three shaded triangles below the long diagonal with the same area represent strips from $x_3$ to $y_1$, $x_1$ to $y_2$, $x_2$ to $y_3$, whose corresponding entries in the expression \eqref{longdiagE} have the same coefficients, but again $(-)$-sign for $y$. 

The negative signs that appear in $E$ and $J$ can be explained as follows.
Our potential for $\bar{L}$ is of the form $W_{3,3,3} = \phi (x^3 -y^3 + z^3 ) +\psi xyz$. Note that if we make a coordinate change for $\mathcal{M}_{weak}^1 (\bar{L})$ by $\tilde{y} = -y$, then our potential is written in terms of the new coordinate as
$$\tilde{W}_{3,3,3} = \phi (x^3 + {\tilde{y}}^3 + z^3 ) - \psi x \tilde{y} z.$$
With these new coordinates,  the matrix factorization $E$ and $J$ becomes
\begin{equation}\label{nelongdiwagJ}
\tilde{J}=
\bordermatrix{
       &  x_1   	&  x_2  & x_3  \cr
y_1 &\beta_1 x^2+\gamma_1 \tilde{y}z &\beta_3 z^2+\gamma_3 x\tilde{y} &\beta_2 \tilde{y}^2+\gamma_2 zx  \cr
\tilde{y}_2 &\beta_2 z^2 +\gamma_2 x\tilde{y} &\beta_1 \tilde{y}^2+\gamma_1 zx &\beta_3 x^2+\gamma_3 \tilde{y}z   \cr
y_3 &\beta_3 \tilde{y}^2+\gamma_3 zx &\beta_2 x^2+\gamma_2 \tilde{y}z &\beta_1 z^2+\gamma_1 x\tilde{y}   \cr}
\end{equation}

\begin{equation}\label{newlongdiagE}
\tilde{E}=
\bordermatrix{
       &  y_1   		 &  y_2             &  y_3               \cr
x_1 &  \alpha_1 x    &  \alpha_2 z  &  \alpha_3 \tilde{y}    \cr
x_2 &  \alpha_3 z    &  \alpha_1 \tilde{y}  &  \alpha_2 x    \cr
x_3 &  \alpha_2 \tilde{y}    &  \alpha_3 x  &  \alpha_1 z    \cr}.
\end{equation}

As $\tilde{J}, \tilde{E}$ give a matrix factorization of $\tilde{W}_{3,3,3}$, we have
\begin{equation}\label{eq:matjecom}
\tilde{J} \tilde{E} = \tilde{E} \tilde{J} = \tilde{W}_{3,3,3} \cdot Id.
\end{equation}
Hence we can write $\beta, \gamma$ in terms of $\alpha$.
\begin{eqnarray*}
\beta_1 = \phi / \alpha_1,& \beta_2 =  \phi / \alpha_2,& \beta_3 = \phi / \alpha_3 \\
\gamma_1 = -\frac{\alpha_1 \beta_3}{\alpha_2} = -\frac{\phi \alpha_1}{\alpha_2 \alpha_3},& \gamma_2 = -\dfrac{\alpha_2 \beta_1}{\alpha_3} = -\dfrac{\phi \alpha_2}{\alpha_3 \alpha_1},&\gamma_3 = -\frac{\alpha_3 \beta_2}{\alpha_1} = -\frac{\phi \alpha_3}{\alpha_1 \alpha_2}
\end{eqnarray*}
We remark that $\alpha_i$'s are invertible series even in $\Lambda^{nov}_0$ because of minimal strips shown in Figure \ref{longdiag} whose area degenerate to $q^0 = 1$ in the limit $\epsilon \to 0$.

Moreover, Equation \eqref{eq:matjecom} implies that 
$\alpha_1 \gamma_1 + \alpha_2 \gamma_2 + \alpha_3 \gamma_3 =  \psi.$
Combining with the above, we get $\phi (\alpha_1^3 + \alpha_2^3 + \alpha_3^3) / \alpha_1 \alpha_2 \alpha_3 = \psi$, or equivalently
$\phi ( \alpha_1^3 + \alpha_2^3 + \alpha_3^3)  - \psi \alpha_1 \alpha_2 \alpha_3 =0.$
This proves that
\begin{prop}
$(\alpha_1, \alpha_2, \alpha_3)$ corresponds to a point in the elliptic curve $\tilde{W}_{3,3,3} = 0$. 
\end{prop}

\begin{remark}
In \cite{B}, the physicists studied the endomorphisms of matrix factorization category of the potential 
%\begin{equation}\label{phyW}
$W = \frac{1}{3}x^3 + \frac{1}{3} y^3 + \frac{1}{3} z^3 - a xyz.$
%\end{equation}
A matrix factorization of $W$ similar to $(E,J)$ was defined, and
they derived from physical arguments that the product structure of endomorphisms of the matrix factorization corresponds to counting of triangles with boundary on (lifts of) long diagonals $L_l$ (without explicitly computing the count).
\end{remark}

\subsubsection{Matrix factorization for a short diagonal $\bar{L}_s$}
The mirror matrix factorization of the short diagonal $L_s$ can be computed in a similar way.
Types of holomorphic strips that contributes to the matrix factorizations are drawn in (b) of Figure \ref{symlongshortdiag}.

%\begin{figure}[htp]
%\begin{center}
%\includegraphics[scale=0.6]{MFshortdiag}
%\caption{The short diagonal brane $L_s$ and holomorphic strips bounding $(\bar{L},b)$ and $L_s$}\label{MFshortdiag} 
%\end{center}
%\end{figure}
The resulting matrix factorization is of the following form and
we leave the detailed check to the reader.
\begin{equation}\label{MFshort}
P_{L_s}:=
\left(\begin{array}{cccc} &  & L_1  & Q_2 \\
 &  & -L_2 & Q_1 \\
 Q_1 & -Q_2 &  &  \\
 L_2 &  L_1 &  & \end{array}\right)
\end{equation}
where linear terms $L_1$, $L_2$ and quadratic terms $Q_1$, $Q_2$ are given by
$$L_1 = \alpha_1 x + \alpha_2 y + \alpha_3 z $$
$$L_2 = \alpha_3 x + \alpha_2 y + \alpha_1 z $$
$$Q_1 = ( \beta_1 x^2 + \gamma_1 yz )+( \beta_2 y^2 + \gamma_2 zx ) +( \beta_3 z^2 + \gamma_3 xy ) $$
$$Q_2 = ( \beta_3 x^2 + \gamma_3 yz )+( \beta_2 y^2 + \gamma_2 zx ) +( \beta_1 z^2 + \gamma_1 xy )$$
and $\alpha_i$, $\beta_i$, and $\gamma_i$ ($i=1,2,3$) are power series in $q$. There are several relations among these power series which come from the equation $P_{L_s}^2 = W_{3,3,3} \cdot Id$, which we also omit.

\begin{remark}
The matrix factorization takes different expression from that in \cite{B} corresponding to the short diagonal.   %but we do not know if the above is equivalent to their matrix factorization.
\end{remark}

\subsubsection{Magical cancellation}\label{ssubsec:magiccancel}
Showing the identity $(m_1^{b,0})^2 = W$ directly by computation is not easy, since
it involves some `magical' cancellations.
Before we show an example of such cancellations, let us review the proof of this equation.
This is done by investigating the moduli space of holomorphic strips of Maslov-Viterbo index two with $\conste^b$ insertions on the upper boundary. Its codimension one boundary components  are
given by either a broken strip consisting of two strips of Maslov-Viterbo index one or a strip of Maslov-Viterbo index zero together with a Maslov index two disc bubble.
%\begin{figure}[htp]
%\begin{center}
%\includegraphics[scale=0.5]{cancelpairlongd}
%\caption{A magical canceling pair}\label{fig:cancelpairlongd} 
%\end{center}
%\end{figure}
\begin{figure}[htp]
\begin{center}
\includegraphics[scale=0.4]{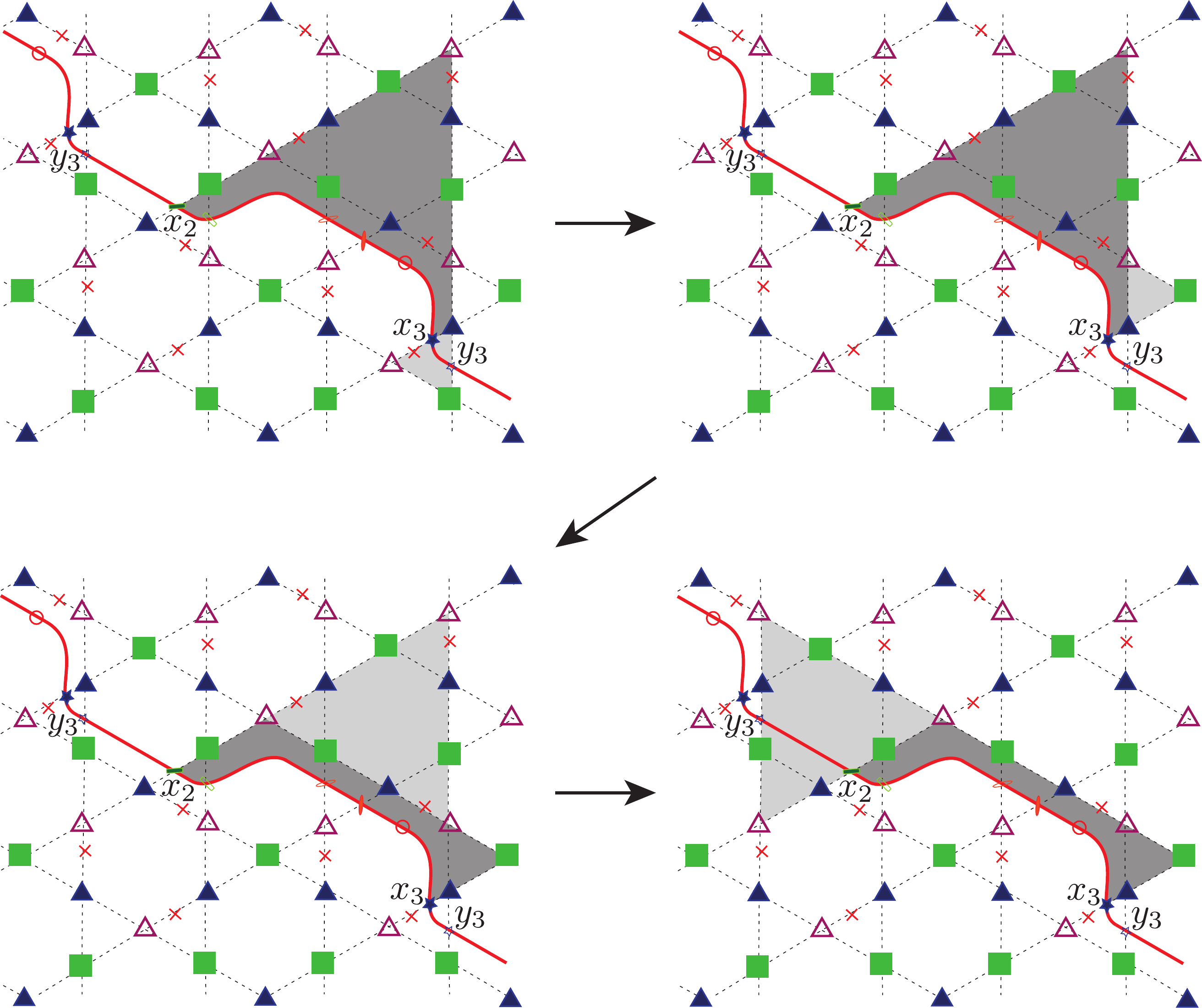}
\caption{A magical cancellation}\label{fig:cancelpairlongd1} 
\end{center}
\end{figure}
In the case of the Seidel Lagrangian, the latter can be divided into two types depending on
whether the nodal point (connecting the disc bubble and the strip)  is immersed or not.
%{\color{red}(?) INSTEAD whether the starting point and the end point of a broken strip coincide or not?}

Non-trivial strips of Maslov-Viterbo index zero exist when the nodal point is immersed.
Such configurations with an immersed nodal point should cancel among themselves
since $b$ is a bounding chain and $m(\conste^b)$ is a multiple of $e$,
and any immersed output contributing to $m(\conste^b)$ in fact cancels out.
In the case that nodal point is not immersed, one can note that the only
non-trivial contribution is from a constant strip, since $m(\conste^b)$ is a multiple of $e$.
This provides the term $W$.

%{\color{red} REVISE : Note that $(m_1^{b,0})^2$ is zero up to holomorphic discs bounding $\bar{L}$ since the long diagonal is unobstructed. Thus, as in Morse theory, broken strips from a intersection point to another intersection point always exist in pair so that they eventually all cancel out. Otherwise, the end points of such strips should appear in the right hand side of the equation $(m_1^{b,0})^2=W$. 
%
%If a broken strip has the same starting and end points, one can see that the only
%non-trivial contribution comes from a constant strip, since $m(\conste^b)$ is a multiple of $e$.
%This provides the term $W$.}

An example of a magical cancellation pair is shown in Figure \ref{fig:cancelpairlongd1}. 
% As in Figure \ref{fig:cancelpairlongd1}, 
These cancellation are connected through the cancellation
pairs of immersed outputs of $m(\conste^b)$ given by reflection along equator.
Hence we are using the fact that all immersed outputs of $m(\conste^b)$ cancel out.

%\subsubsection{Matrix factorization for $\oi$}
%We will  consider the matrix factorization corresponding to $\oi$. 
%The  matrix factorization of $\oi$  is defined as the Floer chain complex  $CF( (\bar{L}, b), \bar{L})$ where $b$ varies in $\mathcal{MC}_{weak}^1 (\bar{L}) \cong \C^3$.
%Note that Floer complex  $CF( (\bar{L}, b), \bar{L})$ has 4 odd generators $p$, $X$, $Y$, $Z$ (which are
%the maximum of the Morse function, and odd immersed generators) and 4 even generators 
% $e$, $\bar{X}$, $\bar{Y}$, $\bar{Z}$ (which are the minimum of the Morse function, and even immersed generators). Here, we are using a sort of Morse-Bott model for $CF( (\bar{L}, b), \bar{L})$.
%
%{\color{red} Before we discuss the full matrix factorization for $W_{(3,3,3)}$, let us consider the minimal potential $W_{(3,3,3)}^{min}$
%and its matrix factorization from the counting of holomorphic strips. WILL BE WRITTEN SOON}

\section{Orbifold sphere $\mathbb{P}^1_{a,b,c}$} \label{sect_dim1}
In this section we apply our mirror construction to the orbifold $\mathbb{P}^1_{a,b,c}$ for all $a,b,c$.  We derive homological mirror symmetry using our functor when $\frac{1}{a}+\frac{1}{b}+\frac{1}{c}\leq 1$.  We use the Seidel Lagrangian $\bar{L}$ depicted in Figure \ref{SeiLagr} for our construction.

The essential ingredients are the following.  Solving the weak Maurer-Cartan equation is the first key step.  We will show that $\bar{L}$ and its odd-degree immersed deformations are weakly unobstructed (Theorem \ref{thm:weak}).  This allows us to construct the mirror of $\mathbb{P}^1_{a,b,c}$ and the functor.  

Second we need to show that the criterions of Theorem \ref{thm:criterion_equiv} are satisfied, and hence our functor derives homological mirror symmetry.  The Seidel Lagrangian split-generates by the generating criterion of Abouzaid applied to a manifold cover of $\mathbb{P}^1_{a,b,c}$ where the Seidel Lagrangian $\bar{L}$ lifts to an embedded smooth curve.  Such a cover always exists for $\frac{1}{a}+\frac{1}{b}+\frac{1}{c}\leq 1$ from the joint work with the authors with Kim \cite{CHKL14} (See Proposition \ref{prop:cover}).  Since the Fukaya category of $\mathbb{P}^1_{a,b,c}$ is defined as the invariant part of that of the cover, it is split generated by $\bar{L}$.

Third we compute the image matrix factorization of $\bar{L}$ under our functor, and show that it split-generates $D^\pi \MF(W)$ (Corollary \ref{cor:cptgen1}) by using the result of Dyckerhoff \cite{Dy}.  We also verify that our functor induces an isomorphism on $HF(\bar{L},\bar{L})$ (Section \ref{subsec:mirrormorlev}).  This proves homological mirror symmetry for $\mathbb{P}^1_{a,b,c}$.

For the construction of mirror $W$, we will work for all $a,b,c > 1$.  However for homological mirror symmetry, we will restrict to the case $\frac{1}{a} + \frac{1}{b} + \frac{1}{c} \leq 1$.  The reason is that the mirror $W$ has more than one critical points for the spherical case $\frac{1}{a} + \frac{1}{b} + \frac{1}{c} > 1$ and exhibits a rather different behavior.  While we believe our functor still derives homological mirror symmetry, we will discuss this case separately elsewhere.

For our purpose here we do not need the explicit expression of $W$, and hence we will only compute the leading-order terms.  $W$ takes the form
$a(q)x^a + b(q)y^b + c(q)z^c + \sigma(q) xyz + \ldots $
which is a quantum-corrected superpotential and is different from the one used by Seidel \cite{Se}.  The explicit computation of the full superpotential $W$ will be given in \cite{CHKL14}.

\subsection{Manifold covers of $\mathbb{P}^1_{a,b,c}$}
The orbifold $\mathbb{P}^1_{a,b,c}$ has a branched covering $\Sigma$ which is a Riemann surface. We recall some relevant facts.
Orbifolds have an analogous covering space theory due to Thurston \cite[Chapter 13]{Thu79}, and the following orbifold Euler characteristic plays an important role in characterizing coverings of orbifolds.

If an orbifold $Q$ has a $CW$-complex structure so that the local group is constant on each open cell $c$ (denote by  $G(c)$ the local group on $c$), then its orbifold Euler characteristic is defined as
$ \chi^{orb} (Q) := \sum_{{\rm cells} \,\, c} \frac{ (-1)^{\dim c} }{ |G(c)| }.$
For $\mathbb{P}^1_{a,b,c}$, we can divide it into two cells, namely the upper and the lower hemispheres, by an equator which passes through the three orbifold points $p_a, p_b, p_c$.  This defines a $CW$-complex structure of $\mathbb{P}^1_{a,b,c}$ satisfying the property above.
Therefore
$$\chi^{orb} (\mathbb{P}^1_{a,b,c} ) = \frac{1}{a} + \frac{1}{b} + \frac{1}{c}  -1.$$

In addition, we consider the reflection $\tau$ along the equator (1-cells) of $\mathbb{P}^1_{a,b,c}$, and let $D_{a,b,c}$ denote the quotient of
$\mathbb{P}^1_{a,b,c}$ by the reflection. Then we have
$\chi^{orb}  (D_{a,b,c}) =\frac{1}{2} \left( \frac{1}{a} + \frac{1}{b} + \frac{1}{c}  -1 \right).
$

Let $M$ be the universal cover of $\mathbb{P}^1_{a,b,c}$.
According to whether $\frac{1}{a} + \frac{1}{b} + \frac{1}{c} > 1$ or $=1$ or $<1$, the universal cover $M$ is $S^2$ or $\R^2$ or $\mathbb{H}^2$. We may assume that the deck transformation group $\pi^1_{orb}  \left(\mathbb{P}^1_{a,b,c} \right)$ acts on $M$ as isometries.

It is well-known (see \cite{D}) that the orbifold fundamental groups of $\mathbb{P}^1_{a,b,c}$ and $D_{a,b,c}$ are
$$\pi_1^{orb} (\mathbb{P}^1_{a,b,c}) = <\rho_a, \rho_b, \rho_c \,\, | \, \left(\rho_a \right)^a = \left(\rho_b \right)^b = \left(\rho_c \right)^c = \rho_c \rho_b \rho_a  =1 >,$$
and 
%\begin{equation}\label{eq:pi1dabc}
$\pi_1^{orb} (D_{a,b,c}) = <\rho_a, \rho_b, \rho_c, \tau_{ab}, \tau_{bc}, \tau_{ca} \,\, |\,\, R>$
%\end{equation}
respectively where the relation set $R$ is generated by
\begin{equation*}
\left\{
\begin{array}{l}
\left(\rho_a \right)^a = \left(\rho_b \right)^b = \left(\rho_c \right)^c = \rho_c \rho_b \rho_a =1\\
\rho_a = \tau_{ca} \tau_{ab} \quad \rho_b = \tau_{ab} \tau_{bc} \quad \rho_c = \tau_{bc} \tau_{ca} \\
\left( \tau_{ab} \right)^2 = \left( \tau_{bc} \right)^2 = \left( \tau_{ca} \right)^2 =1.
\end{array} \right.
\end{equation*}
One can easily see that $\pi_1^{orb} (\mathbb{P}^1_{a,b,c})$ injects to $\pi_1^{orb} (D_{a,b,c})$.

 It is also known that (see \cite[Theorem 2.5]{Sc}) any finitely generated discrete subgroup of the isometry groups of $S^2$, $\R^2$ or $\mathbb{H}^2$ with a compact quotient space has a torsion free subgroup, say $\Gamma$, of finite index. 
Fix a torsion free subgroup  $\Gamma$ of $\pi^1_{orb}  \left(\mathbb{P}^1_{a,b,c} \right)$.  The quotient $ M/\Gamma$ is a compact surface $\Sigma$ with genus $g$ which is a covering space (in the sense of
orbifold covering) of $\mathbb{P}^1_{a,b,c}$. ($M/\Gamma$ is a manifold since $\Gamma$ is torsion-free.)

\begin{remark}\label{orbeuler}
Torsion free subgroup $\Gamma$ is not unique in general.
Suppose that $\Sigma \to \mathbb{P}^1_{a,b,c}$ is a $d$-fold branched covering. 
Then $d$ is always a common multiple of $a, b, c$ and satisfies the relation
$\chi (\Sigma) = 2-2g = d \left(  \frac{1}{a} + \frac{1}{b} + \frac{1}{c}  -1  \right)=d\cdot \chi (\mathbb{P}^1_{a,b,c}).$
\end{remark}
 From covering space theory, 
the deck transformation group $G$ of the branched covering $\Sigma \to \mathbb{P}^1_{a,b,c}$ is
$\pi_1^{orb} \left( \mathbb{P}^1_{a,b,c} \right) / N(\Gamma)$, where $N(\Gamma)$ is the normalizer of $\Gamma$. Note that $G$ is finite. It is often convenient to work on $\Sigma$ $G$-equivariantly instead of working on the orbifold
$\mathbb{P}^1_{a,b,c}$ directly. 

We also give a brief remark on the  universal cover $M \to 
\mathbb{P}^1_{a,b,c}$. 
Endow $M$ with a triangulation by lifting the triangulation of $D_{a,b,c}$.  Since the orbifold fundamental group of $D_{a,b,c}$ acts on its universal cover, we have reflections about each 1-skeleton of the triangulation.  

\begin{remark}
One can show that $\pi_1^{orb} \left(D_{a,b,c} \right) \cong \Z / 2\Z \ltimes \pi_1^{orb} \left(\mathbb{P}^1_{a,b,c} \right)$. %where $\Z /2 \Z$ is generated by a reflection.
\end{remark}

Now we consider the lifts of $\bar{L}$ to a branched covering $\pi:\Sigma \to  \mathbb{P}^1_{a,b,c}$.
%First of all, $\bar{L}$ is an immersion $\iota:S^1 \to \mathbb{P}^1_{a,b,c}$, and  a generator of $\pi_1(S^1)$ maps to $\rho_a\rho_b\rho_c$ which is nontrivial in general. (Recall $\rho_c\rho_b\rho_a=1$.) If $\rho_a\rho_b\rho_c \in \Gamma$, then $\bar{L}$ lifts to an embedding of a circle. This happens for example when the deck transformation group $G$ is abelian. In this case, we may take the $G$-action to one of liftings to get an equivariant immersion into $\Sigma$ from $|G|$-copies of $S^1$.  For the case that $\rho_a\rho_b\rho_c \notin \Gamma$, such lifting does not exist.  Instead, by pre-composing $\iota$ with degree $k$ map $S^1 \to S^1$, we obtain another immersion $\iota':S^1 \to \mathbb{P}^1_{a,b,c}$, corresponding to $(\rho_a \rho_b \rho_c)^k$.  We may take minimal  $k >1$ such that $(\rho_a \rho_b \rho_c)^k \in \Gamma$ since $[\pi_1^{orb} (\mathbb{P}^1_{a,b,c}) : \Gamma] < \infty$.  Hence, there exist a lift $\WT{\iota'}:S^1 \to \Sigma$ of $\iota'$, which is a Lagrangian immersion to $\Sigma$.
In a joint work with Kim \cite[Proposition 2.4]{CHKL14}, we prove that in the hyperbolic case $\frac{1}{a}+\frac{1}{b}+\frac{1}{c}<1$, there exists a compact manifold cover $\Sigma$ of $\mathbb{P}^1_{a,b,c}$  such that the Seidel Lagrangian lifts to a Lagrangian embedding into $\Sigma$.
For an elliptic curve quotient ($\frac{1}{a}+\frac{1}{b}+\frac{1}{c}=1$), by a case-by-case study (where there are only three cases) we can see that the Seidel Lagrangian lifts to an embedded curve in the torus.  When $\frac{1}{a}+\frac{1}{b}+\frac{1}{c}>1$, namely the spherical case, other than $(a,b,c) = (2,2,\textrm{odd})$, we can show that the Seidel Lagrangian lifts to a circle (which is an equator) in the sphere \cite[Lemma 12.1]{CHKL14}.  However for $(a,b,c) = (2,2,\textrm{odd})$, the Seidel Lagrangian lifts to an immersed curve (which is not embedded) in the sphere.  Hence we conclude the following.

\begin{prop}[Proposition 2.4 and Lemma 12.1 of \cite{CHKL14}] \label{prop:cover}
For $\frac{1}{a}+\frac{1}{b}+\frac{1}{c} \leq 1$, there exists a manifold cover of $\bP^1_{a,b,c}$ such that the Seidel Lagrangian lifts as an embedded curve.
\end{prop}

\subsection{Floer generators of the Seidel Lagrangian}\label{FloerSeidel}
We take a manifold cover of $\mathbb{P}^1_{a,b,c}$ in order to define its Fukaya category (Section \ref{sec:orbisurface}).  We will focus on the case $\frac{1}{a}+\frac{1}{b}+\frac{1}{c} \leq 1$.  We take a manifold cover $\Sigma$ of $\mathbb{P}^1_{a,b,c}$ such that the Seidel Lagrangin lifts as an embedded curve in $\Sigma$.  The objects of the Fukaya category are $G$-equivariant collections of embedded curves in $\Sigma$ avoiding orbifold points.

The Lagrangian Floer complex $CF(\bar{L},\bar{L})$ is $\Z/2$-graded and has six immersed generators (two for each immersed points), namely $X$, $Y$, $Z$ of degree 1 and $\bar{X}, \bar{Y}, \bar{Z}$ of degree 2.  It also has generators coming from $H^\ast (S^1)$ where $S^1$ is the domain of the Lagrangian immersion.  
%cohomology classes of $S^1$ which is the domain of the Lagrangian immersion $S^1 \to \mathbb{P}^1_{a,b,c}$ onto $\bar{L}$.  The Floer complex is $\Z/2$-graded, and the parity of the immersed generator is determined by the local intersection number at the corresponding immersed point. (See Figure \ref{figdeg}.).  We will also follow Seidel \cite{Se} to consider certain integer grading (which is not strictly compatible with $\AI$-structure). To simplify the notation, we will use the convention $\overline{\left(\bar{A} \right)} =A$.
We use the Morse complex $CM(f)$ instead of $H^\ast(S^1)$ for some Morse function $f$ on the domain of the immersion $\bar{L}$, which has two critical points $e$ and $p$, namely the minimum and the maximum points of $f$.  Thus $CF(\bar{L},\bar{L})$ is spanned by $X$,$Y$,$Z$, $\bar{X}$, $\bar{Y}$, $\bar{Z}$, $e$ and $p$.

For convenience we may use the $1/3$ grading defined by Seidel \cite{Se}.  Note that in general the $1/3$-grading is not compatible with the orbifold structure of $\bP^1_{a,b,c}$ and does not really define a $\Z$-grading on the Fukaya category.  We will only work over $\Z_2$ grading and hence will only care about the parity of the grading.  It is defined by taking a section $\eta^3$ of $\left(T^\ast S^2\right)^{\otimes 3}$ which has exactly three double poles (and nowhere-vanishing elsewhere), which are identified with the three orbifold points of $\bP^1_{a,b,c}$.  Then we get a map
$S^1 \to S^1$ given by $ \theta \mapsto \frac{\eta^3 \left(\iota_\theta^{\otimes 3} \right) }{ || \eta^3 \left(\iota_\theta^{\otimes 3} \right) || } $
where the first $S^1$ is the domain of the Lagrangian immersion $\iota : S^1 \to \bar{L}$ and $\iota_\theta = \iota_\ast \left( \frac{\partial}{\partial \theta} \right)$.  This map can be lifted to $a : S^1 \to \R$ which is the so called $1/3$-grading.

Under the $1/3$-grading, $e$ has degree $0$, $X,Y,Z$ have degree $1$, $\bar{X},\bar{Y},\bar{Z}$ have degree $2$, and $p$ have degree $3$.  Thus the odd generators are $X,Y,Z,p$ and the even generators are $e,\bar{X},\bar{Y},\bar{Z}$.

%Consider a uniformizing chart around the singular point $p_a$ of order $a$, which is given by the local coordinate change $z = w^a$. Then, $z^{-2} d z^3=  a^3 w^{a - 3 } dw^3$ implies that the induced section $\eta^3$ of $\left(T^\ast \mathbb{P}^1_{a,b,c}\right)^{\otimes 3}$ from $\tilde{\eta}^3$ has vanishing order $a-3$ at $p_a$. Likewise, ${\rm ord\,}(\eta^3, p_b) = b - 3$ and ${\rm ord\,}(\eta^3, p_c) = c - 3$. 

%Without loss of generality let's assume that $\bar{L}$ intersects itself with the angle $\theta = \pi / 3$ at each self intersection points. Grading of each generator is given as $i (e) = 0 , i(q) = 3$ and 
%$$i ( X) = i(Y) = i(Z) = \frac{ 3 \theta  }{ \pi} = 1, i (\bar{X} ) = i ( \bar{Y}) = i (\bar{Z}) = \frac {3 (\pi - \theta) }{\pi} = 2.$$

%With respect to $i$, each $A_\infty$-operation can be decomposed as $m_k = m_k^0 + m_k^{higher}$ where $m_k^{higher}$ comes from polygons which meet one of singular points in positive degree. ($m_k^0$ counts the contributions from polygons in $\mathbb{P}^1_{a,b,c} \setminus \{p_a, p_b, p_c\}$.) From \cite{Se}, a holomorphic disc $u$ contributes to the component of $m_k^{higher}$ with degree 
%\begin{equation}\label{eq:grai}
%6-3k +  2(a-3) {\rm deg}\, (u,p_a) + 2(b-3) {\rm deg}\, (u,p_b) + 2(c-3) {\rm deg}\, (u,p_c).
%\end{equation}
%For $\mathbb{P}^1_{3,3,3}$,  we have ${\rm deg}\, m_k \equiv 6-3k$, and hence we have a $\Z$-grading on our $A_\infty$-algebra.

\subsection{Weak bounding cochains of  the Seidel Lagrangian}
In this section we solve the Maurer-Cartan equation for weakly unobstructed immersed deformations of $\bar{L} \subset \bP^1_{a,b,c}$.  We prove that the linear combinations of all degree $1$ elements, namely $X$, $Y$, $Z$, are weak bounding cochains.

A key ingredient is the anti-symplectic involution $\tau$ on $\bP^1_{a,b,c}$ given by the reflection about the equator which passes through the three orbifold points.  The symplectic structure $\omega$ is taken such that $\tau^*\omega = -\omega$.  Moreover $\bar{L}$ is invariant under $\tau$, namely it is symmetric along the horizontal equator of $\mathbb{P}^1_{a,b,c}$.  The immersed points of $\bar{L}$ lies in the equator which is the fixed locus of $\tau$.  This symmetry will be essential to prove that $\bar{L}$ is weakly unobstructed. 
 
 \begin{lemma}\label{signcancel}
 Suppose that $A_0$,$\cdots$, $A_k$ are all degree 1 immersed generators.  Consider a holomorphic polygon
 $P \in \mathcal{M} (A_0, A_1, \cdots, A_k),$
 and its reflection by $\tau$, $P^{op} \in \mathcal{M} (A_{k-1},A_{k-2} \cdots,A_0, A_k).$
Then the contribution of $P$ to $\left< m_k (A_0, \cdots, A_{k-1}), \bar{A_k} \right>$, which denotes the coefficient of $\bar{A_k}$ in $m_k (A_0, \cdots, A_{k-1})$, cancels with the contribution of $P^{op}$ to $\left<m_k (A_{k-1}, \cdots, A_0), \bar{A_k} \right>$.
\end{lemma}

\begin{figure}[htp]
\begin{center}
\includegraphics[scale=0.4]{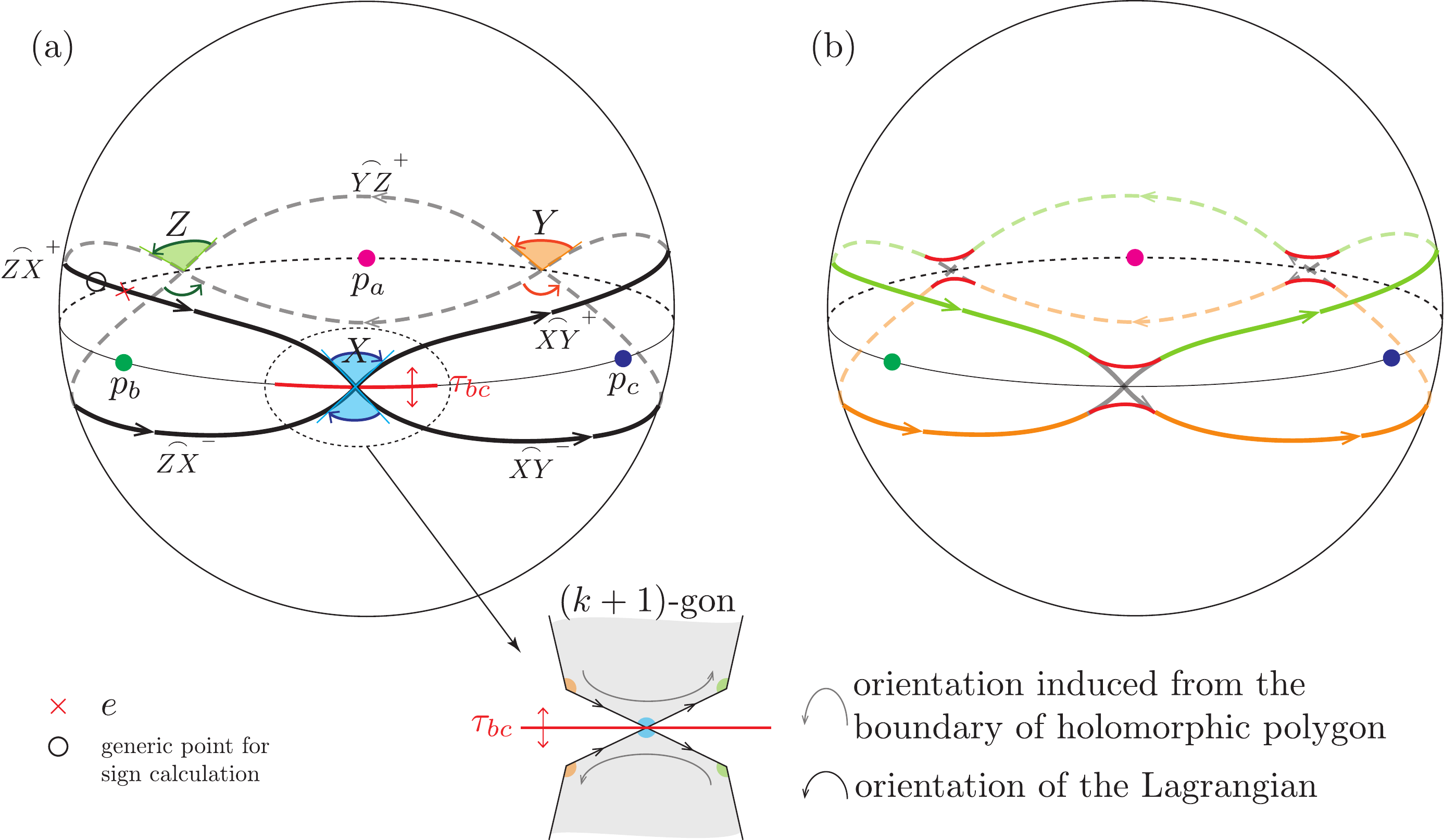}
\caption{(a) The Seidel Lagrangian and (b) its smoothing}\label{SeiLagr}
\end{center}
\end{figure}

%$\left< \, , \, \right>$ is  {\em not} the Poincar\`{e} paring.
 \begin{proof}
Recall that the signs of polygons are determined by two factors. The first is the comparison of
counter-clockwise orientation along the boundary of the polygon with that of the bounding Lagrangians.  The other is the number
times that the boundary of a polygon passes through the generic point which represents the non-trivial spin structure.

Regarding the first factor, observe that the reflection $\tau$ reverses the orientation along the boundary of a holomorphic polygon,
but $\tau$ preserves the orientation of $\bar{L}$.
Hence the sign difference between $P$ and $P^{op}$ from the first factor is given as 
\begin{equation}\label{eq:signfromori}
(-1)^{|A_0| + |A_1| + \cdots + |\bar{A_k}|} = (-1)^{k}.
\end{equation}

Now let us compare the second factor, which is  the number of times that $P$ and $P^{op}$ pass through the given generic point.
Let us consider the union of boundaries $\partial P \cup \partial P^{op}$, and the Seidel Lagrangian $\bar{L}$.
Define $\arc{XY}^\pm$, $\arc{YZ}^\pm$, $\arc{ZX}^\pm$ as in Figure \ref{SeiLagr} and call them {\em minimal arcs} of the Seidel Lagrangian. 
 Suppose that the generic point ``$\circ$" lies on the minimal arc $\arc{ZX}^+$ (see Figure \ref{SeiLagr}). 

Here are a few observations.
The first observation is that  $\partial P \cup \partial P^{op}$ evenly covers $\bar{L}$. (i.e. $\partial P \cup \partial P^{op}$ is a cycle that is an inter multiple of $[\bar{L}]$.) To see this, consider the point $p$
which travels along $\partial P$, and its reflection image $\tau(p)$, which travels along   $\partial P^{op}$.
The pair $( p, \tau(p))$ always travel on  $\pm$( or $\mp$) pair of minimal edges, and ends at the original vertex 
after full rotation along $\partial P$. Hence the number of times that $\partial P \cup \partial P^{op}$ covers a minimal arc is the same for all minimal arcs.

The second observation is that any edge of a holomorphic polygon $P$ for the potential consists of odd number of minimal arcs $\bar{L}$. (In other words, the corners of a holomorphic polygon always lie
on one of the hemisphere.) In fact, if two odd degree corners of such a polygon were connected by an edge consisting of even number of minimal arcs, then one of the odd corner would be non-convex.

Now consider a $(k+1)$-gon $P$ and suppose that  $\partial P \cup \partial P^{op}$ covers $\bar{L}$ $l$-times (i.e. $[\partial P \cup \partial P^{op}] = l [ \bar{L}]$). Then $\partial P \cup \partial P^{op}$ covers $6l$ minimal arcs, and hence $\partial P$ covers $3l$ minimal arcs. 
We claim that $k+1 \equiv l \mod 2$.
First, if $k+1$ is odd, we see from the second observation that $\partial P$ covers odd number 
of minimal arcs. Since  $\partial P$ covers $3l$ minimal arcs, $l$ should be odd.
Similarly, if $k+1$ is even, $\partial P$ covers even number of minimal edges, and hence $l$ is even.

Therefore,
\begin{equation*}
\begin{array}{l}
\mbox{the number of}\,\, ``\circ" \,\, \mbox{on} \,\, P - \mbox{the number of} \,\,``\circ" \,\, \mbox{on} \,\, P^{op}\\
\equiv \mbox{the number of}\,\, ``\circ" \,\, \mbox{on} \,\, P + \mbox{the number of} \,\,``\circ" \,\, \mbox{on} \,\, P^{op} \\
\equiv l \equiv k+1 \mod 2
\end{array}
\end{equation*}
Combined this with Equation \eqref{eq:signfromori}, the total difference of signs is $(-1)^{2k+1} =-1$ and hence contributions of $P$ and $P^{op}$  cancels each other.
 \end{proof}

Now we are ready to prove the main theorem of this subsection.

\begin{theorem}\label{thm:weak}
$b = x X + y Y + z Z$ is a weak bounding cochain for any $x,y,z \in \C$. 
\end{theorem}
\begin{proof}
Let us first consider the elliptic and hyperbolic cases. Then
$m(\conste^b)$ is a linear combination of the terms
$m_k (X_1, X_2, \cdots, X_k)$
where $k \geq 1$ and $X_i$'s are among $X$, $Y$ and $Z$. It suffices to show that the coefficient of $\bar{X}$ in $m_k (X_1, X_2, \cdots, X_k)$ is canceled with that of $\bar{X}$ in $m_k (X_k, X_{k-1}, \cdots, X_1)$.  Indeed, terms in $m (\conste^b)$ always appear in a pair since if a holomorphic polygon $P$ contributes to $m(\conste^b)$, then its reflection image $P^{op}$ also contributes to $m(\conste^b)$. Notice that this is exactly the situation considered in the previous lemma. i.e. two contributions from $P$ and $P^{op}$ are canceled since all $X_i$'s have odd degrees. Therefore we obtain $m ( \conste^b ) = c e$ by considering its degree.
In the spherical case, the above argument still works, but there is an additional contribution of $m_0$. Namely, 
if $(a,b,c) \neq (2,2,\textrm{odd})$, then lifts of Seidel Lagrangian in the universal cover $S^2$ are topological circles and each lift bisect $S^2$. Hence two discs bounded by such a lift gives  $m_0(1) = c' e$ for some $c' \in \Lambda_0$. In \cite{CHKL14} $c'$ is
explicitly computed, as they are the constant terms of the potential in the spherical cases.
  In the case that $(a,b,c) = (2,2,\textrm{odd})$, then lifts in the universal cover is still immersed, and hence does not bound a disc. But the 1-gon in the universal cover could contribute to $m_0(1)$ with an output given by the corresponding immersed generator. But  Lemma \ref{signcancel} for $k=0$ implies that such contributions from this 1-gon  and its reflection cancel out, and hence we get $m_0(1) =0$ in this case. Hence we also have $m ( \conste^b ) = c e$ in the spherical case.

\end{proof}

%\subsection{The generalized SYZ mirror of $\mathbb{P}^1_{a,b,c}$}

\subsection{The mirror superpotential}
We obtain the superpotential $W(b)$ where $m(\conste^b) = W(x,y,z) e$.  It takes the form
$W : \Lambda_0^3 \to \Lambda_0$ with $W(x,y,z)=\sum a_{mnl}(q) x^m y^n z^l,$
 where ${a_{mnl}}$ is a power series in $q$ which counts areas of $(m+n+l)$-gons with $m$ $X$'s, $n$ $Y$'s and $l$ $Z$'s as vertices and the numbers of $e$'s on their edges with appropriate signs. (We write $W$ instead of $W_{a,b,c}$ from now on for notational simplicity.)
  More precisely, consider a $(m+n+1)$-gon $P$ with the counter-clockwise orientation on the boundary which has $m$ $X$'s and $n$ $Y$'s and $l$ $Z$'s as its vertices, which contributes to $W$. Then the pair $\{ P, P^{op} \}$ induces a term 
$$(-1)^{s(P)} \frac{ s(P) + s(P^{op}) }{ r(P) } q^{\omega(P)} x^m y^n z^l$$
where $s(P)$ is the number of $e$'s on the boundary of $P$ and $r(P)$ is the order of the symmetry of $P$. 
Here we are considering an orbifold theory (not a theory on $\Sigma$), and hence such symmetry has to be considered.
In other words, if $P$ has no such symmetry $r(P)=1$, but  if the polygon has a $\Z/k$ symmetry, this gives rise to additional equivalence relation, and hence we set $r(P)=k$.
(In the latter cases, fractional part of these polygons can be understood as orbi-discs, which will be discussed in more detail elsewhere. Unlike in Lemma \ref{signcancel}, $P$ and $P^{op}$ contribute with the same sign here, because $e$ is of even degree.
%$$m (\conste^b ) = \left( \sum a_{mnl} (q) x^m y^n z^l \right) e$$

Recall that we have computed $W(x,y,z)$ explicitly for $\mathbb{P}^1_{3,3,3}$ in Section \ref{ex333}  Theorem \ref{thm:pot333formula}.
But  determining the exact potential $W(x,y,z)$ for general $\mathbb{P}^1_{a,b,c}$ turns out to be a rather non-trivial work, 
and we provide such calculation in \cite{CHKL14} together with S.-H. Kim.

Let us explain a characteristic of $W(x,y,z)$ depending on three cases.
\begin{enumerate}
\item In a spherical case (i.e. $\frac{1}{a} + \frac{1}{b} + \frac{1}{c} >1$), the number of monomials in the potential is finite. Moreover, each coefficient  (as a power series of $q$) has finitely many terms.
\item In an Euclidean case ( $(3,3,3)$, $(2,3,6)$ or $(2,4,4)$), 
 the potential is a  weighted homogeneous polynomial  of $x, y, z$ variables which consists of finitely many monomials. 
 However, there are infinitely many polygons contributing to the same monomial so that each coefficient is an infinite series of $q$.
\item In a hyperbolic case(i.e. $\frac{1}{a} + \frac{1}{b} + \frac{1}{c} <1$), the potential has infinitely many monomials, but there are only finitely many polygons contributing to each monomial.
\end{enumerate}

The first few terms of $W(x,y,z)$ are the following.  The $XYZ$-triangle in the upper hemisphere contributes to $W(x,y,z)$, whereas
that on the lower hemisphere does not since it does not pass through $e$.
We have a holomorphic $a$-gon with corners given by $X^a$ and similarly, we have $b$-gon, $c$-gon with
corners $Y^b$ and $Z^c$.
Hence the leading terms of $W(x,y,z)$ are
%\begin{equation}\label{eq:ldabc}
$a(q)x^a + b(q)y^b + c(q)z^c + \sigma(q) xyz.$
%\end{equation}
%\begin{equation*}
%\begin{array}{l}
%W_{(3,3,3)} (x,y,z) = \phi (q) ( x^3 -  y^3 +  z^3) + \psi (q) xyz \\
%W_{(2,3,6)} (x,y,z) = q^a x^2 + b(q) y^3 + c(q) z^6 + d(q) xyz + e(q) y^2 z^2 + f(q) y z^3\\
%W_{(2,4,4)} (x,y,z) = u(q) x^2 + v(q) y^4 + v(q) z^4 + w(q) xyz
%\end{array}
%\end{equation*}
%\subsubsection{{\color{red}Isomorphism between quantum cohomology and Jacobian ring of LG mirror}}

\subsection{Geometric transform of Seidel Lagrangian in $\mathbb{P}^1_{a,b,c}$} \label{sec:transf-Seidel}
We prove the following theorem in this section.
\begin{theorem}\label{thm:MFfinal}
Consider the Seidel Lagrangian $\bar{L}$, with weak bounding cochains $b = xX + yY + zZ$.
The localized mirror functor $\mathcal{LM}^{\bar{L}, b}$ sends $\bar{L}$ to
the following matrix factorization of $W - \lambda$:
\begin{equation}\label{MFfinal}
\bordermatrix{& p^{new} & X &  Y & Z  & e & \bar{X}^{new} &\bar{Y}^{new} & \bar{Z}^{new} \cr
          p^{new} &      &    &  &  & 0  & x & y & z \cr
                    X &      &   &   &  & x & 0 & w_z & -w_y \cr
                    Y &      &   &   &  & y & -w_z & 0 & w_x \cr
                    Z &      &   &   &  & z & w_y & -w_x & 0 \cr
                    e &  0  &  w_x  & w_y & w_z  &  &   &  &  \cr
 \bar{X}^{new} &   w_x  &   0  &  -z &  y  &  &  &  &  \cr
 \bar{Y}^{new} &   w_y  &   z   & 0 & -x  &  &  &  &  \cr
 \bar{Z}^{new} &   w_z  &   -y   & x & 0  &  &  &  &  \cr                }
\end{equation}
Here $w_x$ is defined in Definition \ref{def:w}, and $\bar{X}^{new}$ is defined in Definition \ref{def:chco}.
And $\lambda$ is the constant term of $W$ and  $\lambda = 0$ in the elliptic , hyperbolic  and $(2,2,\textrm{odd})$ cases, and in the remaining
spherical cases, we have $\lambda = m_0(1) \neq 0$.
\end{theorem}

In order to find the matrix factorization $\mathcal{LM}^{\mathbb{L}} (\bar{L})$ for $\mathbb{L}=(\bar{L},b)$, we need to compute the Floer differential $m_1^{b,0}$ between $((\bar{L},b),\bar{L})$. The differential counts holomorphic strips with upper and lower boundary
mapping to  $\bar{L}$ but upper boundary is allowed to have immersed inputs. Namely, if the differential
is between immersed generators and there is one $b$ appearing in the upper boundary, the holomorphic strip is in fact a 
holomorphic triangle. If two $b$ appear in the upper boundary, the holomorphic strip is a trapezoid.
From now on, we additionally require that the critical points $e$ and $p$  of the Morse function  are symmetric with respect to the reflection about the equator of $\mathbb{P}^1_{a,b,c}$. 

We can count such (boundary deformed) holomorphic strips (or in general {\em pearly trees}) in the following way: \\

\noindent(i) {\bf Strips from $e$ to $X$, $Y$, $Z$} : \\
In this case, only constant strips contribute as we have explained in Section \ref{sec:fuksur}. One can also see this from the unital property of $e$. That is, $m_k(b,\cdots, b,e)$ vanishes except $m_2(b,e) = b$. Thus, 
$\delta(e) = m^b_1(e) = b = xX + yY + zZ.$
%The precise shape of pearly trees for this computation is given in Figure \ref{stripfi1} (a).
%\begin{figure}[htp]
%\begin{center}
%\includegraphics[scale=0.5]{stripfi1}
%\caption{(a) strips for (i) (b) strips for (ii)}\label{stripfi1}
%\end{center}
%\end{figure}

\noindent(ii) {\bf Strips from $\bar{X}, \bar{Y}, \bar{Z}$ to $p$} :\\
This case is similar to (i) in that only constant strips contribute. Namely, the $\AI$ structure with output on $p$ was
defined in such a way that only constant strips contribute so that the degree $3$ parts of $\delta(\bar{X})$, $\delta(\bar{Y})$, $\delta(\bar{Z})$ are $x p$, $y p$, $z p$ respectively.
%The shape of pearly trees for this computation is given in Figure \ref{stripfi1} (b)

\begin{figure}[htp]
\begin{center}
\includegraphics[scale=0.5]{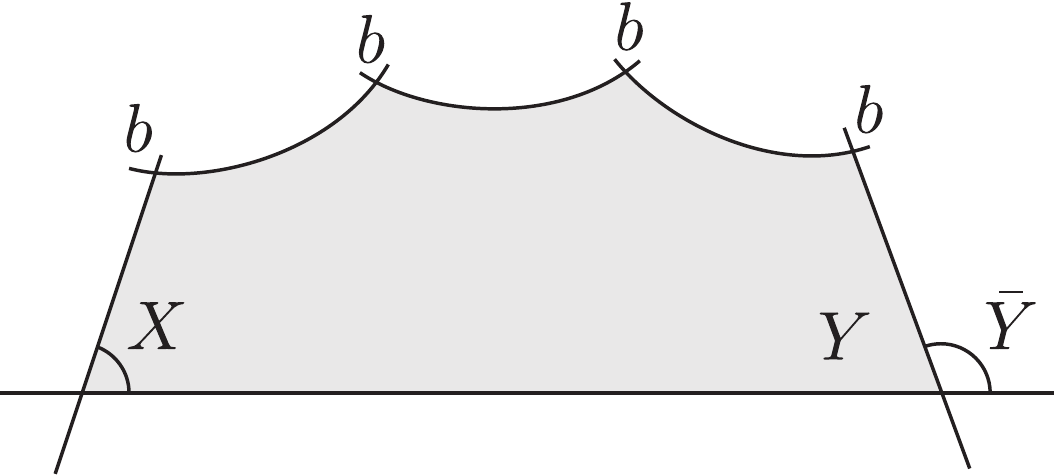}
\caption{strips from $X$ to $\bar{Y}$}\label{fig:XbbbbY}
\end{center}
\end{figure}

\noindent(iii) {\bf Strips from $X$, $Y$, $Z$ to $\bar{X}, \bar{Y}, \bar{Z}$} :\\
As illustrated in Figure \ref{fig:XbbbbY}, holomorphic strips from one of $X,Y,Z$ to one of $\bar{X}, \bar{Y}, \bar{Z}$ come from the same 
holomorphic polygons which were used to compute the potential $W$, but their contributions are different. For example, the number
of times $e$ appears is irrelevant here. We fix a polygon $P$ which were used in computing $W$,
and determine its contribution to this counting of holomorphic strips. 
%Contribution comes from the polygons counted for the superpotential. 
Let P be such a polygon with the counter-clockwise orientation on the boundary Lagrangian. We can write it as $P =  x_1 x_2 \cdots x_{k-1}$ where $x_i$'s are vertices of $P$, arranged in counterclockwise direction. ($x_i =$ $x$ or $y$ or $z$). Here $P$ should be considered as a cyclic word, or one may regard index $i$ of $x_i$
as a number modulo $k$ (for example  $x_{k+1} = x_1$).

 Let $P^{op}$ be the reflection image of $P$ with respect to the equator of $\mathbb{P}^1_{a,b,c}$ which can be written as the word $x_k x_{k-1} \cdots x_1$ also in the counter-clockwise arrangement of vertices. We can use $P$ to consider $m_1^b(X_i)$ contribution to $\overline{X_{i+1}}$,
 or use $P^{op}$ to consider its contribution to $\overline{X_{i-1}}$.
 Then the pair $\{ P,P^{op} \}$ contributes to $\delta$ as
\begin{equation}\label{dPPop}
\delta_{ \{ P,P^{op} \} } : X_i {\mapsto} (-1)^{s(P)} q^{\omega(P)} \left(  x_1 \cdots \check{x_i} \check{x_{i+1}} \cdots x_k \overline{X_{i+1}} - x_k \cdots \check{x_i} \check{x_{i-1}} \cdots x_1 \overline{X_{i-1}}
  \right).
\end{equation}
(the minus sign appearing in the expression \eqref{dPPop} is due to Lemma \ref{signcancel}.)  Combining contribution for all $i$, we get
$$\delta_{ \{ P,P^{op} \}} (X) =  (-1)^{s(P)} q^{\omega (P)} \frac{1}{r(P)} \left( \alpha_{yx} (P)  \frac{M_P}{xy} \bar{Y} +  \alpha_{xz} (M_P)  \frac{P}{xy} \bar{Z} \right),$$
where we regard $M_P$ as the commutative monomial $x_1 \cdots x_k$ and $\alpha_{xy}$ is the number of ``$xy$" appearing in $P$ minus that of ``$yx$" (and similar for $\alpha_{xz}$).

\begin{lemma}
$\delta_{ \{ P,P^{op} \}} (X)$ does not contain nontrivial $\OL{X}$. The same goes for $Y, Z$.
\end{lemma}
\begin{proof}
Suppose the polytope has a boundary edge $XX$. Then contribution from $P$ and from $P^{op}$ cancel out.
\end{proof}
Also, it is easy to check by simple combinatorics that if $P$ consists of less than three variables, say $x_i$ and $x_j$, then $\alpha_{x_i x_j} (P) =0$. Therefore, if we define $P_{xyz}^{+}$ to be the set of polygons contributing to $W$ which contain all three variables and have counter-clockwise orientation on the boundary, then degree 2 part of $\delta(X)$ is given by
\begin{equation}\label{eq:deltaXP}
\sum_{P \in P_{xyz}^+ } (-1)^{s(P)} q^{\omega (P)} \frac{\alpha_{xy} (P)}{r(P)} \frac{M_P}{xy} \bar{Y} + \sum_{P \in P_{xyz}^+ } (-1)^{s(P)}  q^{\omega (P)} \frac{\alpha_{xz} (P)}{r(P)} \frac{M_P}{xz} \bar{Z}.
\end{equation}

Let $\mathcal{P}_{a,m}$ be the set of positively oriented (w.r.t. the orientation on $\bar{L}$) polygons which have the symplectic area $a$, and induces the monomial $m$ in the potential. Define $\alpha_{xy} (\mathcal{P}_{a,m})$ by
$$\alpha_{xy} (\mathcal{P}_{a,m}) := \sum_{P \in \mathcal{P}_{a,m}} \frac{(-1)^{s(P)}}{r(P)} \alpha_{xy} (P).$$
$\alpha_{x_i x_j} (\mathcal{P}_{a,m})$ is similarly defined for $x_i, x_j \in \{x,y,z\}$.
Then Equation \eqref{eq:deltaXP} can be rewritten as
\begin{equation}\label{eq:deltaXPre}
\textnormal{degree 2 part of} \,\,\, \delta(X) = \sum_{a,m} q^{a} \left( \alpha_{xy} (\mathcal{P}_{a,m}) \frac{m}{xy} \bar{Y} + \alpha_{xz} (\mathcal{P}_{a,m}) \frac{m}{xz} \bar{Z} \right).
\end{equation}
Note that negative powers of $x,y,z$ do not appear in the above sum since $\alpha_{x_i x_j} (\mathcal{P}_{a,m}) =0$ for $m$ only consisting of less than two variables among $\{x,y,z\}$.
Let us compute $\left< \delta^2 (X), p \right>$ from the expression \eqref{eq:deltaXPre} which should be zero.
$$\left< \delta^2 (X), p \right>= \sum_{a,m} q^{a}  \left( \alpha_{xy} (\mathcal{P}_{a,m}) + \alpha_{xz} (\mathcal{P}_{a,m}) \right) \frac{m}{x}.$$
Therefore 
we have $\alpha_{xy} (\mathcal{P}_{a,m})  + \alpha_{xz} (\mathcal{P}_{a,m}) = 0$ for all $(a,m)$ or, equivalently $\alpha_{xy} (\mathcal{P}_{a,m}) = \alpha_{zx} (\mathcal{P}_{a,m})$.
So, it makes sense to define a $\Z$-valued function $\theta$ by
$$\theta(a,m) := \alpha_{xy} (\mathcal{P}_{a,m}) = \alpha_{yz} (\mathcal{P}_{a,m}) = \alpha_{zx} (\mathcal{P}_{a,m}).$$
We remark that this identity can be also proved by an easy combinatorial argument.
%\begin{remark}
%We sort polygons by their areas and associated monomials ($\mathcal{P}_{a,m}$) since there might be more than one polygons with the same associated monomial and the same area which are not reflection pairs.
%\end{remark}
%}
%
%Since the coefficient of $p$ in $\delta^2 (X)$ should vanish, we have $\alpha_{xy} (P)  + \alpha_{xz} (P) = 0$ for all $P \in W_{xyz}^+$ or, equivalently $\alpha_{xy} (P) = \alpha_{zx} (P)$ (by arranging the coefficient of $p$ in $\delta^2 (X)$ with respect to the order of $q$'s power). 
%
%So, it makes sense to define $\theta : P_{xyz}^{+} \to \Z$ by
%$$\theta(P) := \alpha_{xy} (P) = \alpha_{yz} (P) = \alpha_{zx} (P).$$
Now we introduce a change of coordinates which will help us to get a nice form of the matrix factorization later.
\begin{defn}\label{def:chco}
We define the following change of  coordinates.
\begin{equation}\label{newcoord}
 \bar{X}^{new} :=\gamma \bar{X}, \quad \bar{Y}^{new} := \gamma\bar{Y}, \quad  \bar{Z}^{new} := \gamma \bar{Z}
\end{equation}
where $\gamma =  \sum_{a,m} q^{a} \theta (a,m) \frac{m}{xyz}$. 
\end{defn}
The above coordinate change is essential to make  the mirror matrix factorization to
have ``contraction and wedge product'' form.

The term $``q^{\epsilon}" xyz$ in $W$ from the minimal triangle with area $\epsilon$ in $\mathbb{P}^1_{a,b,c} \setminus \{p_a, p_b, p_c\}$ guarantees that this coordinate change is, indeed, invertible. With this new basis, degree $2$ parts of $\delta(X)$, $\delta(Y)$, $\delta(Z)$ are respectively,
$z \bar{Y}^{new} - y \bar{Z}^{new}, \quad x \bar{Z}^{new} - z \bar{X}^{new}, \quad y \bar{X}^{new} - x \bar{Y}^{new}.$
Accordingly, letting $p^{new} = \gamma p$,  (ii) implies that the degree 3 parts of $\delta (\bar{X}^{new})$, $\delta (\bar{Y}^{new})$, $\delta (\bar{Z}^{new})$ are $x p^{new}$, $y p^{new}$, $z p^{new}$ respectively. \\

\noindent (iv) {\bf Strips from $X$, $Y$, $Z$ to $e$ and strips from $p^{new}$ to $\bar{X}^{new}$, $\bar{Y}^{new}$, $\bar{Z}^{new}$} :\\
Counting strips from $p^{new}$ to $\bar{X}^{new}$, $\bar{Y}^{new}$, $\bar{Z}^{new}$ is equivalent to counting the strips from $p$ to $\bar{X}$, $\bar{Y}$, $\bar{Z}$ since 
we have multiplied $\gamma$ commonly to all of them to obtain the new coordinates.

 Again consider the polygon $P = x_1 x_2 \cdots x_k$ which contributes to the potential. Roughly speaking, $P$ contributes to $\delta (p^{new})$ by a linear combination of $\frac{M_P}{x_i} \bar{X_i}$ and coefficients depend on the number of $p$'s the edge of $P$ on which the input marked point lies. The similar happens for the $e$-coefficients of $\delta(X)$, $\delta(Y)$, $\delta(Z)$. 

However, it is a little more complicated to find the coefficients precisely since each edge of $P$ can have the different number of $p$ on it. Instead, we will show that the counting strips from $X$, $Y$, $Z$ to $e$ are equivalent to counting $p$ to $\bar{X}$, $\bar{Y}$, $\bar{Z}$.

\begin{defn}\label{def:w}
We define $w_x$, $w_y$, $w_z$ such that $\delta(p) = w_x \bar{X} + w_y \bar{Y} + w_z \bar{Z}.$
($\delta(p)$ does not have an $e$-term since only Morse flows contribute to $\left< \delta(p),e\right>$.)
\end{defn}
Now we claim that 
$\left< \delta{X}, e \right> =w_x, \,\, \left< \delta{Y}, e \right> =w_y, \,\, \left< \delta{Z}, e \right> =w_z.$
This is basically because $e$ is the reflection image of $p$ with respect to the equator. We have a symmetry with positive signs here unlike in Lemma \ref{signcancel}, because the degree difference between $p$ and $e$ is odd. Suppose $P$ has $k$ $e$'s on the edge $\overline{X_i X_{i+1}}$.  By symmetry, its reflection image $P^{op}$ has $k$ $p$'s on the edge $\tau( \overline{ X_i X_{i+1} } )$.
\begin{figure}[htp]
\begin{center}
\includegraphics[scale=0.4]{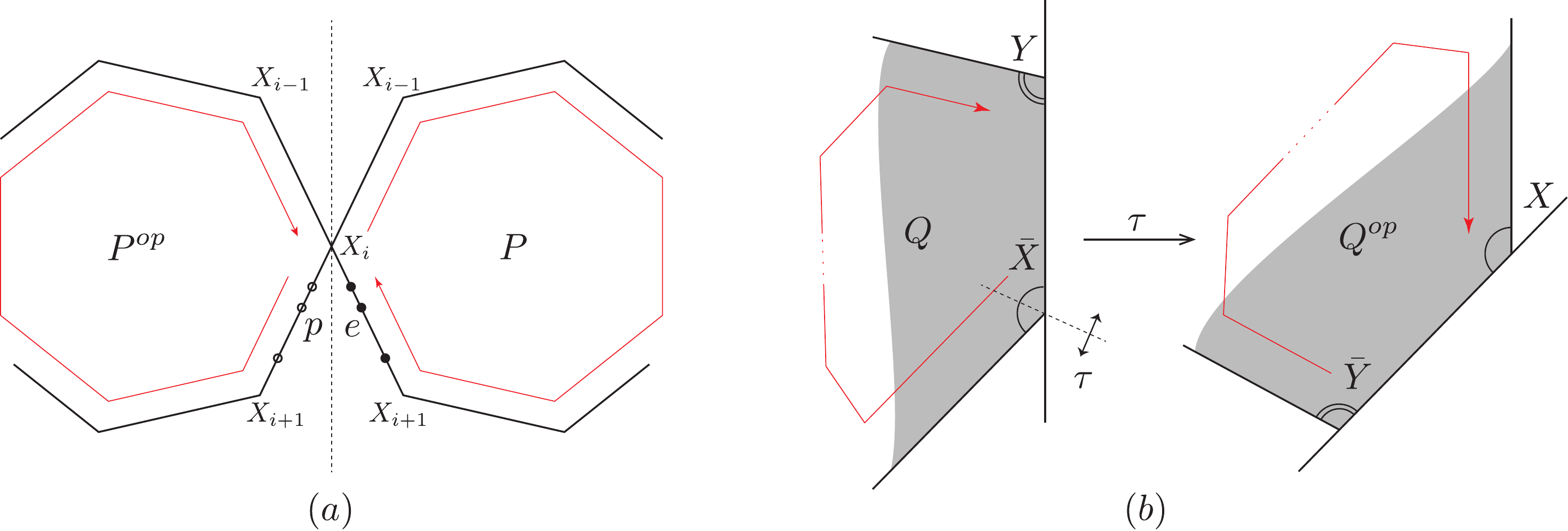}
\caption{Symmetry (a) between  strips from $p$ to $\bar{X_i}$ and from $X_i$ to $e$ (b) between  strips from $\bar{X}$ to $Y$ and from $\bar{Y}$ to $X$}\label{symmcombined} 
\end{center}
\end{figure}
So, $P$ contributes to $\left< \delta (X_i), e \right>$ $k$-times, and $P^{op}$ contributes to $\left< \delta(p), \bar{X_i} \right>$ $k$-times. (It is clear from (a) of Figure \ref{symmcombined} that the strip from $p$ to $X_i$ corresponds to the strip from $X_i$ to $e$ by the reflection and vice versa.) This implies $\left<\delta (p) , \bar{X_i} \right> = \left<\delta (X_i), e \right>$ for all $X_i = X, Y, Z$.\\

\noindent(v) {\bf Further properties from the reflection symmetry} :
\begin{lemma}
We have $\left< \delta(\bar{X}^{new}), Y\right>=-\left< \delta(\bar{Y}^{new}), X\right>$, $\left< \delta(\bar{Y}^{new}), Z\right>=-\left< \delta(\bar{Z}^{new}), Y\right>$, $\left< \delta(\bar{Z}^{new}), X\right>=-\left< \delta(\bar{X}^{new}), Z\right>$.
\end{lemma}
\begin{proof}
Without loss of generality, it suffices to prove the first identity. Consider a strip $Q$ from $\bar{X}$ to $Y$.  Its reflection image $Q^{op}$ is a strip from $\bar{Y}$ to $X$ as shown in (b) of Figure \ref{symmcombined}.
%\begin{figure}[htp]
%\begin{center}
%\includegraphics[scale=0.5]{barXtoY}
%\caption{A strip from $\bar{X}$ to $Y$ and its reflection image}\label{barXtoY} 
%\end{center}
%\end{figure}
Considering the orientation of the Lagrangian along the boundaries of $Q$ and $Q^{op}$, we get $\left<\delta(\bar{X}),Y \right> = -\left<\delta(\bar{Y}),Z \right>$ by Lemma \ref{signcancel} and hence the first identity follows.
\end{proof}

\begin{corollary}
We have
$\left< \delta(\bar{X}^{new}), X\right> = \left< \delta(\bar{Y}^{new}), Y\right> = \left< \delta(\bar{Z}^{new}), Z\right> =0.$
\end{corollary}
\begin{proof}
We can apply the same argument as in the proof of the previous lemma for the map from $\bar{X}$ to $X$. If a strip from $\bar{X}$ to $X$ has other immered input $b$'s, then
the reflection image of such a strip is different from the strip itself, and hence the same cancellation argument works. 
If a strip from $\bar{X}$ to $X$ is a bigon (without other immersed inputs), then its two corners are both $\bar{X}$, but this bigon is counted twice as a strip from $\bar{X}$ to $X$
(in opposite directions), and they cancel out each other. 
\end{proof}

%\begin{remark}
%It is not easy to find all the strips from $\bar{X}$ to $Y$ for arbitrary weights $(a,b,c)$, but we can always find the very first one which gives the zeroth order term in $\delta (\bar{X})$ with respect to $q$. It is the pearly tree of the shape shown in Figure \ref{cluster}.
%\begin{figure}[htp]
%\begin{center}
%\includegraphics[scale=0.5]{cluster}
%\caption{A strip from $\bar{X}$ to $Y$}\label{cluster} 
%\end{center}
%\end{figure}
%
%This strip corresponds to the term $1 \cdot xy Y$ in $\delta (\bar{X})$ and similar pearly trees exist for other pairs of variables.
%%(or $1 \cdot zx X$ in $\delta (\bar{Y})$ or $1 \cdot yz Z$ in $\delta (\bar{Z})$ depending to the position of $e$ and $p$.)
%\end{remark}

In summary, we have the matrix factorization of $W - \lambda$ shown below:
\begin{equation}
\bordermatrix{& p^{new} & X &  Y & Z  & e & \bar{X}^{new} &\bar{Y}^{new} & \bar{Z}^{new} \cr
          p^{new} &      &    &  &  & 0  & x & y & z \cr
                    X &      &   &   &  & x & 0 & -f & h \cr
                    Y &      &   &   &  & y & f & 0 & -g \cr
                    Z &      &   &   &  & z & -h & g & 0 \cr
                    e &  0  &  w_x  & w_y & w_z  &  &   &  &  \cr
 \bar{X}^{new} &   w_x  &   0  &  -z &  y  &  &  &  &  \cr
 \bar{Y}^{new} &   w_y  &   z   & 0 & -x  &  &  &  &  \cr
 \bar{Z}^{new} &   w_z  &   -y   & x & 0  &  &  &  &  \cr                }
 \end{equation}
where we write holomorphic strips from $j$ to $i$ in the $(i,j)$-th entry and we define 
$$f := \left< \delta(\bar{X}^{new}), Y\right>, g:=\left< \delta(\bar{Y}^{new}), Z\right>, h:=\left< \delta(\bar{Z}^{new}), X\right>.$$
In fact, we can identify $f$, $g$, $h$ as follows. We look at the product of the second row of the bottom left 4 by 4 matrix and the third column of the upper right 4 by 4 matrix, which should be zero. Therefore
%\begin{equation}\label{funceq}
$y w_x + y g = 0$ and $g =- w_x.$
%\end{equation}
By the same manner, we get $h =- w_y$ and $f =- w_z$.  Since this is a matrix factorization of $W- \lambda$, it satisfies $ W-\lambda = x w_x + y w_y + z w_z.$
%Recall that $\lambda$ is a constant term of $W$. This proves the theorem.

\subsubsection{Wedge-contraction formulation}\label{sec:cw}
Let us express the matrix factorization \eqref{MFfinal} in the ``contraction-wedge'' form.
The Floer homology  $CF^*(\bar{L},\bar{L})$ has 8 generators (as a module),
which are given by $e, X, Y, Z, \bar{X},\bar{Y},\bar{Z},p$. This can be identified
with the exterior algebra if we define 
$ Y \wedge Z = \bar{X}, Z \wedge X = \bar{Y}, X \wedge Y = \bar{Z}, 1=e, X\wedge Y \wedge Z = p.$
But this exterior algebra structure \emph{require quantum corrections} in order to obtain the matrix
factorization \eqref{MFfinal} (see also Definition \ref{def:chco}).

%For the change of variables in Definition \ref{def:chco}, we need to modify the exterior algebra structure as follows.
\begin{defn}\label{def:newext}
We define the new exterior algebra $\bigwedge\nolimits_{new}^*\langle X,Y,Z \rangle$  over Novikov field as follows.
$$X \wedge^{new} Y := \bar{Z}^{new}, \quad Y \wedge^{new} Z := \bar{X}^{new}, \quad Z \wedge^{new} X := \bar{Y}^{new},$$
where $\bar{X}^{new}, \bar{Y}^{new},\bar{Z}^{new}$ are defined in Definition \ref{def:chco}.
We also define
$$X \wedge^{new} \bar{X}^{new} = Y \wedge^{new} \bar{Y}^{new} = Z \wedge^{new} \bar{Z}^{new} =p^{new}.$$
Here we define $1 \wedge^{new} (\cdot)$ to be an identify operation.
Remaining $\wedge$-operations are defined to be trivial.

Contraction $\iota^{new}$ is defined in a similar way.
Namely, 
$\iota^{new}_X  \bar{Z}^{new} = Y, \quad \iota^{new}_Y  \bar{X}^{new} = Z, \quad \iota^{new}_Z  \bar{Y}^{new} = X. $
\end{defn}
Recall from  Definition \ref{def:chco} that leading term of $\gamma$ comes from
the minimal triangle $\Delta XYZ$ (the shaded triangle in Figure \ref{Seidel_Lag}).
Hence the leading term of $\wedge^{new}$ is the usual exterior algebra structure
(multiplied by the exponentiated area of $\Delta XYZ$). This may be understood as a quantum corrected exterior algebra structure.  Then we have

\begin{corollary}\label{cor:seiext}
The matrix factorization in Theorem \ref{thm:MFfinal} is isomorphic to
\begin{equation}\label{eq:cw}
\left(\bigwedge\nolimits_{new}^*\langle X,Y,Z \rangle, x X \wedge^{new} (\cdot) + y Y \wedge^{new} (\cdot) + z Z \wedge^{new} (\cdot) + w_x \, \iota_X^{new} + w_y \, \iota_Y^{new} + w_z \, \iota_Z^{new} \right).
\end{equation}
\end{corollary}

We directly compute the mirror matrix factorization for $\bar{L}$ explicitly for $(a,b,c) = (3,3,3)$.

\begin{prop} \label{prop:MF333}
For $(a,b,c) = (3,3,3)$, the entries $w_x$, $w_y$ and $w_z$ in the matrix factorization in Theorem \ref{thm:MFfinal} are given as follows:
%\begin{eqnarray*}
%\gamma &=& q+ \sum_{k=1}^{\infty} \left(  (-1)^{3k+1} q^{(6k+1)^2} -  (-1)^{3k} q^{(6k-1)^2} \right) \\
%&=& -q + \sum_{k=1}^{\infty} (-1)^{k+1} \left( \psi_k^+ (q_\alpha) +  \psi_k^- (q_\alpha) \right)
%\end{eqnarray*}
%\begin{equation}\label{wx}
%\begin{array}{rl}
%w_x \!\!\!\! &=   x^2 \displaystyle\sum_{k=0}^{\infty}  (-1)^{k+1}(2k+1) \phi_k (q_\alpha^3) -yz q_\alpha   \\
% &\quad +yz \displaystyle\sum_{k=1}^{\infty}  (-1)^{k+1} \left(  (2k+1) \psi_k^{+} (q_\alpha) - (2k -1) \psi_k^{-} (q_\alpha) \right) 
%\end{array}
%\end{equation}
$$w_x =  x^2 \displaystyle\sum_{k=0}^{\infty}  (-1)^{k+1}(2k+1) q_\alpha^{ (6k+3)^2} +yz \displaystyle\sum_{k=1}^{\infty}  (-1)^{k+1} \left(  (2k+1) 
q_\alpha^{(6k+1)^2} - (2k -1) q_\alpha^{(6k-1)^2} \right) -yz q_\alpha,$$
$$w_y = y^2 \sum_{k=0}^{\infty}  (-1)^{k} (2k+1) q_\alpha^{ (6k+3)^2} +  xz \sum_{k=1}^{\infty}  (-1)^{k+1} \left(  2k q_\alpha^{(6k+1)^2} -  2k 
q_\alpha^{(6k-1)^2} \right), $$
$$w_z = z^2 \sum_{k=0}^{\infty}  (-1)^{k+1} (2k+1) q_\alpha^{ (6k+3)^2}+  xy \sum_{k=1}^{\infty}  (-1)^{k+1} \left(  2k q_\alpha^{ (6k+1)^2} -  2k q_\alpha^{ (6k-1)^2} \right). $$
where $q_\alpha$ is the area of a minimal $xyz$-triangle.
%$$w_x =  x^2 \displaystyle\sum_{k=0}^{\infty}  (-1)^{k+1}(2k+1) \phi_k (q_\alpha^3) +yz \displaystyle\sum_{k=1}^{\infty}  (-1)^{k+1} \left(  (2k+1) \psi_k^{+} (q_\alpha) - (2k -1) \psi_k^{-} (q_\alpha) \right) -yz q_\alpha,$$
%$$w_y = y^2 \sum_{k=0}^{\infty}  (-1)^{k} (2k+1) \phi_k (q_\alpha^3) +  xz \sum_{k=1}^{\infty}  (-1)^{k+1} \left(  2k \psi_k^{+} (q_\alpha) -  2k \psi_k^{-} (q_\alpha) \right), $$
%$$w_z = z^2 \sum_{k=0}^{\infty}  (-1)^{k+1} (2k+1) \phi_k (q_\alpha^3) +  xy \sum_{k=1}^{\infty}  (-1)^{k+1} \left(  2k \psi_k^{+} (q_\alpha) -  2k \psi_k^{-} (q_\alpha) \right). $$
The coordinate change $\gamma$ in Equation \eqref{newcoord} is
$\gamma = -q + \sum_{k=1}^{\infty} (-1)^{k+1} \left( \psi_k^+ (q_\alpha) +  \psi_k^- (q_\alpha) \right)$.
\end{prop}

We will consider $W$ as an element of $\Lambda [[x,y,z]]$.  In a joint work \cite{CHKL14} with Kim, we prove that $W$ is convergent over $\C$ for $\bP^1_{a,b,c}$.

\begin{corollary}\label{cor:cptgen1}
The matrix factorization \eqref{eq:cw} is a compact generator of the matrix factorization category $\mathcal{MF}(W)$ if $\frac{1}{a} + \frac{1}{b} + \frac{1}{c} \leq 1$.
\end{corollary}
\begin{proof}
This follows from the result of Dyckerhoff \cite{Dy}. We check that our setting fulfills the condition of \cite[Theorem 4.1]{Dy}. First of all, $\Lambda [[x,y,z]]$ is a regular local ring (of Krull dimension $3$ over $\Lambda$) since $\Lambda$ is a field. In Theorem 13.1 of \cite{CHKL14}, we  show that the hypersurface $W^{-1} (0)$ has an isolated singularity only at the origin under the assumption 
 $\frac{1}{a} + \frac{1}{b} + \frac{1}{c} < 1$.
 
 If $\frac{1}{a} + \frac{1}{b} + \frac{1}{c} =1$, we have explicit expressions of the potentials from \cite[Section 9,10]{CHKL14}, and since they are polynomials one can check that they have the desired property by hands. Indeed, there are coordinate changes which transform these polynomials into much simpler forms which makes computation easier. See \cite[Proposition 13.2]{CHKL14}. 
 
\end{proof}
In the spherical cases ($\frac{1}{a} + \frac{1}{b} + \frac{1}{c} >1$), the computation of $W$ shows that $W$ appear to be  a Morsification of the leading term $x^a+ y^b+z^c - \sigma xyz$ and
hence  has several critical points other than the origin. See \cite{CHKL14} for the complete computation of $W$.

In the paper of Dycherhoff \cite{Dy}, the matrix factorization of wedge and contraction type (like Equation \eqref{eq:cw}) corresponds to the skyscraper sheaf at the critical point, and denoted as  $\Lambda^{{\rm stab}} [1]$.

\subsection{Mirror functor on the morphism level}\label{subsec:mirrormorlev}
Now we prove that our functor induces an isomorphism between $\Hom_{{\Fuk} (\mathbb{P}^1_{a,b,c})} (\bar{L}, \bar{L}) = CF (\bar{L}, \bar{L})$ and $\Hom_{\mathcal{MF} (W)} (P_{\bar{L}}, P_{\bar{L}})$. From the construction in Section \ref{sect:geom_fctor}, we already have an $A_\infty$-morphism 
$\{ \Phi_k \}_k : CF (\bar{L}, \bar{L}) \to \Hom (P_{\bar{L}}, P_{\bar{L}})$
which is the $A_\infty$-functor $\mathcal{LM}^{\mathbb{L}}$ from ${\Fuk} ( \mathbb{P}^1_{a,b,c} )$ to ${\mathcal{MF}} (W)$ restricted to a single object $\bar{L}$.  

We first show that $m_1$ is identically zero on $CF(\bar{L}, \bar{L})$ for $\frac{1}{a} + \frac{1}{b} + \frac{1}{c} \leq 1$, hence implies that Floer homology has rank 8 in these cases.
Here we consider $m_1$, which is not the boundary deformed one $m_1^{b,0}$. Let us first describe the elliptic and hyperbolic cases.
\begin{lemma}\label{lem:xleq1}
 Suppose $\frac{1}{a} + \frac{1}{b} + \frac{1}{c} \leq 1$. Then $m_1$ on $CF(\bar{L},\bar{L})$ vanishes.
\end{lemma}
\begin{proof}
If $a,b,c \geq 3$, then ${\rm deg \, } m_k^{higher} \geq 6 - 3k + 2 \min\{ a-3 ,b-3, c-3\} \geq 6-3k$  and hence
$m_1^{higher}$ (and hence $m_1$) has degree bigger than or equal to $3$, and hence is zero.

Hence it is enough to consider the case when one of $a,b,c$, say $c$, is $2$.
Note that $m_1 (e)$ as well as the degree-$3$ part of $m_1(\bar{B})$ for $B \in \{X,Y,Z\}$ vanishes, since the operation with an input $e$ or output $p$ vanishes unless it is related to $m_2$-product with a unit. Also  $m_1$ between $p$ and $e$ vanishes which are Morse differentials.

For $A,B \in \{ X,Y,Z\}$, the differential $m_1$ from $A$ to $e$ or from $p$ to $\bar{B}$ could be non-zero if there is a polygon with only one corner (1-gon).
But in Lemma 4.5 (1) of \cite{CHKL14}, we prove that Seidel Lagrangian in the universal cover is topologically a line, and there cannot by such a 1-gon
in the elliptic and hyperbolic cases. Hence such differential vanishes. 

For $A,B \in \{ X,Y,Z\}$, the differential $m_1$ from $A$ to $\bar{B}$ or  from $\bar{B}$ to $A$ could be non-zero if there is a bigon.
In Corollary 4.6 of \cite{CHKL14}, we prove that there is no bigon in the elliptic and hyperbolic cases, except the minimal bigon with two corners given by $C$.
Hence  $m_1$ from $Z$ to $\bar{Z}$ could be non-trivial, and the other differentials vanish. But such a bigon is counted twice as a strip from $Z$ to $\bar{Z}$
travelling in opposite directions. And one can check that these two contributions always cancel out.
This proves the claim.
\end{proof}
Now let us consider the spherical cases. It turns out that in the cases of $(a,b,c) = (2,2,\textrm{odd})$, $m_1$ does not vanish, and in  all other cases, $m_1$ vanishes.
\begin{lemma}
For $(a,b,c) \neq (2,2,\textrm{odd})$ with $\frac{1}{a} + \frac{1}{b} + \frac{1}{c} >1$, 
$m_1$ on $CF(\bar{L},\bar{L})$ vanishes.
\end{lemma}
\begin{proof}
The lifts of Seidel Lagrangians in the universal cover $S^2$ are (topological) circles, which bisect the area of $S^2$ (Lemma 12.1 \cite{CHKL14}). If a bi-gon connecting $p$ and $q$ appears, by extending each edge, such a bigon can be realized as an intersection of two bisecting circles, which are lifts of Seidel Lagrangians in $S^2$.  There are two possibilities from the deck
transformation group $G(\cong \pi_1^{orb} (\mathbb{P}^1_{a,b,c}))$ action on $S^2$ whose quotient is the
orbi-sphere.

If $p$ and $q$ are in the same orbit of the $G$-action, then
it defines $m_1$ as a map  $X_i \to \bar{X}_i$ or $\bar{X}_i \to X_i$ for
some $i \in \{ 1,2,3\}$.
However, such a bigon contributes twice as a map from $p$ to $q$ and as a map
from $q$ to $p$, with opposite signs. Since $p$ and $q$ are in the same orbit, they contribute to $m_1$ to the same map with opposite sign. Thus, their contributions to $m_1$ is zero after cancellation.

If $p$ and $q$ are not in the same orbit of $G$-action, then
since the bigon is on the sphere, there is another bigon (on the opposite hemisphere) connecting $p$ to $q$ with the same area. There cannot
be an action of the group $G$ which maps one bigon to the other one,
since such a group action has to be a rotation by $\pi$ at the intersection point $p$ or $q$, but the $G$-action does not have $p$ or $q$ as a fixed point.
Therefore, one can show that the total contribution to $m_1$ vanishes.
\end{proof}
%$m_2^{higher}$ has degree $\geq 0$. {\color{red} Assume $a,b,c \geq 4$. Then, ${\rm deg\,} m_2^{higher} \geq 1$. So only $m_2^{higher} (A,B)$-type term can survive where $A,B \in \{X,Y,Z\}$. Moreover, such a term should be a multiple of $q$ by the degree reason. But, all contributions to the coefficient of $q$ come from constant discs. So, $m_2^{higher} \equiv 0$. } 
In the $(2,2,\textrm{odd})$ cases, a holomorphic 1-gon gives a non-trivial $m_1$ from  one of immersed generators to $e$. Since the unit $e$ is a coboundary, $m_1$ homology vanishes in these cases.
The potential $W$ in these cases have a linear term (see \cite{CHKL14}), so that $x=y=z=0$ is not a critical point of $W$.

\begin{remark}
We remark that if $a,b,c \geq 4$, then the degree argument for $m_2$ shows that $m_2^{higher} \equiv 0$. Hence the contribution to $m_2$ comes only from the $XYZ$ triangles in the upper and lower hemisphere.
\end{remark}

For the rest of the paper, we assume that $\frac{1}{a} + \frac{1}{b} +\frac{1}{c} \leq 1$.
We next show that our $A_\infty$-functor on the morphism level induces an isomorphism between $CF (\bar{L}, \bar{L})$ and $\Hom (P_{\bar{L}}, P_{\bar{L}})$ as $A_\infty$-algebras.
Recall that $\Phi_1$ is given by
$\Phi_1 :  CF (\bar{L}, \bar{L}) \to \Hom (P_{\bar{L}}, P_{\bar{L}}) \quad p \mapsto m_2^{b,0} ( \cdot, p).$
We claim that $\Phi_1$ induces an isomorphism on the cohomology level. 
In order to see this, we find $\Psi : \Hom (P_{\bar{L}}, P_{\bar{L}}) \to CF (\bar{L}, \bar{L})$ which is a right inverse of $\Phi_1$. Then it induces a right inverse of the cohomology level map 
\begin{equation}\label{cohomPhi}
[\Phi_1] :  HF (\bar{L}, \bar{L}) \to H \left( \Hom (P_{\bar{L}}, P_{\bar{L}})  \right).
\end{equation}

Let ${\bf 1}_{\bar{L}}$ be the identity morphism (unit) in $CF (\bar{L}, \bar{L})$, which is nothing but the minimum $e$ of the chosen Morse function on $\bar{L}$.  Then we define $\Psi$ by $\Psi (\phi) = \phi ({\bf 1}_{\bar{L}} )|_{b=0}$
for $\phi \in \Hom (P_{\bar{L}}, P_{\bar{L}})$. (Recall that as a module $P_{\bar{L}}$ is isomorphic to $\Lambda [[x,y,z]] \otimes_{\Lambda_0} CF(\bar{L}, \bar{L})$, which corresponds to a trivial vector bundle over $Spec \,\Lambda [[x,y,z]]$.)

\begin{lemma}\label{rightinv}
$\Psi$ is a chain map and $\Psi \circ \Phi_1 (p) = p$ for all $p \in CF(\bar{L}, \bar{L})$.
\end{lemma}

\begin{proof}
To see $\Psi$ is a chain map, we have to show $\Psi (d \phi) =0$ for $\phi \in \Hom (P_{\bar{L}}, P_{\bar{L}})$ since the differential on $CF(\bar{L}, \bar{L})$ is identically zero. (See (iii) of Subsection \ref{FloerSeidel}.) But,
\begin{eqnarray*}
 \Psi (d \phi) &=& \left(\delta \circ \phi ({\bf 1}_{\bar{L}} ) + (-1)^{{\rm deg}\, (\phi)} \phi \circ \delta ({\bf 1}_{\bar{L}} ) \right)|_{b=0} \\
 &=& m_1^{b=0,0} (\phi ({\bf 1}_{\bar{L}} )) + (-1)^{{\rm deg}\, (\phi)} \phi \left(m_1^{b=0,0} ({\bf 1}_{\bar{L}}) \right) =0.
 \end{eqnarray*}
Here $m_1^{b=0,0}=0$ since we turn off the boundary deformation $b$.

The second statement is almost direct from the definition of the unit of $A_\infty$-algebra.
\begin{eqnarray*}
\Psi \circ \Phi_1 (p) &=& m_2^{b,0} ( {\bf 1}_{\bar{L}}, p)|_{b = 0} \\
&=& m_2 ({\bf 1}_{\bar{L}}, p) = p
\end{eqnarray*}
by the property of the unit.
\end{proof}

From the lemma, $\Phi_1$ is injective, and so is $[\Phi_1]$. Note that both sides of Equation \eqref{cohomPhi} has rank $2^3 =8$ over $\Lambda$. (See \cite[Chapter 4]{Dy} for details about $H (\Hom (P_{\bar{L}},P_{\bar{L}})$). This proves that $[\Phi_1]$ is an isomorphism and hence $\{ \Phi_k \}_k$ is an $A_\infty$-isomorphism.

\begin{remark}
This argument actually proves that the $A_\infty$-functor obtained from the method developed in Section \ref{GMir} is faithful whenever the Floer differential of a Lagrangian vanishes.
\end{remark}

\subsection{Homological mirror symmetry of $\mathbb{P}^1_{a,b,c}$}
We now prove the homological mirror symmetry for $\mathbb{P}^1_{a,b,c}$ under the assumption that $\frac{1}{a} + \frac{1}{b} + \frac{1}{c} \leq 1.$
%We have constructed an $A_\infty$-functor $\mathcal{LM}^\BL$ from ${\Fuk} (\mathbb{P}^1_{a,b,c})$ to ${\mathcal{MF}} (W)$ which becomes an isomorphism when restricted to $\Hom (\bar{L}, \bar{L})$. Moreover, $\mathcal{LM}^\BL (\bar{L})$ is a split-generator of ${\mathcal{MF}} (W)$ (Corollary \ref{cor:cptgen1}).
%
%Therefore, it only remains to know 
%
We first show that the Seidel Lagrangian $\bar{L}$ split-generates the Fukaya category of $\mathbb{P}^1_{a,b,c}$ in order to apply Theorem \ref{thm:criterion_equiv} later. We use the generation criterion in \cite{AFOOO}, where  the following corollary was
already known in principle in \cite{A10}.
\begin{theorem}[\cite{AFOOO}]
If $(M, \omega)$ is a compact $2d$-dimensional symplectic manifold, $\mathcal{L}$ a full subcategory of ${\Fuk} (M)$ with some finite set of objects, and the map
$$ \mathcal{OC} : HH_{2d} (\mathcal{L}) \to QH_{2d} (M)$$
hits the unit, then $\mathcal{L}$ split-generates ${\Fuk} (M)$.
\end{theorem}
If $M$ is a surface of genus $g \geq 1$, the condition of the theorem is of a particularly simple form.
\begin{corollary}[\cite{AFOOO}] \label{cor:surfacegene}
If $M$ is a surface with genus $\geq1$, and $\mathcal{L}$ is a full subcategory of $\mathcal{F}uk (M)$ consisting of curves which divide $M$ into polygons, then $\mathcal{L}$ split-generates $\mathcal{F}uk (M)$.
\end{corollary}

Recall that there is a manifold covering $\Sigma  \to \mathbb{P}^1_{a,b,c}$ of $\mathbb{P}^1_{a,b,c}$ whose deck transformation group is $G$, and that we have fixed one of such $\Sigma$ and defined $\mathcal{F}uk (\mathbb{P}^1_{a,b,c})$ as the $G$-invariant part of $\mathcal{F}uk (\Sigma)$.
 Let $\mathcal{L}:=\{L_1, \cdots, L_d \}$ be the family of all Lagrangians in $\Sigma$ each of which lifts $\bar{L}$.  Then $\mathcal{L}$ is preserved by the $G$-action.  In our definition of $\mathcal{F}uk (\mathbb{P}^1_{a,b,c})$, the Seidel Lagrangian $\bar{L}$ corresponds to $\mathcal{L}$ viewed as an object of the $G$-invariant part of $\mathcal{F}uk (\Sigma)$.
 %(Previously, the union $\cup_i L_i$ was regarded as a $G$-equivariant Lagrangian immersion into $\Sigma$.) 

\begin{prop}\label{prop:seigen7}
The Seidel Lagrangian $\bar{L}$ split-generates $\mathcal{F}uk (\mathbb{P}^1_{a,b,c})$.
\end{prop}

\begin{proof}
From the above discussion, we have to prove that $\mathcal{L}$ split-generates the $G$-invariant part of $\mathcal{F}uk (\Sigma)$ 
%(by taking iterated cones over $G$-invariant morphisms). 
The condition  $\frac{1}{a} + \frac{1}{b} + \frac{1}{c} \leq 1$ implies $g \geq 1$, and hence we can apply Corollary \ref{cor:surfacegene} to $M=\Sigma$ to see that the family $\mathcal{L}$ split-generates $\mathcal{F}uk (\Sigma)$. i.e. it generates the derived Fukaya category of $\Sigma$ is given as the idempotent completion of the $H^0$ of the $\AI$-category of twisted complexes over $\mathcal{L}$.

In general, the $G$-orbit of a twisted complex over $\mathcal{L}$ (obtained as the direct sum of $G$-images of a twisted complex over $\mathcal{L}$) is again a twisted complex, and such a sum may be considered as a twisted complex over $G \cdot \mathcal{L}$. Hence
it is not hard to see that if $\mathcal{L}$ split-generates derived Fukaya category of $\Sigma$, then $G \cdot \mathcal{L}$ 
split-generates the $G$-invariant part of the derived Fukaya category of $\Sigma$ (cf. \cite[Corollary 4.7]{Wu}), which we took as a definition of the
derived Fukaya category of the quotient orbifold $\mathbb{P}^1_{a,b,c}=[\Sigma/G]$. More precisely, $G \cdot \mathcal{L}$ generates the $G$-invariant part of $D^\pi \mathcal{F}uk (\Sigma)$ which is tautologically equivalent to the split closure of $\mathcal{F}uk (\mathbb{P}^1_{a,b,c})$ (cf. \cite[Lemma 2.19]{Wu}).
$G \cdot \mathcal{L}$ is simply a direct sum of several copies of $\mathcal{L}$, which implies $\mathcal{L}$ split-generates the $G$-invariant part of $\mathcal{F}uk (\Sigma)$.
%
%
%----------------old-----------------------------------------------
%
%Let $\mathbf{x}$ be a class in $HH_{2} (\mathcal{L})$ whose image under $\mathcal{OC}$ is exactly the unit $1_{\Sigma}$ in $QH_2 (\Sigma)$. Observe that $HH_\ast (\mathcal{L})$ admits a diagonal action of $G$ since $\mathcal{L}$ is closed under the $G$-action. Moreover, it is easy to see from construction that the map $\mathcal{OC} : HH_2 (\mathcal{L}) \to QH_2 (\Sigma)$ is $G$-equivariant. It follows that
%\begin{equation}\label{eq:hitunit1}
%\frac{1}{|G|} \sum_{g \in G} g \cdot \mathbf{x} \stackrel{\mathcal{OC}}\longmapsto \frac{1}{|G|} \sum_{g \in G} g_\ast (1_{\Sigma} ) = 1_{\Sigma},
%\end{equation}
%since $g_\ast (1_{\Sigma} ) = 1_{\Sigma}$.
%Let $\mathbf{y}$ denote $\frac{1}{|G|} \sum_{g \in G} g \cdot \mathbf{x}$. As $\mathbf{y}$ is $G$-invariant, $\mathbf{y}$ can be thought of as an element of $HH_2 ( CF(\bar{L}, \bar{L}))$. Note that the quotient map $\Sigma \to \mathbb{P}^1_{a,b,c}$ sends $1_{\Sigma}$ to $|G| \cdot 1_{\mathbb{P}^1_{a,b,c}}$. Therefore, \eqref{eq:hitunit1} tells us that $\mathbf{y} /|G| \in HH_2 (CF(\bar{L}, \bar{L})$ goes precisely to $1_{\mathbb{P}^1_{a,b,c}}$ via the downstair open-closed map
%$$ \mathcal{OC} : HH^2 ( CF(\bar{L}, \bar{L}) ) \to QH_2 (\mathbb{P}^1_{a,b,c} ),$$
%which implies the proposition.
\end{proof}
%
%\begin{corollary}
%The Seidel Lagrangian $\bar{L}$ spilt generates $\mathcal{F}uk (\mathbb{P}^1_{a,b,c})$. 
%\end{corollary}
%\begin{proof}
%Recall that $\mathcal{F}uk (\mathbb{P}^1_{a,b,c})$ is defined by the orbifolding of the Fukaya category of $\Sigma$. Thereofore, it suffices to show that the family $\{L_0, \cdots, L_k\}$ split-generates the Fukaya category of $\Sigma$. {\color{red} WILL BE WRITTEN.}
%\end{proof}

In summary, both $\bar{L}$ and $P_{\bar{L}}$ are split-generators (Proposition \ref{prop:seigen7}, Corollary \ref{cor:cptgen1}), and our functor induces an isomorphism between (cohomogies of) their endomorphism spaces (Section \ref{subsec:mirrormorlev}). Therefore, Theorem \ref{thm:criterion_equiv} implies that

\begin{theorem}\label{thm:mirsymmstate}
There is an equivalence of triangulated categories
\begin{equation} \label{eq:HMS-abc}
D^{\pi} (\mathcal{F}uk (\mathbb{P}^1_{a,b,c})) \stackrel{\cong}{\longrightarrow}  D^{\pi} ({\mathcal{MF}} (W)).
\end{equation}
\end{theorem}

As discussed in Subsection \ref{sec:dualgponmir}, both $\mathcal{F}uk (\mathbb{P}^1_{a,b,c})$ and ${\mathcal{MF}} (W)$ admit a natural action of the character group of the deck transformation group $G$ associated to the covering $\Sigma \to \mathbb{P}^1_{a,b,c}$. Since liftings of $\bar{L}$ generate the Fukaya category of $\Sigma$, one can deduce the ``upstair homological mirror symmetry" from Theorem \ref{thm:mirsymmstate} if $G$ is abelian. Namely we have the equivalence between $D^\pi (\mathcal{F}uk (\Sigma))$ and $D^\pi ( {\mathcal{MF}}_{\WH{G}} )$. See Proposition \ref{prop:gendualact}.

\begin{remark}
The character group action on $\mathcal{F}uk (\mathbb{P}^1_{a,b,c})$  in the case of $(a,b,c)= (5,5,5)$ agrees with the one appearing in \cite[Section 9]{Se}. 
\end{remark}

Recall that the left hand side of Equation \eqref{eq:HMS-abc} depends on the choice of a compact smooth cover in which the Seidel Lagrangian lifts as an embedded curve.  On the other hand the right hand side is independent of such a choice.  Thus we have

\begin{corollary} \label{cor:indep-cover}
The derived Fukaya category $D^{\pi} (\mathcal{F}uk (\mathbb{P}^1_{a,b,c}))$ is independent of choice of a compact smooth cover of $\bP^1_{a,b,c}$ in which the Seidel Lagrangian lifts as an embedded curve.
\end{corollary}

\section{Fermat hypersurfaces}\label{ndim}

In the previous section we apply our mirror construction to the orbifold $\bP^1_{a,b,c}$ and its manifold cover.  On the other hand the construction of mirror functor in this paper is completely general and can be applied to K\"ahler manifolds of arbitrary dimensions.  In this section, we discuss Fermat hypersurfaces and give two conjectural statements about the mirror map and our functor.  We study our construction for compact toric manifolds in \cite{CHL-toric}, and will study other classes of geometries such as toric Calabi-Yau orbifolds, rigid Calabi-Yau manifolds and higher-genus orbifold surfaces in a series of forthcoming works.

Sheridan \cite{Sh} proved homological mirror symmetry for Fermat-type Calabi-Yau hypersurfaces
$$ \tilde{X} = \left\{[z_0:\ldots:z_{n+1}] \in \proj^{n+1}: z_0^{n+2} + \ldots + z_{n+1}^{n+2} = 0 \right\} $$
using a specific immersed Lagrangian $\bar{L}$ in the quotient $X = \tilde{X} / G$
where
$$G = \frac{(\Z/{(n+2)})^{n+2}}{\{(\lambda,\ldots,\lambda):\lambda \in \Z/{(n+2)}\}} \cong (\Z/{(n+2)})^{n+1}$$
acts on $\tilde{X}$ by $(\lambda_0,\ldots,\lambda_{n+1}) \cdot [z_0:\ldots:z_{n+1}] := [\lambda_0 z_0:\ldots:\lambda_{n+1} z_{n+1}]$.  

$X$ is an orbifold $\C\proj^{n}$ with orbifold strata located at the hyperplanes $\{z_j = 0\}$ for $j=0,\ldots,n$.  For $n=1$ $\bar{L}$ is the Lagrangian introduced by Seidel, and we have used it to construct a mirror functor for the elliptic curve (and its $\Z/3$ quotient) in Section \ref{ex333}.  In general dimensions Sheridan \cite{Sh} showed that $\bar{L}$ lifts to Lagrangian spheres which split-generate the Fukaya category $\Fuk(\tilde{X})$ by using the generating criterion of \cite{AFOOO}.

We refer the detailed construction of the immersed Lagrangian $\bar{L}$ to \cite[Section 6.1]{Sh}.  $\bar{L}$ is a perturbation of the real Lagrangian $\R\proj^{n} \subset \C\proj^{n}$ by using Weinstein Lagrangian neighborhood theorem: the standard double cover $ S^{n} \to \R\proj^{n}$
can be extended to an immersion $ D^*_\eta S^n \to X $
where $D^*_\eta S^n$ is the cotangent disc bundle of radius $\eta$.  Then one considers the graph of an exact one form $\epsilon d f$ in $D^*_\eta S^n$, where $\epsilon \in \R_+$ is small enough, and $f: S^n \to \R$
is a Morse function with $\left(\begin{array}{cc} n+2 \\ k+1\end{array}\right)$
critical points of index $k$ for $k=0,\ldots,n$ with $f(-x) = -f(x)$ for all $x \in S^n$.  Then $\bar{L} \subset X$ is defined to be the image of this graph, and $f$ is chosen such that the immersed points (whose pre-images are the critical points of $f$) are disjoint from the orbifold strata of $X$.

The critical points of $f$ can be labelled by non-empty proper subsets $I \subset \{1,\ldots,n+2\}$, where $|I|-1$ is the index of the critical point $X_I$.  $X_I$ is also regarded as an immersed point of $\bar{L}$, and we have $\bar{X}_I = X_{\bar{I}}$, where $\bar{I}$ denotes the complement of $I$.  Then the deformation space $H$ of $\bar{L}$ is spanned by $X_I$'s, the point class and fundamental class of $S^n$.  Define $V \subset H$ to be the subspace spanned by $X_1,\ldots,X_{n+2}$.  Elements in $V$ are of the form
$ b = \sum_{i=1}^{n+2} x_i X_i. $

We have the anti-symplectic involution
$ \Phi: [z_0: \ldots :z_n] \mapsto [\bar{z}_0: \ldots :\bar{z}_n]$
on $\C\proj^{n}$ under which $\R\proj^{n} \subset \C\proj^{n}$ is the fixed locus.  The immersed Lagrangian $\bar{L}$ is fixed under this involution, and so holomorphic discs bounded by $\bar{L}$ form pairs $u$ and $\Phi \circ u$ (with the conjugate complex structure in the domain so that the map remains to be holomorphic).  In dimension $n=1$, we use this involution to prove that all elements in $V$ are weakly unobstructed.  We also show that the generalized SYZ map equals to the mirror map (see Section \ref{sec:syzmirror}).
We conjecture that these statements still hold in general dimensions:

\begin{conjecture} \label{conj-mm}
\begin{enumerate}
\item \label{conj1} The elements in $V$ defined above are weak bounding cocycles.  Thus we can apply our construction to obtain the mirror superpotential $W$, which is written in terms of (virtual) counts of polygons.

\item \label{conj2} There exists a change of coordinates on $(x_1,\ldots,x_{n+2})$ such that $W$ is equivalent to
$ \sum_{i=1}^{n+2} x_i^{n+2} + \check{q}(q) x_1\ldots x_{n+2} $
where $\check{q}(q)$ is the inverse mirror map for Fermat hypersurfaces $\tilde{X}$.
\end{enumerate}
\end{conjecture}

Statement (\ref{conj2}) of the above conjecture would give an enumerative meaning of mirror maps of Fermat hypersurfaces.  It is motivated by the principle that (generalized) SYZ mirror should be automatically written in flat coordinates.  Gross-Siebert \cite{GS07} made a conjecture of this type in the tropical world, namely the mirror they construct using tropical discs and scattering is automatically written in flat coordinates.  The corresponding statement for toric compact semi-Fano case and toric Calabi-Yau case has been completely proved in \cite{CLLT12,CCLT13}.  We have verified Conjecture \ref{conj-mm} for $n=1$, and we leave the higher dimensional cases for future research.

\begin{remark}
The mirror maps $q(\check{q})$ for Fermat hypersurfaces are well-known in existing literatures.  In Section \ref{sec:syzmirror} the elliptic curve case is explicitly written down.  For the quintic case ($\dim X=3$) which is of the most interest to physicists, the mirror map is
$$ q = - \check{q} \exp \left( \frac{5}{y_0(\check{q})} \sum_{k=1}^\infty \frac{(5k)!}{(k!)^5}\left[\sum_{j=k+1}^{5k}\frac{1}{j}\right] (-1)^k \check{q}^k \right) $$
where
$ y_0 = \sum_{k=0}^\infty \frac{(5k)!}{(k!)^5} (-1)^k x^k$, see for instance \cite{CK-book}.
\end{remark}

\begin{remark}
The leading order terms of $W$ were calculated by Sheridan \cite[Lemma 6.5.4]{Sh} (up to sign in the coefficients), which is
$ W = q_\alpha x_1 \ldots x_{n+2} + q_\beta \sum_{i=1}^{n+2} x_i^{n+2} + o(q_\alpha,q_\beta) $
where $q_\alpha = \conste^{- A_\alpha}$, $q_\beta = \conste ^{- A_\beta}$, $A_\alpha$ and $A_\beta$ are the symplectic areas of the minimal polygons $\alpha$ with corners $X_1, \ldots, X_{n+2}$ and $\beta$ with corners being $n+2$ copies of $X_1$'s respectively.
\end{remark}
%\begin{prop}
%The elements in $V$ are weak bounding cochains.
%\end{prop}
%\begin{proof}
%The stable polygons contributing to $m_0$ have Maslov-index two.

%Since $\bar{L}$ is fixed under the involution, holomorphic discs form pairs $u$ and $\Phi \circ u$ (with the conjugate complex structure in the domain so that the map remains to be holomorphic).  Let $\beta$ be the class containing $u$, and $\bar{\beta}$ be the class containing $\Phi \circ u$.  

%We have the evaluation maps $\ev^\beta_0: \cM(\beta) \to \bar{L}$ and $\ev^{\bar{\beta}}_0: \cM(\beta) \to \bar{L}$, with one of them being orientation preserving and one being orientation reversing.  Thus the contribution of these two discs to the coefficients of $X_L$ in $m_0^b$ cancel, and thus $m_0^b$ is a multiple of $1_L$.
%\end{proof}

Assuming Statement (\ref{conj1}) of Conjecture \ref{conj-mm} that elements in $V$ are weak bounding cochains, by the construction of Section \ref{GMir} we have the mirror functor
$ \Fuk(\tilde{X}) \to \mathcal{MF}_{(\Z /{(n+2)} )^n}(\tilde{W})$
which is obtained by taking $G$-product of the ``downstair mirror functor"
$ \Fuk (X) \to \mathcal{MF} (W).$
Under the downstair mirror functor, $\bar{L}' \in \Fuk (X)$ is transformed to the matrix factorization $d = m_1^{b,0}$ on $CF^*((\bar{L},b),\bar{L}')$.

In Section \ref{sec:cw}, we have shown for $\bP^1$-orbifolds $\bP^1_{a,b,c}$ that there exists a new exterior algebra structure
which enables us to write down the mirror matrix factorization of the Seidel Lagrangian in the contraction-wedge form. This new product structure accommodates the
change of variables in Definition \ref{def:chco}, which is given by counting of specific
$J$-holomorphic polygons.  Such an expression of the mirror matrix factorization implies that it split generates $\mathcal{MF} (W)$ by \cite{Dy}, and hence leads to homological mirror symmetry.

We expect that in higher dimensions there should be a new product structure $\wedge^{new}$ on $CF(\bar{L},\bar{L})$, which makes the matrix factorization $(CF^*((\bar{L},b),\bar{L}),m_1^{b,0})$ into the following ``wedge-contraction" form.

\begin{conjecture} \label{conj-MF}
Under our mirror functor $\LM^{\BL}$, Sheridan's immersed Lagrangian $\bar{L}$ is transformed to the matrix factorization
$\left(\bigwedge\nolimits_{new}^*\langle X_1,\ldots,X_{n+2} \rangle, \sum_{i=1}^{n+2} x_i X_i \wedge^{new} (\cdot) + w_i \, \iota_{X_i}^{new} \right)$
where $\sum_{i=1}^{n+2} x_i w_i = W$.  Hence $\LM^{\BL}$ derives an equivalence of triangulated categories
$ D^{\pi} ({\Fuk} (X)) \cong  D^{\pi} ( \mathcal{MF}(W)).$
\end{conjecture}

\bibliographystyle{amsalpha}
\bibliography{geometry}

\end{document}